\documentclass[10pt]{article}
\usepackage{graphicx,amssymb,amsmath,amsbsy,eucal,amsfonts,mathrsfs,amscd,bm,amsthm}
\usepackage[all]{xy}
\usepackage{tikz}
\usetikzlibrary{arrows,shapes}
\usetikzlibrary{arrows, decorations.markings,fit}
\usetikzlibrary{calc,3d}
\usepackage{pgfplots,pgfplotstable}
\pgfplotsset{compat=newest}

\usepackage{epigraph}

\usepackage{mathtools}

\usepackage{color}
\usepackage{xcolor,colortbl}

\usepackage{authblk}

\numberwithin{equation}{section}

\allowdisplaybreaks[4]

\newtheorem{definition}[equation]{Definition}
\newtheorem{theorem}[equation]{Theorem}
\newtheorem{proposition}[equation]{Proposition}
\newtheorem{example}[equation]{Example}
\newtheorem{lemma}[equation]{Lemma}
\newtheorem{remark}[equation]{Remark}
\newtheorem{corollary}[equation]{Corollary}

\newcommand{\eps}{\varepsilon}
\newcommand{\Ome}{{\Omega}}

\newcommand{\p}{{\partial}}
\newcommand{\Del}{{\Delta}}
\newcommand{\nab}{\nabla}

\newcommand{\bl}{\bigl\langle}
\newcommand{\br}{\bigr\rangle}

\newcommand{\mct}{\mathcal{T}_h}

\newcommand{\avg}[1]{\bigl\{\hspace{-0.045in}\bigl\{#1\bigr\}\hspace{-0.045in}\bigr\}}

\newcommand{\jump}[1]{\left[\hspace{-0.025in}\left[#1\right]\hspace{-0.025in}\right]}

\DeclareMathOperator{\cof}{cof}

\newcommand{\bld}[1]{\boldsymbol{#1}}

\newcommand{\bchi}{\bld{\chi}}

\newcommand{\ff}{{\bld{f}}}

\newcommand{\bv}{\bld{v}}
\newcommand{\bw}{\bld{w}}

\newcommand{\bp}{\bld{p}}
\newcommand{\bn}{\bld{n}}
\newcommand{\bu}{\bld{u}}

\newcommand{\bDi}{\mathcal{D}}

\newcommand{\bq}{\bld{q}}
\newcommand{\blam}{\bld{\lambda}}

\newcommand{\bz}{\boldsymbol{z}}

\newcommand{\bgamma}{{\bld{\gamma}}}

\newcommand{\be}{\bm e}

\newcommand{\calF}{\mathcal{F}}
\newcommand{\calT}{\mathcal{T}}
\newcommand{\bbH}{\mathbb{H}}
\newcommand{\bbR}{\mathbb{R}}

\newcommand{\bbP}{\mathbb{P}}
\newcommand{\calC}{\mathcal{C}}
\newcommand{\calM}{\mathcal{M}}

\newcommand{\polS}{\mathbb{S}}
\newcommand{\calV}{\mathcal{V}}

\newcommand{\bnu}{\bld{\nu}}
\newcommand{\calG}{\mathcal{G}}
\newcommand{\Q}{\mathbb{Q}}
\newcommand{\N}{\mathbb{N}}
\newcommand{\calO}{\mathcal{O}}
\newcommand{\bg}{\bld{g}}
\newcommand{\MAop}[3]{\textup{MA}_{#2}^{\textup{#3}}[#1]}
\DeclareMathOperator{\Span}{span}
\DeclareMathOperator{\Vor}{Vor}
\newcommand{\calI}{\mathcal{I}}
\newcommand{\diff}{\textup{d}}
\newcommand{\calK}{\mathcal{K}}
\newcommand{\polF}{\mathbb{F}}
\DeclareMathOperator{\PRoj}{Proj}
\newcommand{\eremk}{\hfill\rule{1.2ex}{1.2ex}}
\newcounter{numexample}
\newcommand{\calS}{{\mathcal{S}}}

\newcommand{\MA}{Monge-Amp\`{e}re }
\newcommand{\dm}{\Omega}
\newcommand{\bdry}{\partial \Omega}
\newcommand{\gradv}{\nab}
\newcommand{\R}{\mathbb{R}}
\newcommand{\Z}{\mathbb{Z}}
\newcommand{\Nhb}{\Omega_h^b}
\newcommand{\Nh}{\Omega_h}
\newcommand{\Nhi}{\Omega_h^i}
\newcommand{\Th}{\mathcal{T}_h}
\newcommand{\abs}[1]{|{#1}|}
\newcommand{\norm}[1]{\| #1\|}
\newcommand{\dist}{\textrm{dist}}
\newcommand{\vf}[1]{\bld{#1}}

\title{The \MA equation}

\begin{document}

%

\author[1]{Michael Neilan\thanks{neilan@pitt.edu}}

\author[2]{Abner J. Salgado\thanks{asalgad1@utk.edu}}

\author[3]{Wujun Zhang\thanks{wujun@math.rutgers.edu}}

\affil[1]{Department of Mathematics, University of Pittsburgh}
\affil[2]{Department of Mathematics, The University of Tennessee}
\affil[2]{Department of Mathematics, Rutgers University}

\date{\today}

\maketitle

\maketitle

\begin{abstract}{\normalsize We review recent advances in the numerical analysis 
of the Monge-Amp\`ere equation.  Various computational
techniques are discussed including wide-stencil finite difference schemes, two-scaled methods,
finite element methods, and methods based on geometric considerations.
Particular focus is the development of appropriate stability and consistency
estimates which lead to rates of convergence of the discrete approximations.
Finally we present numerical experiments which highlight each method
for a variety of test problem with different levels of regularity.}\end{abstract}


\section{Introduction}

\epigraph{\tiny{\it All exact science is dominated by the idea of approximation. When a man tells you that he knows the exact truth about anything, you are safe in inferring that he is an inexact man.}}{\tiny{B.~Russell \cite{RussellBook}}}

In this chapter we review recent
progress in the numerical treatment
of  \MA type equations.
In its simplest form, and assuming Dirichlet boundary
conditions, the problem we consider 
is to seek a scalar function $u$ satisfying the {partial differential equation} (PDE)
\begin{subequations}
\label{eqn:MA}
\begin{align}
\label{eqn:MA1}
  \det D^2 u(x) = &\; f(x) \quad {x \in } \Omega
  \\
  \label{eqn:MA2}
  u(x) = &\; g(x) \quad {x \in } \partial \Omega
\end{align}
\end{subequations}
Here, $D^2 u$ denotes the Hessian matrix
of $u$,  $f\ge 0$
and $g$ are given functions, and $\Omega\subset \bbR^d$
is a bounded, convex domain.
Problem \eqref{eqn:MA} is a 
prototypical second--order, fully
nonlinear PDE, and it arises in several broad applications
in differential geometry, meteorology, cosmology, economics, and optimal
mass transportation problems. {Some of these applications are briefly described below}.

Despite its growing list of applications,
and in contrast to its extensive and mature PDE theory,
the construction and analysis of computational methods
for \eqref{eqn:MA} is still a relatively new and emerging 
field in numerical analysis.
Numerical algorithms, based on geometric considerations,
 for the two--dimensional problem ($d=2$)
first appeared in 1988 in \cite{Oliker88},
and the extension to practical three-dimensional
schemes were not introduced until some 20 years
later \cite{FroeseOberman11,FengNeilan09Z,MR2745457,BrennerNeilan12}. {Other early attempts that deserve mention are the least squares and augmented Lagrangian approaches of \cite{MR1989280,MR2111728,MR2169156,GlowinskiDean06B,GlowinskiDean06A}, and we refer the reader to \cite{FGN13} for more details on these schemes}.

The reasons for this delayed development in numerical methods
are plentiful.  The most evident obstacle
is the full non-linearity of the problem.
However, this is arguably a secondary 
difficulty, as black-box nonlinear  solvers
can, at least heuristically, be applied to 
algebraic systems resulting from  discretizations
 of \eqref{eqn:MA}.  {A rather} fundamental
 difficulty to construct, and especially to analyze,
 computational methods for \MA type equations 
 is the variety of solution
 concepts and, correspondingly, the low regularity
 solutions generically possess.  As we explain
 below, weak solutions are not based on variational principles,
 but rather on either geometric considerations
 or by monotonicity conditions of test functions that touch the graph of the solution from above or below. 
 These solution concepts are difficult to mimic at the discrete level,
 and as a result, the construction of convergent
 schemes is an arduous task.
{Finally, as if these complications were not enough, the \MA equation \eqref{eqn:MA} is usually supplemented by the constraint that the solution $u$ is convex. This is not only because of geometric applications, but in many cases a necessary condition for uniqueness, and for the existence of a well--developed PDE theory. As convexity is a global constraint, it is very difficult to enforce it in a discrete setting.}
  
Nonetheless, an explosion of results and new techniques to develop them
in computational methods for \eqref{eqn:MA}
have occurred  during the last $10$ years.
These include 
the construction of monotone,
wide-stencil finite difference schemes, semi-Lagrangian
methods, and finite element methods.
Within only the past few years, significant progress
has been made in the convergence analysis
with an emphasis on the rates of convergence for various
discretization schemes.

The main goal of this chapter is to highlight these
recent advances in the numerical analysis
of the Monge-Amp\`ere problem \eqref{eqn:MA}.
To this end, we organize the paper as follows.
After stating some geometric applications
and a brief PDE theory of the Monge-Amp\`ere
problem in this section, 
we discuss wide stencil finite difference schemes
in Section \ref{sec:FD}.  There we 
introduce the monotone finite difference
schemes \cite{Oberman08,FroeseOberman11,MR2745457} and the corresponding
filtered schemes \cite{MR3033017},  lattice reduction schemes \cite{BCM16}, 
methods based on power diagrams \cite{Mirebeau15}, and the so-called
two scale methods \cite{NochettoNtogkasZhangfirsttwoscale,NochettoNtogkasZhang,NochettoNtogkasFilter}.  Of particular focus will be the rates of convergence
of these schemes if available.
Next, in Section \ref{sec:OP},
we review the {original method of \cite{Oliker88}, which in honor of its proponents henceforth we shall call the
Oliker-Prussner scheme. This method} is based on geometric interpretations
of the \MA operator and the notion of Alexandrov 
solutions.  Again, the emphasis of the discussion
is on consistency error and pointwise rates 
of convergence recently established
in \cite{NochettoZhangMA}.
Section \ref{sec:FEM} discusses
finite element methods for both
smooth and singular solutions.
Finally in Section \ref{sec:Numerics}
we perform some numerical experiments
using some of the methods we discuss
in this review for a variety of test problems
with different levels of regularity.

\subsection{Geometric applications}
\label{sub:Geoapplications}
To draw connections with the theme
of the {current volume in the Handbook of Numerical Analysis}, and to further emphasize 
the prevalence of the \MA problem,
in this section we briefly summarize 
some applications with a geometric
flavor where the \MA problem plays
an essential role.

\subsubsection{Gauss curvature problem}\label{subsub:Gauss}

The classic Gauss curvature problem (cf.~\cite{MR1305147,Oliker84,GuanSpruck93}) seeks a manifold $\mathcal{M}\subset \mathbb{R}^{n+1}$ 
with prescribed boundary and Gauss curvature $\mathcal{K}$.  We recall that Gauss curvature
is the product of the principal curvatures, which themselves are the eigenvalues
of the shape operator (or Weingarten map).  One may reduce this problem to a PDE problem of Monge-Amp\`ere type
if one assumes that the manifold
is the graph of a function, i.e.,
\[
\mathcal{M} = \{(x,u(x)):\ u:\Omega\to \mathbb{R}\}.
\]
The shape operator is given by $s =  \mathrm{I}^{-1}\ \mathrm{I\!I}$,
where $\mathrm{I}$ and $\mathrm{I\!I}$ denote, respectively, the first and second fundamental forms.
In the case that $\mathcal{M}$ is the graph of the function $u$, we have
${\mathrm{I}} =I +\nab u\otimes \nab u$ 
and ${\mathrm{I\!I}} = \frac{D^2 u}{\sqrt{1+|\nab u|^2}}$, {where $I$ denotes the $d\times d$ identity matrix}. 
Therefore the Gauss curvature
is given by 
\[
\mathcal{K} = \det(s) = \frac{\det(\mathrm{I\!I})}{\det(\mathrm{I})}= \frac{\det D^2 u}{(1+|\nab u|^2)^{(d+2)/2}}.
\]

Thus, the problem is  to find a scalar function $u: \bar\Omega \to \R$ satisfying
\begin{subequations}
\label{eq:Gausscurvature}
  \begin{alignat}{2}
  \label{eq:Gausscurvaturein}
    \det D^2 u(x) &= \mathcal{K}(x)(1+|\nab u(x)|^2)^{(d+2)/2} \quad &&\mbox{in } \Omega, \\
  \label{eq:Gausscurbaturebc}
    u(x) &= g(x) \quad &&\mbox{on } \p\Omega.
  \end{alignat}
\end{subequations}
In conclusion the Gauss curvature problem, in this setting,
seeks solutions of a \MA type problem with lower--order terms.

\subsubsection{Reflector design problem}\label{subsub:Reflector}
The reflector design problem \cite{NorrisWestcoff76,OlikerWaltman86,XWang96,Oliker87Reflect} 
{can be posed as follows}: Let $S^2$ be the unit
sphere in $\bbR^3$ centered at the origin, 
and let $\Omega,\calO$ be two disjoint
domains on $S^2$.  Let $f$ be a positive
function defined on $\calO$, and suppose
that rays originate from the origin with
density $\rho$.  We then seek a surface, called $\Gamma$,
whose radial projection onto $S^2$ equals $\Omega$,
such that the directions of the reflected rays cover $\calO$
with distributed density equal to $f$.  

To formulate a PDE {model for} this problem,
we set $\Gamma = \{x m(x):\ x\in \Omega\}$,
so that if a ray radiates from the origin with
direction $x$, then it is reflected at the point $xm(x)$.
This will create a reflected ray in the direction $T(x)\in \calO$.
Now if we denote by $\bn$ the unit normal of $\Gamma$ at $x m(x)$,
then we have $T(x)-x = -2 (x\cdot \bn)\bn$, and calculations
 show that $\bn = (\nab m-m x)/\sqrt{{m}^2 + |\nab m|^2}$.
Here, $\nab = e^{ij} \p_i x\p_j$,
where $x$ is a smooth parametrization of $S^2$, $e = e_{ij}\diff t^i \diff t^j$ is the first
fundamental  form of $S^2$,  $e^{ij} = (e_{ij})^{-1}$, and $\p_j = \p/\p t^j$.
 Combining these two identities we find that
the direction $T$ is related to $m$ via
\begin{align}\label{eqn:TuRelation}
T(x) = \frac{2 m \nab m + (|\nab m|^2-m^2) x}{m^2 +|\nab m|^2}.
\end{align}
Next, if the directions of the reflected light
do not overlap and if no loss of energy
occurs in the reflection, then we have 
the energy conservation property
\[
\int_E \rho(x)\diff x= \int_{T(E)} f(y)\diff y = \int_E f(T(x)) \frac{|\p_1 T(x)\times \p_2 T(x)|}{\det(e_{ij})}\diff x
\]
for all Borel sets $E\subset \Omega$.
Thus we have, at least formally,
\[
\frac{|\p_1 T(x)\times \p_2 T(x)|}{\det(e_{ij})} = \frac{\rho(x)}{f(T(x))}.
\]
Finally, we set $u(x) = 1/m(x)$, and substitute
\eqref{eqn:TuRelation} into this last equation to get
the {following problem of \MA type}
(see \cite{Oliker93,XWang96} for details)
\begin{align*}
\frac{\det(D^2 u+(u-\eta)e_{ij})}{\eta^2 \det(e_{ij})} = \frac{\rho(x)}{f(T(x))}\qquad x\in \Omega,
\end{align*}
where $T$ is given by \eqref{eqn:TuRelation}
and $\eta = (|\nab u|^2 + u^2)/(2u)$.
%

\subsubsection{Affine plateau problem}
Following \cite{TrudWang05,TrudWang08,Calabi90}, we consider
the following problem.
Let $\calM_0 \subset \bbR^{d+1}$ be a bounded
and connected hypersurface with smooth boundary
that is locally uniformly convex
We denote by $S[\calM_0]$ the set of locally
uniformly convex hypersurfaces that can be smoothly
deformed from $\calM_0$ within the family of locally
uniformly convex hypersurfaces and whose Gauss map
images lie in that of $\calM_0$.
As in Section~\ref{subsub:Gauss}, for a manifold $\calM$
we denote by $\mathrm{I\!I}$ its second fundamental form and
by $\calK$ its Gauss curvature.  Associated with $\calM$
is the {\it Berwald--Blaschke metric}
\[
g = \calK^{-1/(d+2)} \mathrm{I\!I},
\]
which is an affine invariant Riemmannian metric on the surface.
The {\it affine Plateau problem} is then to determine 
the maximizer of the  {\it affine area functional}
\[
A(\calM) = \int_{\calM} \calK^{1/(d+2)} \diff \calM .
\]
over $S[\calM_0]$. 

Recall that if $\calM = \calM_u$ is the graph of a function $u:\Omega\to \bbR$,
with $\Omega\subset \bbR^n$, then the Gauss curvature
is $\calK = \det(D^2 u)/(1+|\nab u|^2)^{(d+2)/2}$, and so, we have by a change of variables,
\[
A(\calM_u) = \int_\Omega (\det D^2 u(x))^{1/(d+2)} \diff x.
\]
Thus if $\calM_0$ is the graph of a locally uniformly
convex $g$, then in the graph case, $S[\calM_0]$
consists of the graphs of locally uniformly convex functions $v\in C^2(\Omega)\cap C^0(\bar\Omega)$
satisfying $v = g$ on $\p\Omega$ and $\nab v(\Omega)\subset \nab g(\Omega)$.
In this setting the affine Plateau problem
seeks $u$ such that
\[
A(\calM_u) = \sup\{A(\calM_v):\ \calM_v\in S[\calM_0]\}.
\]
Formally taking the Euler--Lagrange equation yields
the {\it affine maximal surface equation}
\begin{align*}
 \cof D^2 u:D^2 w = 0,\qquad w = \big(\det D^2 u\big)^{-(d+1)/(d+2)}.
 \end{align*}

\subsubsection{Optimal mass transport problem}\label{subsub:OT}

\epigraph{\tiny{\it This problem appeared as a generalization of an earlier considered practical problem of assigning production locations on a railway network to consumption locations with minimum total transportation expenses.}}{\tiny{L.V. Kantorovich \cite{MR2117877}}}

The optimal mass transport problem was originally proposed
by G.~Monge in the $18$th century to find the optimal way
to move oil to an excavation with minimal transportation cost.
In general, the mass transport problem seeks, 
for two given sets and densities, the optimal mass--preserving mapping
between them.

In further detail, given bounded $\Omega,\calO\subset \bbR^d$
and 
 measures $\rho_\Omega:\Omega \to \bbR$,
$\rho_{\calO}:\calO\to \bbR$, the optimal transport
problem with quadratic cost
seeks a map $T:\Omega\to \calO$
such that $T_{\#} \rho_\Omega = \rho_{\calO}$ that
minimizes the functional
\begin{align}\label{eqn:OTFunc}
\frac12 \int_\Omega |x-T(x)|^2\diff \rho_\Omega(x)
\end{align}
over all mass preserving maps.  
Here, we assume that
 the measures
are absolutely continuous with respect to Lebesgue measure,
with $\diff \rho_\Omega = f_\Omega \diff x$ and $\diff \rho_{\calO} = f_{\calO} \diff x$,
and that the measures satisfy the mass balance condition
\[
\int_\Omega f_\Omega(x) \diff x  = \int_{\calO} f_{\calO}(x) \diff x.
\]
Above, we denoted by
$T_{\#} \rho_{\Omega}$ the pushforward of the measure
$\rho_{\Omega}$ under the mapping $T$, i.e., under the given
assumptions, we have
\[
\int_E f_{\calO}(x)\diff x = \int_{T^{-1}(E)} f_{\Omega}(x)\diff x.
\]
Thus, by a change of variables, we have,
at least formally, 
\begin{align}
\label{eqn:measureP}
\det(\nab T(x))f_{\calO}(T(x)) = f_\Omega(x)\qquad x\in \Omega,
\end{align}
with $T(\Omega) \subset \calO$.
Thus in summary, we seek a mapping $T$ that minimizes
\eqref{eqn:OTFunc} with the constraint \eqref{eqn:measureP}.
One of the fundamental {results in the theory of optimal transport} \cite{MR1100809,MR1009457,MR1035606,MR1062553}
is that there exists a unique solution to this problem
and that this optimal mapping is characterized as the gradient
of some convex function $u$:
\[
T(x) = \nab u(x).
\]
Hence, by substituting this relation into \eqref{eqn:measureP}, we
see that  the problem reduces to a \MA type PDE
\begin{align}
\label{eqn:OPeqn}
f_{\calO}(\nab u(x)) \det D^2u(x)= f_\Omega(x)\qquad x\in \Omega.
\end{align}
with the constraint $\nab u(\bar\Omega) \subset \bar\calO$.
Thus we find that, with quadratic cost, the optimal mass transport
problem reduces to a \MA equation
with transport boundary boundary conditions.

\subsection{Solution concepts for the \MA equation}
\label{sub:PDE}

\epigraph{\tiny{\it It is impossible to understand an unmotivated definition [...]}}{\tiny{V.I.~Arnold \cite{ArnoldBadAss}}}

In order to properly analyze the numerical schemes that we present below, it is important to understand in which sense a function $u:\bar\Omega \to \R$ must satisfy the equation and boundary conditions in \eqref{eqn:MA} to be a solution. It is not our intention here to give a survey of the PDE theory regarding the \MA  equation, and we refer the reader to \cite{Gutierrez01,MR3617963,MR1305147} for an in-depth presentation.

\subsubsection{Classical solutions}
\label{subsub:classicalsol}

The first definition of a solution to \eqref{eqn:MA} is that of a classical solution. Essentially we require that \eqref{eqn:MA} holds at every point of $\bar\Omega$.

\begin{definition}[classical solution]
\label{def:classicalsol}
A function $u \in C^2(\Omega) \cap C(\bar\Omega)$ is called a classical solution of \eqref{eqn:MA} if these identities hold for every $x \in \bar\Omega$.
\end{definition}

Notice that this necessarily implies that the right hand side $f:\Omega \to \R$ is continuous. Regarding the existence of classical solutions we have the following result; see \cite[Section 3.1]{MR3617963} for a detailed presentation.

\begin{theorem}[existence of classical solutions]
\label{thm:existenceclassical}
Let $\alpha \in (0,1)$. Assume that $\Omega$ is a bounded and uniformly convex domain, whose boundary is of class $C^3$, $f \in C^\alpha(\bar\Omega)$ with $f \geq f_0 >0$, and $g \in C^3(\p \Omega)$. Then problem \eqref{eqn:MA} has a unique solution $u \in C^{2,\alpha}(\bar\Omega)$.
\end{theorem}

It is important to notice that classical solutions may not always exist, see for instance the counterexample given in \cite[Section 3.2]{MR3617963}. This motivates us to introduce weaker notions of solutions.

\subsubsection{Viscosity solutions}
\label{subsub:viscosols}

The \MA operator $w \mapsto \det D^2w$ is a fully  nonlinear second order operator, that is it depends nonlinearly on the highest (in this case second) order derivatives that appear in the expression. For this reason, the theory regarding fully nonlinear operators can guide us to develop a notion of solution (viscosity solution) that is weaker than classical. We refer the reader to \cite[Chapter 17]{GT}, \cite{MR1351007}, \cite{CIL} and \cite[Section 2]{NSWActa} for additional details.

We begin with a definition that encodes the type of admissible nonlinearities that will allow for the development of the theory of viscosity solutions. Here and in what follows we denote by $\polS^d$ the collection of symmetric $d \times d$ matrices. {The set $\polS^d$ is endowed with a partial order: if $M,N \in \polS^d$ then we say that $M \leq N$ if $\bv \cdot M \bv \leq \bv \cdot N \bv$ for every $\bv \in \R^d$.}

\begin{definition}[elliptic operator]
\label{def:FLelliptic}
Let $F : \bar \Omega \times \bbR \times \polS^d \to \bbR$ be locally bounded. We say that $F$ is {\it elliptic} if it satisfies the following {\it monotonicity condition}: For all $x \in \bar\Omega$, $r,s \in \bbR$ and $M,N \in \polS^d$ with $r \geq s$ and $M \leq N$ then
\[
  F(x,r,M) \leq F(x,s,N).
\]
Moreover, we say $F$ is {\it uniformly elliptic} if for all $r,s \in \bbR$ and $M \in \polS^d$ with $r \geq s$ we have
\[
  F(x,r,M)\le F(x,s,M),
\]
and, in addition, there are constants $0<\lambda \leq \Lambda$ such that for all $M \in \polS^d$ we have
\[
  \lambda \|N\|_2 \leq F(x,r,M+N) - F(x,s,M) \leq \Lambda \|N\|_2, \quad \forall N \geq 0.
\]
\end{definition}

Letting $F : \bar \Omega \times \bbR \times \polS^d \to \bbR$ be an elliptic operator as defined above, we consider the fully nonlinear elliptic problem
\begin{equation}
\label{eq:BVPBarlesSouganidisStyle}
  F(x,u(x),D^2u(x)) = 0 \quad \mbox{in } \bar \Omega.
\end{equation}
To be able to properly describe the notion of viscosity solutions we need to recall the following.

\begin{definition}[upper and lower semicontinuous envelopes]
\label{def:semicontenvelopes}
Let $w: \bar \Omega \to \bbR$. By $w^\star \in USC(\bar\Omega)$ and $w_\star \in LSC(\bar\Omega)$, we denote the upper and lower semicontinuous envelopes, respectively, of the function $w$. In other words
\[
  w^\star(x) = \limsup_{y \to x} w(x), \qquad w_\star(x) = \liminf_{y \to x} w(x).
\]
Finally, by $USC(\bar\Omega)$ and $LSC(\bar\Omega)$ we denote, respectively, the sets of upper and lower semicontinuous functions.
\end{definition}

We are now ready to introduce the notion of viscosity solution.

\begin{definition}[viscosity solution]
\label{def:viscosol}
Let $F$ be elliptic in the sense of Definition~\ref{def:FLelliptic}. We say that the locally bounded function $u :\bar \Omega \to \bbR$ is:
\begin{enumerate}
  \item A {\it viscosity subsolution} of \eqref{eq:BVPBarlesSouganidisStyle} if whenever $x_0 \in \bar\Omega$, $\varphi \in C^2(\bar\Omega)$ and $u^\star-\varphi$ has a local maximum at $x_0$ we have that
  \[
    F_\star(x_0,\varphi(x_0),D^2 \varphi(x_0)) \geq 0.
  \]
  
  \item A {\it viscosity supersolution} of \eqref{eq:BVPBarlesSouganidisStyle} if whenever $x_0 \in \bar\Omega$, $\varphi \in C^2(\bar\Omega)$ and $u_\star-\varphi$ has a local minimum at $x_0$ we have that
  \[
    F^\star(x_0,\varphi(x_0),D^2 \varphi(x_0)) \leq 0.
  \]
  
  \item A {\it viscosity solution} if it is a sub- and supersolution.
\end{enumerate}
\end{definition}

The condition ``$u^\star - \varphi$ has a local maximum at $x_0$'' is usually phrased as ``{\it $\varphi$ touches the graph of $u$ from above at $x_0$}''. The reader is encouraged to draw a picture to see why these two have the same meaning.
Similarly, ``$u_\star - \varphi$ has a local minimum at $x_0$'' is: ``{\it $\varphi$ touches the graph of $u$ from below at $x_0$}''.

One of the main technical tools in asserting existence and uniqueness of viscosity solutions is a comparison principle.

\begin{definition}[comparison principle]
\label{def:comparisonVisco}
We say that problem \eqref{eq:BVPBarlesSouganidisStyle} satisfies a comparison principle if whenever $\overline{w} \in USC(\bar\Omega)$ and $\underline{w} \in LSC(\bar\Omega)$ are sub and supersolutions, respectively, then we must have
\[
  \overline{w} \leq \underline{w}.
\]
\end{definition}

Notice now that if we define
\begin{equation}
\label{eq:MaasFNLelliptic}
  F_{MA}(x,r,M) = \begin{dcases}
                      \det M - f(x), & x \in \Omega, \\
                      g(x) - r, & x \in \p\Omega,
                    \end{dcases}
\end{equation}
this operator satisfies the monotonicity conditions given in Definition~\ref{def:FLelliptic} {\it only} if we restrict the third argument to the set of positive semidefinite matrices 
which we denote by $\polS^d_+$. Consequently, we need to restrict the class of admissible functions, that define a viscosity solution to \eqref{eqn:MA} to the set of convex functions.

\begin{definition}[viscosity solution]
\label{def:viscosolMA}
Let $u \in C(\bar\Omega)$ be a convex function. We say that $u$ is:
\begin{enumerate}
  \item A {\it viscosity subsolution} of \eqref{eqn:MA} on the set of convex functions if $u \leq g$ on $\partial\Omega$ and, whenever $x_0 \in \Omega$, $\varphi \in C^2(\Omega)$, and $u-\varphi$ has a local maximum at $x_0$ we have that
  \[
    \det D^2 \varphi(x_0) \geq f(x_0).
  \]
  
  \item A {\it viscosity supersolution} of \eqref{eqn:MA} on the set of convex functions if $u \geq g$ on $\p\Omega$ and, whenever $x_0 \in \Omega$, $\varphi \in C^2(\Omega)$ {is convex}, and $u-\varphi$ has a local minimum at $x_0$ we have that
  \[
    \det D^2 \varphi(x_0) \leq f(x_0).
  \]
  
  \item A {\it viscosity solution} if it is a sub- and supersolution on the set of convex functions.
\end{enumerate}
\end{definition}

The reader may wonder why these definitions are asymmetric. The concept of supersolution requires convexity of the test functions, whereas subsolutions do not. This is due to the fact that, as noted in \cite[Remark 1.3.2]{Gutierrez01}, if $u$ is convex and $u-\varphi$ has a local maximum at $x_0$, then $\varphi$ is  (locally) convex.

The existence and uniqueness of viscosity solutions will be a consequence of Theorems~\ref{thm:existenceuniquenessAlexandrov} and \ref{thm:Aleximpliesvisco} below. Here we mention a remarkable property of viscosity solutions, namely their stability. The following result can be found, for instance, in \cite[Lemma 5.3]{NochettoNtogkasZhangfirsttwoscale}.

\begin{proposition}[continuous dependence]
\label{prop:contdependencevisco}
Let $f_1,f_2 \in C(\bar\Omega)$ with $f_1,f_2\geq 0$ and $g_1,g_2 \in C(\p\Omega)$ and denote by $u_1,u_2 \in C(\bar\Omega)$ the corresponding convex viscosity solutions to \eqref{eqn:MA}. Then we have
\[
  \| u_1 - u_2 \|_{L^\infty(\Omega)} \leq C \| f_1 - f_2 \|_{L^\infty(\Omega)}^{1/d} + \| g_1 - g_2 \|_{L^\infty(\p\Omega)}.
\]
In addition, if $f_1 \geq f_2 \geq 0$ and $g_1 \leq g_2$ we have that $u_1 \leq u_2$.
\end{proposition}

Finally we comment that viscosity solutions can be approximated by classical ones over larger, but smooth, domains; see \cite[Lemma 5.4]{NochettoNtogkasZhangfirsttwoscale}.

\begin{proposition}[smooth approximation]
\label{prop:smoothapproxvisco}
Let $\Omega$ be uniformly convex, $f,g \in C(\bar\Omega)$ with $f\geq 0$, and $u$ the convex viscosity solution to \eqref{eqn:MA}. There exists:
\begin{enumerate}
  \item A decreasing (in the sense of inclusion) sequence of uniformly convex smooth domains $\Omega_n$ such that
  \[
    \dist_H(\Omega_n,\Omega) \to 0, \quad n \to \infty,
  \]
  where by $\dist_H(A,B)$ we mean the $d$--dimensional Hausdorff distance between the sets $A$ and $B$.
  
  \item A decreasing sequence of smooth functions $f_n : \bar\Omega_n \to \R$ with $f_n >0$ such that
  \[
    \| f_n - f \|_{L^\infty(\Omega)} \to 0, \quad n \to \infty.
  \]
  
  \item A sequence of smooth functions $g_n : \bar\Omega_n \to \R$ such that
  \[
    \| g_n - g \|_{L^\infty(\Omega)} \to 0, \quad n \to \infty.
  \]
\end{enumerate}
Moreover, if $u_n \in C(\bar\Omega_n)$ denotes the convex viscosity solution to \eqref{eqn:MA} over the domain $\Omega_n$ and with data $f_n$ and $g_n$, then
\[
  \| u_n - u \|_{L^\infty(\Omega)} \to 0, \quad n \to \infty.
\]
\end{proposition}

\subsubsection{Alexandrov solutions}
\label{subsec:Alexandroff}

Besides the concept of solution in the viscosity sense, another type of weak solution to the \MA equation is the Alexandrov solution, which is based on a geometric interpretation. To motivate it, let $w \in C^2(\dm)$ be convex so that the gradient map $\gradv w: \dm \to \R^d$ is well defined and monotone. In this case, an interesting observation is that $\det D^2 w$ is actually the determinant of the Jacobian of the gradient map. Therefore, for any open (or Borel) subset $E \subset \dm$, we have
\[
  \int_E \det D^2 w(x) \diff x = \int_{\gradv w (E)} \diff y = |\gradv w (E)|
\]
where $|\cdot|$ denotes the $d$-dimensional Lebesgue measure.

What is remarkable about this simple observation is that to make sense of $\det D^2 u$, we only require $\gradv w (E)$ to be well defined for any Borel set $E$. This enables us to make sense of the previous identity even if $w \notin C^2(\dm)$. To define the weak (Alexandrov) solution, we first introduce the subdifferential of a convex function.

\begin{definition}[subdifferential]
Let $\dm$ be convex and $w:\dm \to \R$ be a convex function. The subdifferential of $w$ at point $x \in \dm$ is the set 
\[
  \partial w(x) := 
  \{ \bp \in \R^d, w(x) + \bp \cdot (y - x) \leq w(y) \quad \forall y \in \dm \}.
\]
For any Borel set $E \subset \dm$, we define
\[
  \partial w(E) = \cup_{x \in E} \partial w(x).
\]
\end{definition}

In other words, the subdifferential is the collection of slopes of all affine functions that touch the graph of $w$ at $(x,w(x))$ and bound the graph from below. From this observation, it is easy to see that if $w$ is strictly convex and smooth, then $\p w(x) = \{\nab w(x)\}$. Here we give an example of subdifferential of a convex (but not strictly convex) function.

\begin{example}[subdifferential]
\label{example1}
Let $\dm = B_1 (0) \subset \mathbb R^2$ and
\[
  w(x) = |x|. 
\]
Then at the origin $x=0$, we note that 
\[
  w(0) + \bp \cdot y \leq w(y) \quad \forall y \in \dm
\]
provided that the norm of the vector $|\bp| \leq 1$. Hence, by definition, the subdifferential of $w$ at $x=0$ is the {closed} unit ball centered at $0$, i.e.
\[
  \partial w(0) = {\overline{B_1(0)}}.
\]
At any other point $x \in \dm$, since the function $w$ is differentiable, we note that the inequality 
\[
  w(x) + \bp \cdot (y-x) \leq w(y) \quad \forall y \in \dm
\]
holds if and only if $\bp = \gradv w(x)$. Hence, for all $x \in \dm\setminus\{0\}$, 
\[
  \partial w(x) = \{ \gradv w(x) \}.
\]
\eremk\end{example}

{With this motivation at hand we can introduce the so--called \MA measure, which will be essential in defining Alexandrov solutions.

\begin{definition}[\MA measure]
\label{def:MAmeasure}
Let $\Omega \subset \R^d$ be convex and $w : \Omega \to \R$ be a convex function. The {\it \MA measure} associated to $w$ is
\[
  \mu_w(E) = \left| \p w(E) \right|.
\]
\end{definition}

It can be shown, see \cite[Theorem 2.3]{MR3617963} that this is indeed a locally finite Borel measure on $\Omega$. With this, we are ready to define Alexandrov solutions.
 
\begin{definition}[Alexandrov solution]
\label{def:AlexSolution}
Let $f$ be a Borel measure defined in $\Omega$. A convex function $u\in C(\bar{\Omega})$ is an {\it Alexandrov solution} to the \MA equation \eqref{eqn:MA} if $u = g$ on $\p\Omega$ and $\mu_u = f$, that is,
\begin{align}
\label{Alek}
  |\p u(E)| =  f(E).
\end{align}
for all Borel sets $E\subset \Omega$.
\end{definition}

To illustrate the definition of the Alexandrov solution, we consider Example~\ref{example1}. Let $E \subset \Omega$ be Borel, if the set contains the origin, we have the subdifferential
\[
  \partial u(E) = \cup_{x \in E} \partial u(x) = \overline{B_1(0)},
\]
which yields
\[
  |\partial u(E) | = |B_1(0)| = \pi \quad \mbox{if $x \in E$.}
\]
On the other hand, if the set does not contain the origin, then the subdifferential 
\[
  \p u(E) = \cup_{x \in E} \{\gradv u(x)\} \subset \p B_1(0)
\]
Hence, we get 
$
|\partial u(E) | = 0 
$
if $0 \notin E$. 
Finally, we conclude that $u$ is an Alexandrov solution of Monge-Amp\`{e}re equation 
\[
\det D^2 u (x) = \pi \delta_{\{x=0\}},
\]
where $\delta_{\{x=0\}}$ is the Dirac measure at origin. It is worth mentioning that  $u$ is not a viscosity solution because the right hand side is not a (continuous) function. Also note that the continuity of the source term $f$ is no longer required for \eqref{Alek} to be well defined.

The existence and uniqueness of Alexandrov solutions is summarized in the next theorem, see \cite[Theorem 1.6.2]{Gutierrez01} and \cite[Theorem 2.14]{MR3617963}.

\begin{theorem}[existence and uniqueness]
\label{thm:existenceuniquenessAlexandrov}
Let $\Omega \subset \R^d$ be a strictly convex domain, let $g \in C(\p\Omega)$ and $f$ be a nonnegative Borel measure on $\Omega$ with $f(\Omega)<\infty$. Then there exists a unique convex function $u \in C(\bar\Omega)$ that is a solution of \eqref{eqn:MA} in the sense of Definition~\ref{def:AlexSolution}.
\end{theorem}

An important property of Alexandrov solutions is their stability with respect to weak convergence. We refer the reader to \cite[Lemma 1.2.3]{Gutierrez01} for a proof of the following result.

\begin{lemma}[weak convergence]
\label{weakconvergence}
Let $\{w_k\}_{k=1}^\infty, w$ be convex functions on $\Omega$ and assume that, as $k \to \infty$, we have $w_k \to w$ uniformly over compact subsets of $\dm$. Then,  the associated \MA measures $\mu_{w_k}$ tend to $\mu_w$ weakly, that is, 
\[
  \int_{\dm} \phi(x) \diff \mu_{w_k}(x) \rightarrow \int_{\dm} \phi(x) \diff \mu_w(x),
\]
for every $\phi$ continuous with compact support in $\dm$. 
\end{lemma}

The relation between viscosity and Alexandrov solutions is given in the following result \cite[Propositions 1.3.4 and 1.7.1]{Gutierrez01}. Notice that this result not only shows, as we have already pointed out, that the notion of Alexandrov solution is strictly weaker than that of viscosity solutions but, on the basis of Theorem~\ref{thm:existenceuniquenessAlexandrov}, shows existence and uniqueness of viscosity solutions.

\begin{theorem}[equivalence]
\label{thm:Aleximpliesvisco}
Let $u \in C(\bar\Omega)$ be an Alexandrov solution of \eqref{eqn:MA}. If $f \in C(\Omega)$, then $u$ is also a viscosity solution in the sense of Definition~\ref{def:viscosolMA}. Conversely, if $u$ is a viscosity solution of \eqref{eqn:MA} and $f \in C(\bar\Omega)$ with $f>0$, then $u$ is an Alexandrov solution.
\end{theorem}

Since it will be useful in the sequel, we introduce here  the convex envelope of a function, which is the largest convex function that is bounded above by the given one.

\begin{definition}[convex envelope]
\label{def:convexenvelope}
Let $\Omega \subset \R^d$ be convex and $w : \bar \Omega \to \R$. The {\it convex envelope} of $w$, denoted by $\Gamma w$, is the largest convex function whose graph lies below the graph of $w$. It can be computed by 
\[
  \Gamma w(x) = \sup\left\{ L(x): \mbox{$L$ affine function and } L(y) \leq w(y)\ \ \forall y\in \bar\Omega \right\}.
\]
\end{definition}

To conclude our preliminary discussion we recall the Brunn-Minkowski inequality, a celebrated result in convex geometry. Given two compact sets $A,B$ of $\mathbb{R}^d$, we define their {\it Minkowski sum}
\begin{equation}\label{minkowski-sum}
A + B := \{ v + w \in \mathbb{R}^d:  v \in A \text{ and } w \in B \}.
\end{equation}
The Brunn-Minkowski inequality relates the Lebesgue measures of compact subsets $A,B$ of Euclidean space $\mathbb{R}^d$ with that of their Minkowski sum $A+B$.

\begin{lemma}[Brunn-Minkowski inequality]\label{BM}
Let $A$ and $B$ be two nonempty compact subsets of $\mathbb R^d$ for $d \geq 1$. Then the following inequality holds:
\[
  |A + B|^{1/d} \ge |A|^{1/d} + |B|^{1/d}.
\]
\end{lemma}

\section{Wide Stencil Finite Differences}
\label{sec:FD}

\epigraph{\tiny{\it Problems involving the classical linear partial differential equations of mathematical physics can be reduced to algebraic ones of a very much simpler structure by replacing the differentials by difference quotients on some (say rectilinear) mesh.}}{\tiny{R.~Courant, K.~Friedrichs and H.~Lewy \cite{MR0213764}}}

In this section we will study finite difference schemes that aim to approximate the viscosity solution, in the sense of Definition~\ref{def:viscosolMA}, of \eqref{eqn:MA}.

\subsection{A general framework for approximation schemes}
\label{sub:BarlesSouganidis}

Let us describe a general framework under which convergence of approximation schemes can be shown. Let 
$F : \bar\Omega \times \bbR \times \polS^d \to \R$
be elliptic in the sense of Definition~\ref{def:FLelliptic} and assume we wish to approximate the viscosity solution to \eqref{eq:BVPBarlesSouganidisStyle}. To do so, we introduce a family of approximation schemes, which are described by the collection of maps $\{ F_\eps \}_{\eps>0}$, where $F_\eps : \bar\Omega \times \bbR \times B(\bar\Omega) \to \bbR$, and $B(\bar\Omega)$ denotes the space of bounded functions on $\bar\Omega$. The parameter $\eps$ can be understood as a discretization parameter. With this family at hand, we seek for $u_\eps \in B(\bar\Omega)$ such that
\begin{equation}
\label{eq:BarlesSougScheme}
  F_\eps(x,u_\eps(x),u_\eps) = 0, \quad \mbox{in } \bar\Omega.
\end{equation}
We assume that the approximation schemes satisfy the following assumptions:
\begin{enumerate}
  \item {\it Monotonicity}: For all $\eps>0$, $x \in \bar \Omega$, $t \in \bbR$, and $u,v \in B(\bar\Omega)$ such that $u \geq v$ we have that
  \begin{equation}
  \label{eq:BarlesSougMonotone}
    F_\eps(x,t,u) \geq F_\eps(x,t,v).
  \end{equation}
  \item {\it Stability}: There is $\eps_0>0$ such that if $\eps<\eps_0$, the scheme \eqref{eq:BarlesSougScheme} has a unique solution and there is a constant, independent of $\eps$, such that
  \begin{equation}
  \label{eq:BarlesSougStable}
    \| u_\eps \|_{L^\infty(\Omega)} \leq C.
  \end{equation}

  \item {\it Consistency}: For all $x_0 \in \bar\Omega$ and $\varphi \in C^2(\bar\Omega)$ we have
  \begin{subequations}
  \label{eq:BarlesSougConsistent}
    \begin{equation}
    \label{eq:BarlesSougConsistentUSC}
      \limsup_{\eps\downarrow 0, y \to x_0, \xi \to 0} F_\eps(y,\varphi(y) + \xi, \varphi + \xi) \leq F_\star(x_0,\varphi(x_0),D^2 \varphi(x_0)) 
    \end{equation}
    \begin{equation}
    \label{eq:BarlesSougConsistentLSC}
      \liminf_{\eps\downarrow 0, y \to x_0, \xi \to 0} F_\eps(y,\varphi(y) + \xi, \varphi + \xi) \geq F^\star(x_0, \varphi(x_0), D^2 \varphi(x_0)).
    \end{equation}
  \end{subequations} 
\end{enumerate}

The main convergence result in this framework is the following; see \cite[Theorem 2.1]{BarlesSoug91}.

\begin{theorem}[Barles--Souganidis]
\label{thm:BarlesSouganidis}
Assume that the family of approximation schemes \eqref{eq:BarlesSougScheme} is monotone, stable and consistent, in the sense of \eqref{eq:BarlesSougMonotone}, \eqref{eq:BarlesSougStable}, and \eqref{eq:BarlesSougConsistent}, respectively. Assume, in addition, that problem \eqref{eq:BVPBarlesSouganidisStyle} has a comparison principle {in the sense of Definition~\ref{def:comparisonVisco}}. Then, as $\eps \downarrow 0$, the functions $u_\eps$, solution of \eqref{eq:BarlesSougScheme} converge locally uniformly to $u$, solution of \eqref{eq:BVPBarlesSouganidisStyle}.
\end{theorem}
\begin{proof}
Define $\overline{u}, \underline{u} \in B(\bar\Omega)$ by
\[
  \overline{u}(x) = \limsup_{y \to x, \eps \downarrow 0} u_\eps(y), \qquad \underline{u}(x) = \liminf_{y \to x, \eps \downarrow 0} u_\eps(y).
\]
Notice that, by stability, we obtain that these functions are well defined and bounded. In addition, we have that $\overline{u},\underline{u}$ are upper and lower semicontinuous, respectively.

The idea now  is to show that $\overline{u}$ is a subsolution and $\underline{u}$ is a supersolution of \eqref{eq:BVPBarlesSouganidisStyle}, for if that is the case we can invoke the comparison principle to see that $\overline{u} \leq \underline{u}$, and so that these must coincide with the viscosity solution of \eqref{eq:BVPBarlesSouganidisStyle}. This, in turn, implies the local uniform convergence of $u_\eps$ to $u$.

Let us then show that $\overline{u}$ is a subsolution. Let $\varphi \in C^2(\bar\Omega)$ and assume that $\overline{u}-\varphi$ has a local maximum at $x_0 \in \bar\Omega$ with $\overline{u}(x_0) = \varphi(x_0)$. It can be shown then that there are sequences $\{\eps_n\}_{n=1}^\infty \subset \bbR^+$ and $\{y_n\}_{n=1}^\infty \subset \bar\Omega$ such that $\eps_n \downarrow 0$, $y_n \to x_0$, $u_{\eps_n}(y_n) \to \overline{u}(x_0)$ and the sequence of functions $u_{\eps_n} - \varphi$ attains its maximum at $y_n$.

Notice now that, upon denoting $\xi_n = u_{\eps_n}(y_n) - \varphi(y_n)$, we get that $\xi_n \to 0$ and $u_{\eps_n}(x) - \varphi(x) \leq \xi_n$ locally.
Monotonicity then implies that
\begin{align*}
  0 = F_{\eps_n}( y_n, u_{\eps_n}(y_n), u_{\eps_n} ) 
  &= F_{\eps_n}(y_n, \varphi(y_n) + \xi_n, \varphi + (u_{\eps_n}-\varphi) )\\
  &\leq F_{\eps_n}(y_n, \varphi(y_n) + \xi_n, \varphi + \xi_n ),
\end{align*}
which by the consistency condition \eqref{eq:BarlesSougConsistentUSC} yields
\[
  F_\star(x_0,\varphi(x_0),D^2\varphi(x_0)) \geq 0,
\]
so that $\overline{u}$ is a subsolution.
\end{proof}

\begin{remark}[limitations]
\label{rem:BarlesSougLimitations}
We must remark that, although Theorem~\ref{thm:BarlesSouganidis} seems sufficiently general:
\begin{enumerate}
  \item It only provides sufficient conditions for convergence. There is no guideline towards the construction of monotone, consistent and stable finite difference schemes.
  
  \item This result, as is, {\it cannot} be applied to approximate viscosity solutions of the \MA equation \eqref{eqn:MA} directly. This is because, as pointed out in Section~\ref{subsub:viscosols}, the \MA operator is only elliptic over $\bar \Omega \times \bbR \times \polS^d_+$.
    
  \item The existence of a comparison principle in the sense of Definition~\ref{def:comparisonVisco} is assumed. Notice that, in \cite[Proposition 2.1]{JensenSmearsComparison} it is shown that, for a reformulation of the \MA problem as a Hamilton Jacobi Bellman equation (which will be discussed in Section~\ref{subsub:FengJensen}), if $f \equiv 0$, there cannot be a comparison principle for this problem. In other words, this is a highly nontrivial assumption.
\end{enumerate}
\eremk\end{remark}

Although not applicable to the \MA equation \eqref{eqn:MA}, one of the messages of Theorem~\ref{thm:BarlesSouganidis} is that monotonicity of a numerical scheme is a highly desirable property. Thus, it is necessary to explore how to construct monotone approximation schemes. In the context of finite difference schemes it was realized as early as in \cite{MotzkinWasow53} that, even for linear problems, monotonicity of a numerical scheme requires the use of wide stencils, which is rather problematic at points near the boundary. We refer the reader to \cite[Section 3.2]{NSWActa} for more details, and to \cite{MR3504992} for the construction of minimal stencils in two dimensions. For this reason, in the remaining of this section, we will consider wide stencil finite difference schemes to approximate the viscosity solution of \eqref{eqn:MA}.

\subsection{A variational characterization of the determinant}
\label{sub:vardet}

Let us provide a variational characterization of the determinant that will motivate most of the constructions which will come below. This was originally shown in \cite[Lemma 2]{FroeseOberman11}.

\begin{lemma}[characterization of the determinant]
\label{lem:detisvar}
Let $A$ be a symmetric positive definite $d\times d$ matrix and let
\[
  \calV= \left\{ \{\bw_i\}_{i=1}^d \subset \bbR^d : \bw_i\cdot \bw_j = \delta_{i,j} \right\},
\]
be the set of all orthonormal bases of $\bbR^d$. Then we have that
\[
  \det A = \min_{ \{\bw_i\}_{i=1}^d \in \calV } \prod_{i=1}^d \bw_i \cdot A\bw_i.
\]
\end{lemma}
\begin{proof}
To shorten notation, let $M = \min_{ \{\bw_i\}_{i=1}^d \in \calV } \prod_{i=1}^d \bw_i \cdot A\bw_i$. Then let $\{\bv_j\}_{j=1}^d$ be an orthonormal set of eigenvectors of $A$ so that
\[
  \det A = \prod_{i=1}^d \bv_i \cdot A\bv_i \geq M.
\]

On the other hand, for $\{\bw_i\}_{i=1}^d \in \calV$, we can represent them in the basis of eigenvectors $\bw_i = \sum_{k=1}^d (\bw_i \cdot \bv_k) \bv_k$. We have
\begin{align*}
  -\log \prod_{i=1}^d \bw_i \cdot A\bw_i&=  -\sum_{i=1}^d \log (\bw_i \cdot A\bw_i) \\
  &= - \sum_{i=1}^d \log \left( \sum_{m=1}^d (\bw_i \cdot \bv_m) \bv_m \cdot  \sum_{k=1}^d (\bw_i \cdot \bv_k) A \bv_k\right) \\
  &= - \sum_{i=1}^d \log \left( \sum_{k=1}^d \lambda_k (\bw_i \cdot \bv_k)^2 \right),
\end{align*}
where $\sigma(A) = \{\lambda_k\}_{k=1}^d$ is the spectrum of $A$. Since $| \bw_i| = 1$ the term $\sum_{k=1}^d \lambda_k (\bw_i \cdot \bv_k)^2$ is a convex combination of the elements of $\sigma(A)$. Owing to the convexity of $t \mapsto -\log t$ we can apply Jensen's inequality to obtain that
\[
  -\log \prod_{i=1}^d \bw_i \cdot A\bw_i  \leq - \sum_{k=1}^d \log \lambda_k \sum_{i=1}^d (\bw_i \cdot \bv_k)^2 = - \sum_{k=1}^d \log \lambda_k = - \log \prod_{i=1}^d \lambda_i.
\]
As the function $t \mapsto -\log t$ is decreasing, we conclude that
\[
  \det A \leq \prod_{i=1}^d \bw_i \cdot A\bw_i,
\]
which since $\{\bw_i\}_{i=1}^d \in \calV$ was arbitrary implies $\det A \leq M$ and this concludes the proof.
\end{proof}

The previous result allows us to conclude that, if $\varphi \in C^2(\Omega)$ is convex, we can express the determinant of its Hessian at a point in terms of second directional derivatives, that is, if $x_0 \in \Omega$ we have
\[
  \det D^2 \varphi(x_0) = \min_{ \{\bw_i\}_{i=1}^d \in \calV } \prod_{i=1}^d \bw_i \cdot D^2\varphi(x_0) \bw_i = \min_{ \{\bw_i\}_{i=1}^d \in \calV } \prod_{i=1}^d \frac{\p^2 \varphi}{\p \bw_i^2}(x_0).
\]
Recall, in addition, that a solution to \eqref{eqn:MA} must be convex. To enforce convexity we then introduce the following operator
\begin{equation}
\label{eq:MAVariationalwithconvex}
  \MAop{\varphi}{}{}(x_0) = \min_{ \{\bw_i\}_{i=1}^d \in \calV } \left[ \prod_{i=1}^d\left(\frac{\p^2 \varphi}{\p \bw_i^2}(x_0) \right)^+ - \sum_{i=1}^d \left( \frac{\p^2 \varphi}{\p \bw_i^2}(x_0) \right)^- \right],
\end{equation}
where $x^+=\max\{x,0\}$ and $x^- = (-x)^+$ denote the positive and negative parts of $x$, respectively. Notice that, if $\varphi \in C^2(\bar\Omega)$ is convex, $\MAop{\varphi}{}{}=\det D^2\varphi$. The idea behind \eqref{eq:MAVariationalwithconvex} is that, if $D^2\varphi(x_0)$ has a negative eigenvalue, then there is $V \in \calV$ and $\bw \in V$ for which $\bw\cdot D^2\varphi(x_0)\bw < 0$. Thus,
\[
  \MAop{\varphi}{}{}(x_0) \leq 0 - (\bw\cdot D^2\varphi(x_0)\bw)^{-}<0.
\]
Consequently, $\varphi$ cannot be a solution to \eqref{eqn:MA} since, at $x_0$ we have
\[
  \det D^2\varphi(x_0) = f(x_0) \geq 0.
\]
These ideas are made rigorous in \cite[Lemma 5.6]{NochettoNtogkasZhangfirsttwoscale}.

\begin{proposition}[equivalence of operators]
\label{prop:equiv}
Let $f \in C(\Omega)$ with $f \geq 0$. The function $u \in C(\bar\Omega)$ is a convex viscosity solution of \eqref{eqn:MA} in the sense of Definition~\ref{def:viscosolMA} if and only if it is a viscosity solution, in the sense of Definition~\ref{def:viscosol}, of the following problem
\begin{equation}
\label{eq:FvMA}
  F_{vMA}(x, u(x), D^2u(x)) = 0
  \end{equation}
  with
  \begin{align*}
   F_{vMA}(x, u(x), D^2u(x)) = \begin{dcases}
                                                                      \MAop{u}{}{}(x)-f(x), & x \in \Omega,\\
                                                                      g(x) - u(x), & x \in \p\Omega.
                                                                    \end{dcases}
\end{align*}
\end{proposition}

One of the advantages of formulation \eqref{eq:FvMA} is that it has a comparison principle.

\begin{proposition}[comparison principle for the $F_{vMA}$ operator]
\label{prop:detisvarhascompare}
The operator $F_{vMA}$, defined in \eqref{eq:FvMA} has a comparison principle in the sense of Definition~\ref{def:comparisonVisco}.
\end{proposition}
\begin{proof}
It follows from the fact that the operator $F_{vMA}$ satisfies the structural assumptions given, for instance, in \cite[Theorem 3.3]{CIL}.
\end{proof}

The characterization of the determinant given in Lemma~\ref{lem:detisvar} will be the basis of many of the wide stencil schemes we will describe below.

\subsection{Wide stencil finite difference schemes}
\label{sub:Oberman}

Let us describe the first class of methods that exploit the characterization described in Lemma~\ref{lem:detisvar} via the operator introduced in \eqref{eq:MAVariationalwithconvex} as originally proposed in \cite{FroeseOberman11}. Let $h>0$ be a (spatial) discretization parameter and assume that, up to a linear change of variables, our domain $\Omega$ is discretized on a Cartesian grid. In other words, we let
\[
  \bar\Omega_h = \bar \Omega \cap \Z_h^d, \quad \Z_h^d = \left\{ h \be : \be \in \Z^d \right\}, \quad \partial\Omega_h = \partial\Omega \cap \Z_h^d, \quad \Omega_h = \bar\Omega_h \setminus \partial\Omega_h.
\]
We set $X_h$ as the space of {\it grid functions}, that is the collection of functions $w_h : \bar\Omega_h \to \bbR$. 

Given $\be \in \Z^d$ we call the point $x_h \in \Omega_h$ {\it interior} with respect to $\be$ if $x_h \pm h \be \in \bar\Omega_h$. We will also say that a point is interior with respect to a subset of $S \subset \Z^d$ if it is interior with respect to all elements of $S$.

Given $\be \in \Z^d$ and an interior point $x_h$, we define the {\it second difference} in the direction $\be$ to be the operator
\begin{equation}
\label{eq:2nddiff}
  \Del_{\be} w_h(x_h) = \frac1{|e|^2h^2} \left( w_h(x_h + h \be) - 2w_h(x_h) + w_h(x_h - h\be) \right).
\end{equation}
When $x_h$ is not interior with respect to $\be$, it essentially means that $x_h$ is close to $\p\Ome$. Owing to the convexity of $\Omega$, there are unique $\rho_\pm \in (0,1]$ such that $x_h \pm \rho_\pm h \be \in \p\Ome$. Thus, we can use the boundary condition \eqref{eqn:MA2} to extend this definition as
\begin{equation}
\label{eq:2nddiffbdry}
\begin{aligned}
  \Del_e w_h(x_h) =
  \frac2{(\rho_+ + \rho_-)|\be|^2h^2} &\left( \frac{\tilde{g}(x_h+\rho_+h\be) - w_h(x_h)}{\rho_+} \right. \\ &-\left. \frac{w_h(x_h) - \tilde{g}(x_h+\rho_-h\be)}{\rho_-}\right),
\end{aligned}
\end{equation}
where $\tilde{g}$ is either the boundary condition, or an interpolant of $w_h$ based on neighboring nodes. With these notions at hand, we would like to define the discretization of the operator $\MAop{\cdot}{}{}$, introduced in \eqref{eq:MAVariationalwithconvex}, as
\[
  \MAop{w_h}{h}{WS}(x_h) = \min_{\{\bw_i\}_{i=1}^d \in \calV } \prod_{i=1}^d\left(\Del_{\bw_i} w_h(x_h) \right)^+ .
\]
Notice, however, that the given expressions may not be defined for all $\calV$, as the points $x_h \pm h \bw_i$ may not belong to $\bar\Omega_h$. {Even if they did, it may be very computationally expensive to compute these directional differences at all the nodes. For these reasons}, we also need to introduce a discretization of $\calV$. To this end we introduce a finite subset $\calG_\theta \subset (\Z^d)^d$ such that, if $\{\bnu_i\}_{i=1}^d \in \calG_\theta$ then the vectors $\bnu_i$ are pairwise orthogonal. We call this the {\it directional} discretization of the \MA operator and parametrize it by $\theta>0$. Thus we define the operator
\begin{equation}
\label{eq:MAhtheta}
  \MAop{w_h}{h,\theta}{WS}(x_h) = \min_{\{\bnu_i\}_{i=1}^d \in \calG_\theta} \prod_{i=1}^d\left(\Del_{\bnu_i} w_h(x_h) \right)^+ .
\end{equation}

With this notation at hand, we define the {\it wide stencil} finite difference approximation scheme of \eqref{eqn:MA} as: Find $u_h \in X_h$ such that
\begin{subequations}
\label{eq:WSFD}
  \begin{align}
    \label{eq:WSFDin}
    \MAop{u_h}{h,\theta}{WS}(x_h) &= f(x_h), \quad \forall x_h \in \Omega_h, \\
  \label{eq:WSFDbc}
    u_h(x_h) &= g(x_h), \quad \forall x_h \in \p\Omega_h.
  \end{align}
\end{subequations}

\begin{remark}[variant]
\label{rem:variantWS}
We could have also introduced another wide stencil operator via
\[
  \MAop{w_h}{h,\theta}{WS}(x_h) = \min_{\{\bnu_i\}_{i=1}^d \in \calG_\theta} \left[\prod_{i=1}^d\left(\Del_{\bnu_i} w_h(x_h) \right)^+ - \sum_{i=1}^d \left( \Del_{\bnu_i} w_h(x_h) \right)^- \right],
\]
see \eqref{eq:MAVariationalwithconvex}.
\eremk\end{remark}

\begin{remark}[a regularized version]
Notice that, owing to the presence of the $\min$ and $\max$ operator in the definition of \eqref{eq:MAhtheta}, this operator is not differentiable. This may make it difficult to efficiently solve the ensuing nonlinear systems, since Newton methods are not directly applicable. One could, instead, use semismooth Newton methods \cite{MR1972219} since these operators are slant differentiable; see \cite[Lemma 3.1]{MR1972219}. However, if we insist in dealing with smooth operators, \cite[Section 3.5]{FroeseOberman11} introduces a regularized version of $\MAop{\cdot}{h,\theta}{WS}$ given by
\[
  \MAop{w_h}{h,\theta,\delta}{WS}(x_h) = \underset{\{\bnu_i\}_{i=1}^d \in \calG_\theta}{{\min}^\delta} \prod_{i=1}^d\left(\Del_{\bnu_i} w_h(x_h) \right)^{+,\delta},
\]
where {
\begin{align*}
  {\max}^\delta\{x,y\} &= \frac12 \left( x+y + \sqrt{(x-y)^2 + \delta^2} \right), \\
  {\min}^\delta\{x,y\} &= \frac12 \left( x+y - \sqrt{(x-y)^2 + \delta^2} \right), \\
  {\min}^\delta\{x_1, \ldots, x_n \} &= {\min}^\delta\{{\min}^\delta\{x_1,\ldots,x_{n-1}\},x_n\},
\end{align*}
}and $x^{+,\delta} = \max^\delta\{x,0\}$. The properties of this operator are similar to those of $\MAop{\cdot}{h,\theta}{WS}$.
\eremk\end{remark}

\begin{remark}[two dimensions]
Given $A \in \polS^d$ we have the classical Rayleigh--Ritz relations
\[
  \lambda_m(A) = \min_{\bw \in \bbR^d} \frac{ \bw\cdot A\bw}{|\bw|^2} = \min\sigma(A),  \qquad \lambda_M(A) = \max_{\bw \in \bbR^d} \frac{ \bw\cdot A\bw}{|\bw|^2} = \max\sigma(A),
\]
so that, if $d=2$, we have that
\[
  \det A = \min_{\bw \in \bbR^2} \frac{ \bw\cdot A\bw}{|\bw|^2} \max_{\bw \in \bbR^2} \frac{ \bw\cdot A\bw}{|\bw|^2}.
\]
This relation was used in \cite{Oberman08} to introduce a two dimensional scheme via
\[
  \MAop{w_h}{h,\theta}{WS,2d}(x_h) = \min_{\bnu_i \in \{\bnu_j\}_{j=1}^d \in \calG_\theta} \left(\Del_{\bnu_i} w_h(x_h) \right)^+ \max_{\bnu_i \in \{\bnu_j\}_{j=1}^d \in \calG_\theta} \left(\Del_{\bnu_i} w_h(x_h) \right)^+.
\]
Note that, although similar to \eqref{eq:WSFD}, these operators are different. This was illustrated in \cite[Section 3.4]{FroeseOberman11} with the following example: Let
\[
  w(x_1,x_2) = x_1^2+x_2^2+x_1^2x_2^2,
\]
which is convex in a neighborhood of the origin, and 
\[
  \calG_\theta = \left\{ \left\{ \begin{pmatrix} 1 \\ 0\end{pmatrix}, \begin{pmatrix} 0 \\ 1\end{pmatrix}\right\}, 
      \left\{ \begin{pmatrix} 1 \\ 1\end{pmatrix}, \begin{pmatrix} -1 \\ 1\end{pmatrix}\right\}
  \right\}.
\]
Computing each of the operators over these directions yields
\[
  \MAop{w}{h,\theta}{WS,2d}(0,0) = 4+2h^2, \qquad \MAop{w}{h,\theta}{WS}(0,0) = 4.
\]
Notice however, that since both operators are consistent with order $\calO(h^2)$ we have that, {for a convex function $v$},
\[
  \left| \MAop{v}{h,\theta}{WS,2d}(x_h) - \MAop{v}{h,\theta}{WS}(x_h) \right| = \calO(h^2), \quad \forall x_h.
\]
\eremk\end{remark}

The analysis of method \eqref{eq:WSFD} will be a particular case of the methods and analyses presented in Section~\ref{sub:twoscale}. We just comment that, even for smooth solutions, wide stencils are required in this scheme to assert consistency. Let us illustrate this in a simple case where there is no boundary conditions and in two dimensions ($d=2$). In other words, given $x_0 \in \Omega$ we assume that it is an interior point for any $\be \in \Z^2$. Let now $\varphi(x) = \tfrac12 x\cdot M x$ be a convex quadratic, so that
\[
  \Del_{\be} \varphi(x_0) = \frac1{|\be|^2} \be \cdot M \be,
\]
and therefore
\[
 \MAop{\varphi}{h,\theta}{WS}(x_0) = \min_{\{\bnu_1,\bnu_2\} \in \calG_\theta} \frac1{|\bnu_1|^2|\bnu_2|^2} \left(\bnu_1 \cdot M \bnu_1\right)\left( \bnu_1 \cdot M \bnu_2 \right),
\]
independently of $x_0$ and the mesh size. At this point we need to recall that there is $\{\bw_1,\bw_2\} \in \calV$, namely the normalized eigenvectors of $M$, for which
\[
  \det D^2 \varphi = \det M = \left( \bw_1 \cdot M \bw_1 \right)\left( \bw_2 \cdot M \bw_2 \right).
\]
Notice finally, that once $\bw_1$ is determined, $\bw_2 = \bw_1^\perp$ is obtained by a rotation. In conclusion, to assert consistency, given a $\bw \in \R^2$ in the unit sphere, for every $\delta>0$ we must be able to find $\be \subset \Z^2$ such that
\begin{equation}
\label{eq:eigenhasrationalapprox}
  \left| \bw - \frac1{|\be|}\be \right|  < \delta.
\end{equation}
Indeed, if we denote by $\be_1$ the vector that satisfies this property with respect to $\bw_1$, then $\be_2 = \be_1^\perp$ does so for $\bw_2$. Let now $\bnu_i = \tfrac1{|\be_i|} \be_i$ for $i=1,2$. Then we have that
\[
  \left| \det M - \left(\bnu_1 \cdot M \bnu_1\right)\left( \bnu_2 \cdot M \bnu_2 \right) \right| \leq C(\Lambda) \delta,
\]
where $C(\Lambda)$ is a constant that depends polynomially on $\Lambda$, the maximal eigenvalue of $M$.

Notice that, since $\be_i \in \Z^2$, then $\bnu_i \in \Q^2$, so finding points that satisfy \eqref{eq:eigenhasrationalapprox} is the problem of {\it rational approximation in the sphere}. While how to actually find such points is beyond our discussion here, what we are interested in is the size of $|\be|$, which would serve as an estimate of the stencil size that guarantees convergence. The following result is a specialization of \cite[Lemma 2.1]{MR2425007} to the two dimensional case; we refer the reader to this reference a proof, its generalization to $d > 2$, and to the case of rational approximation orthogonal matrices {which is of interest when finding elements of $\calG_\theta$}.

\begin{proposition}[rational approximation]
\label{prop:diophantine}
Let $\bw \in \R^2$ be such that $|\bw|=1$. Then, for every $\delta>0$, there exists $\bnu \in \Q^2$ such that $|\bnu|=1$ and
\[
  | \bw - \bnu | < \delta.
\]
Moreover, if $\bnu = (p_1/q_1, p_2/q_2)^\intercal$ with $p_1,p_2 \in \Z$ and $q_1,q_2 \in \N$ then we have that
\[
  0 < q_i \leq \frac{64}{\delta^2}.
\]
\end{proposition}

Now, for a given $\bw \in \R^2$, let $\bnu$ be as in Proposition~\ref{prop:diophantine}. This means that $\be = \textrm{hcf}(q_1,q_2) \bnu \in \Z^2$ is the smallest vector parallel to $\bnu$ that satisfies \eqref{eq:eigenhasrationalapprox}
(here,  $\textrm{hcf}(q_1,q_2)$ denotes the highest common factor of $q_1$ and $q_2$).  Consequently, we have that, generically
\[
  |\be| \leq C \ \textrm{hcf}(q_1,q_2) \leq C\max\{q_1,q_2\} \leq \frac{C}{\delta^2}.
\]
In conclusion, the size of the stencil must grow unboundedly
if we restrict ourselves to Cartesian meshes.

\subsection{Filtered schemes}
\label{sub:Filtered}

The estimates on the stencil size of the previous section are rather pessimistic. This is because they are not assuming anything but convexity of the solution. On the other hand, say in the two dimensional case ($d=2$), a standard nine point stencil finite difference approximation can be proposed
\begin{equation}
\label{eq:MAFD}
  \MAop{w_h}{h}{FD}(x_h) = \Del_{(1,0)}w_h(x_h) \Del_{(0,1)}w_h(x_h) - \left( \mathring{\Del}_{(1,1)}w_h(x_h) \right)^2,
\end{equation}
where, if $z_h = (x_1,x_2)^\intercal$, then
\begin{align*}
  \mathring{\Del}_{(1,1)}w_h(z_h) &= \frac1{2h} \Big( \frac{w_h(x_1 + h,x_2 + h) - w_h(x_1 - h,x_2 + h)}{2h}\\
  &\qquad\qquad - \frac{w_h(x_1 + h,x_2 - h) - w_h(x_1 - h,x_2 - h)}{2h}\Big).
\end{align*}
This formula easily extends to higher dimensions.

It is not difficult to see that $\MAop{\cdot}{h}{FD}$ has second order consistency, even for nonconvex functions. However, it is not monotone, even if one forgets about boundary conditions. Thus, it does not perform well when used to discretize problems that have singular solutions.

Reference \cite{MR2745457} takes advantage of the simplicity of \eqref{eq:MAFD} and the robustness of a wide stencil scheme by proposing a {\it hybrid} scheme.
Locally, it is a convex combination of each one of these schemes, where the weighting is chosen depending on the expected behavior of the solution. At points where the solution should be smooth 
the simple scheme \eqref{eq:MAFD} is used, whereas if the solution is {expected to be} singular
the robustness of \eqref{eq:MAhtheta} is better suited to capture this behavior. Summing up, the following discretization is used
\begin{equation}
\label{eq:MAhybrid}
  \begin{aligned}
    \MAop{w_h}{h}{H}(x_h) &= \omega(x_h) \MAop{w_h}{h}{FD}(x_h) \\ &+ (1-\omega(x_h))\MAop{w_h}{h,\theta}{WS}(x_h).
  \end{aligned}
\end{equation}
Here $\omega \in C(\bar \Omega, [0,1])$ is a weighting function  defined {\it a priori} from the data as follows: For $\epsilon >0$ we let $\Omega_\epsilon$ be a neighborhood of the set where the solution $u$ may be singular, that is,
\[
  \Omega_\epsilon = \left\{ x \in \Omega : 0 \leq f(x) < \epsilon \right\} \cup \left\{ x \in \partial\Omega: g \notin C^{2,\alpha}(U_x), \text{ or } U_x \cap \partial\Omega \text{ is flat} \right\},
\]
where $U_x$ is a neighborhood of the point $x$.
We then set $\omega \equiv 0$ in $\Omega_\epsilon$ and one away from it. This scheme was tested in \cite{MR2745457} for a series of cases, ranging from smooth to singular solutions, and computational experiments suggested that this method is robust and accurate.

This method, however, has a major drawback. The tunable function $\omega$ must be described by the user, and its values depend on the behavior of the problem data. For this reason in \cite{MR3033017} it was proposed that instead the difference
\[
  \left| \MAop{w_h}{h,\theta}{WS}(x_h) -\MAop{w_h}{h}{FD}(x_h) \right|,
\]
be used as an {\it a posteriori} indicator of accuracy. In regions where this difference is small, it is expected that the solution is smooth, whereas when this is large one expects singularities. On the basis of this, we can choose which scheme to apply. The way to measure this difference is by introducing a filter.

\begin{definition}[filter]
\label{def:filter}
A {\it filter} is a function $S\in C_0(\R)$ such that $S(t) = t$ in a neighborhood of the origin.
\end{definition}

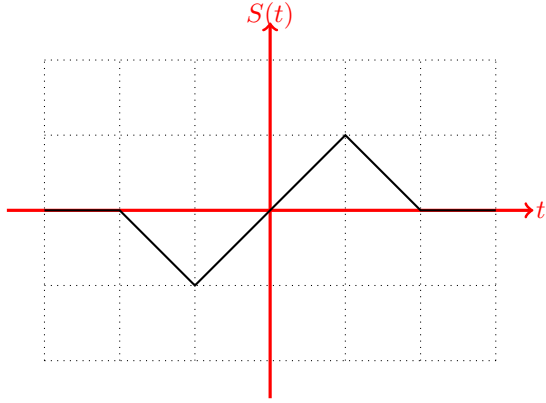
\begin{figure}
\begin{center}
  \begin{tikzpicture}
    \draw[step = 1, dotted] (1,1) grid (7,5);
    \draw[red, very thick, ->] (0.5,3) -- (7.5,3);
      \node[red] at (7.6,3) {$t$};
    \draw[red, very thick, ->] (4,0.5) -- (4,5.5);
      \node[red] at (4,5.6) {$S(t)$};
    \draw[thick] (1,3) -- (2,3) -- (3,2) -- (5,4) -- (6,3) -- (7,3);
  \end{tikzpicture}
\end{center}
\caption{The function $S$ defined in \eqref{eq:exoffilter} is a filter.}
\label{fig:AFilter}
\end{figure}

For instance, the function
\begin{equation}
\label{eq:exoffilter}
  S(t) = \begin{dcases}
           x, & |x| \leq 1 \\
           0, & |x| \geq 2, \\
           2-x, & 1< x< 2, \\
           -x-2, & -2 < x < -1
         \end{dcases}
\end{equation}
depicted in Figure~\ref{fig:AFilter} is a possible filter, see \cite[Figure 1.1 and (1.3)]{MR3033017}. With this at hand, a {\it filtered} operator can be defined via
\begin{equation}
\label{eq:Filterscheme}
  \begin{aligned}
    \MAop{w_h}{h}{F}(x_h) &= \MAop{w_h}{h,\theta}{WS}(x_h) \\ &+ h^\alpha S\left( \frac{\MAop{w_h}{h}{FD}(x_h) - \MAop{w_h}{h,\theta}{WS}(x_h)}{h^\alpha} \right),
  \end{aligned}
\end{equation}
where $\alpha \in (0,2]$ is to be chosen by the user. A filtered scheme seeks  $u_h \in X_h$ such that
\begin{subequations}
\label{eq:FilterFD}
  \begin{align}
    \label{eq:Filterin}
    \MAop{u_h}{h}{F}(x_h) &= f(x_h), \quad \forall x_h \in \Omega_h, \\
  \label{eq:Filterbc}
    u_h(x_h) &= g(x_h), \quad \forall x_h \in \p\Omega_h.
  \end{align}
\end{subequations}

\begin{remark}[consistency]
Recall that \cite{MR3416386,Oberman06} a monotone scheme cannot be more than second order accurate. Notice, in addition, that by construction we have 
\[
  \left| \MAop{w_h}{h}{F}(x_h) - \MAop{w_h}{h,\theta}{WS}(x_h) \right| \leq h^\alpha,
\]
so that a filtered scheme is also consistent, up to second order. Moreover, if the parameter $\alpha$ is chosen smaller than the consistency order of both the wide stencil, and the finite difference scheme, and the mesh size $h$ is sufficiently small, it can be shown that
\[
  \MAop{\varphi}{h}{F}(x) = \MAop{\varphi}{h}{FD}(x),
\]
whenever $\varphi$ is sufficiently smooth. These two observations serve as a guideline for the choice of $\alpha$.
\eremk\end{remark}

\begin{remark}[motivation]
The construction of a filtered scheme seems to be motivated by similar constructions for conservation laws and first order Hamilton Jacobi equations. For instance, \cite{MR1314597} shows the convergence of filtered finite difference schemes (constructed in a similar way), for Hamilton Jacobi equations. {In the realm of hyperbolic conservation laws, several types of limiters or artificial viscosity methods \cite{MR2787948,MR3167449,MR3204837,MR3860123} have been derived from these ideas.}
\eremk\end{remark}

As a step towards the analysis of schemes like \eqref{eq:FilterFD}, reference \cite{MR3033017} introduced a class of schemes called {\it nearly monotone}, and showed that the theory of Section~\ref{sub:BarlesSouganidis} also applies to them. To show this, we begin with a definition.

\begin{definition}[nearly monotone]
\label{def:nearmonotone}
The family of approximation schemes\\ $\{F_\eps\}_{\eps>0}$ where $F_\eps : \bar\Omega \times \R \times B(\bar\Ome)$ is called {\it nearly monotone}, if every $F_\eps$ can be written as
\[
  F_\eps = F_\eps^M + F_\eps^P,
\]
where $F_\eps^M$ is monotone in the sense of \eqref{eq:BarlesSougMonotone}, and the function $F_\eps^P$, called a perturbation, satisfies
\[
  \lim_{\eps \downarrow 0}|F^P_\eps(x,t,v)| = 0,
\]
uniformly on bounded subsets of $\bar\Omega \times \R \times B(\bar\Ome)$.
\end{definition}

The convergence of nearly monotone schemes closely follows that of monotone schemes.

\begin{corollary}[convergence]
\label{cor:BSnearmonotone}
Let $\{F_\eps\}_\eps$ be a family of approximation schemes, that is nearly monotone, in the sense of Definition~\ref{def:nearmonotone}; consistent, in the sense of \eqref{eq:BarlesSougConsistent}; and stable in the sense of \eqref{eq:BarlesSougScheme}. Assume, in addition, that problem \eqref{eq:BVPBarlesSouganidisStyle} has a strong comparison principle. In this setting we have that, as $\eps \downarrow 0$, the functions $u_\eps$, solutions of $F_\eps(x,u_\eps(x),u_\eps) = 0$ converge locally uniformly to $u$, solution of \eqref{eq:BVPBarlesSouganidisStyle}.
\end{corollary}
\begin{proof}
The proof is a small variation on the proof of Theorem~\ref{thm:BarlesSouganidis}. Indeed, with the notation of this proof, we have
\begin{align*}
  0 &= F_{\eps_n}( y_n, u_{\eps_n}(y_n), u_{\eps_n} )\\ 
&  = F_{\eps_n}^M(y_n, \varphi(y_n) + \xi_n, \varphi + (u_{\eps_n}-\varphi) ) + F_{\eps_n}^P( y_n, u_{\eps_n}(y_n), u_{\eps_n} ) \\
    &\leq F_{\eps_n}^M(y_n, \varphi(y_n) + \xi_n, \varphi + \xi_n ) + F_{\eps_n}^P( y_n, u_{\eps_n}(y_n), u_{\eps_n} ).
\end{align*}
The stability of the scheme allows us to invoke the fact that the perturbation vanishes in the limit. Consequently, we still have that $\overline{u}$ is a subsolution.
\end{proof}

Notice that the same considerations made in Remark~\ref{rem:BarlesSougLimitations} apply in this setting.

\subsection{Lattice basis reduction scheme}
\label{sub:LBR}

Let us now discuss a two dimensional method, which was introduced in \cite{BCM16} and is termed the {\it lattice basis reduction} scheme. The aim of this scheme is, for a given stencil, to obtain a different way to compute the determinant, so that the scheme is more accurate.
We begin with a definition.

\begin{definition}[superbasis]
\label{def:supberbasis}
We will say that a {\it basis} of $\Z^2$ is a pair of vectors $(\be_1,\be_2) \in (\Z^2)^2$ that satisfy $\left| \det (\be_1, \be_2) \right| = 1$. A {\it superbasis} of $\Z^2$ is a triple $(\be_0,\be_1,\be_2) \in (\Z^2)^3$ such that $(\be_1,\be_2)$ is a basis and $\be_0+\be_1+\be_2 = 0$.
\end{definition}

We will call a {\it stencil} a finite subset of $\Z^2\setminus\{0\}$ that is symmetric around the origin. To a stencil $S$ we associate the set of superbases
\[
  Y(S) = \left\{ (\be_0,\be_1,\be_2) \in S^3: \left|\det(\be_1,\be_2)\right|=1, \ \be_0+\be_1+\be_2 = 0 \right\}.
\]

With these notations at hand, we define the {\it lattice basis reduction} \MA operator
\begin{align}
\label{eq:LBRop}
&  \MAop{w_h}{h,S}{LBR}(x_h)\\
&\nonumber\qquad = \min_{(\be_0,\be_1,\be_2) \in Y(S)} \gamma \left( \left(\Del_{\be_0} w_h(x_h) \right)^+, \left(\Del_{\be_1} w_h(x_h) \right)^+, \left(\Del_{\be_2} w_h(x_h) \right)^+ \right),
\end{align}
where
\[
  \gamma(\delta_0, \delta_1, \delta_2) = \begin{dcases}
                                           \delta_{i+1}\delta_{i+2}, & \delta_i \geq \delta_{i+1} + \delta_{i+2},\\ 
                                           \frac12(\delta_0\delta_1+ \delta_1\delta_2 + \delta_0 \delta_2) - \frac14( \delta_0^2 + \delta_1^2 + \delta_2^2 ), & \text{otherwise}.
                                         \end{dcases}
\]
This allows us to introduce the following scheme: Find $u_h \in X_h$ such that
\begin{subequations}
\label{eq:LBR}
  \begin{align}
    \label{eq:LBRin}
    \MAop{u_h}{h,S}{LBR}(x_h) &= f(x_h), \quad \forall x_h \in \Omega_h, \\
  \label{eq:LBRbc}
    u_h(x_h) &= g(x_h), \quad \forall x_h \in \p\Omega_h.
  \end{align}
\end{subequations}

The motivation for this, at first glance obscure, definition of the operator $\MAop{\cdot}{h,S}{LBR}$ is given in \cite[Remark 1.10]{BCM16}. Let $Y = (\be_0,\be_1,\be_2) \in Y(S)$ and notice that for any point $x_h$ that is interior interior with respect to $Y$, we have that the convex hull of $\{ x_h \pm h \be_i \}_{i=0}^2$ is a hexagon. Given a function $w_h \in X_h$ we can associate to it its {\it local convex envelope}, that is the maximal convex function $\Gamma_{x_h,Y}w_h$ that is bounded from above by $w_h$ at the points $\{ x_h \pm h \be_i \}_{i=0}^2$. It is then possible to show that $\Gamma_{x_h,Y}w_h$ is a piecewise linear function over a particular triangulation of the aforementioned hexagon. Then we have that
\begin{align}\label{LBR_subdifferential}
  \gamma \left( \left(\Del_{\be_0} w_h(x_h) \right)^+, \left(\Del_{\be_1} w_h(x_h) \right)^+, \left(\Del_{\be_2} w_h(x_h) \right)^+ \right) = \left| \partial \Gamma_{x_h,Y}w_h(x_h)\right|,
\end{align}
which is consistent with the definition of the \MA operator in the sense of Alexandrov given in Definition~\ref{def:AlexSolution} and hints at the consistency of this scheme. 

The consistency analysis of the operator \eqref{eq:LBRop} hinges on the following definition.

\begin{definition}[$M$--obtuseness]
\label{def:Mobtuse}
Let $M \in \polS_+^2$.  We say that the superbasis $(\be_0,\be_1,\be_2)$ of $\Z^2$ is $M$--obtuse if and only if
\[
  \be_j \cdot M\be_i \leq 0, \quad \forall 0 \leq i < j \leq 2.
\]
\end{definition}

From this definition, a necessary and sufficient condition for consistency follows \cite[Theorem 1.9]{BCM16}.

\begin{theorem}[consistency]
Let $\varphi = \tfrac12 x \cdot M x$ be a convex quadratic polynomial. We have that
\[
  \MAop{\varphi}{h,S}{LBR}(x) = \det M, \quad \forall x
\]
if and only if $Y(S)$ contains a $M$--obtuse superbasis.
\end{theorem}
\begin{proof}
We will follow \cite[Section 2.1]{BCM16}. To simplify the discussion, we set
\begin{align*}
  \bDi &= \left\{ (a_0,a_1,a_2) \in \R^3: a_i \leq a_{i+1} + a_{i+2},\  i = 0,1, 2, \mod 3 \right\}, \\
  \gamma_1(a_0,a_1,a_2) &= \frac12(a_0a_1+ a_1a_2 + a_0 a_2) - \frac14( a_0^2 + a_1^2 + a_2^2 ).
\end{align*}
Notice that $\gamma(a_0,a_1,a_2) = \gamma_1(a_0,a_1,a_2)$ if and only if $(a_0,a_1,a_2) \in \bDi$, and that if that is not the case, then $\gamma(a_0,a_1,a_2) - \gamma_1(a_0,a_1,a_2) = \frac14(a_0-a_1-a_2)^2 >0$. In conclusion, we have that
\begin{equation}
\label{eq:setKandh1}
  \begin{dcases}
    \gamma(a_0,a_1,a_2) \geq \gamma_1(a_0,a_1,a_2), \\ 
    \gamma(a_0,a_1,a_2) = \gamma_1(a_0,a_1,a_2) \Leftrightarrow (a_0,a_1,a_2) \in \bDi.
  \end{dcases}
\end{equation}

Given a superbasis $(\be_0,\be_1,\be_2)$ define $\delta_i = \be_i \cdot M \be_i =\big(\Delta_{\be_i} \varphi(x_h)\big)^+$. For a permutation $(i,j,k)$ of $(0,1,2)$ we have
\[
  \delta_i - \delta_j - \delta_k = (\be_j+\be_k)\cdot M(\be_j+\be_k) - \be_j \cdot M \be_j - \be_k M \be_k = 2 \be_j \cdot M\be_k.
\]
Consequently, $(\delta_0,\delta_1,\delta_2) \in\bDi$ if and only if the superbasis $(\be_0,\be_1,\be_2)$ is $M$--obtuse.

Let $A$ be the linear transformation that maps $\be_1$ and $\be_2$ to $\ff_1 = (1,0)^\intercal$ and $\ff_2 = (0,1)^\intercal$, respectively. Then we must have that $\ff_0 = A \be _0 = (-1,-1)^\intercal$. 
Thus, $\delta_i = \be_i \cdot M\be_i = A^{-1} \ff_i \cdot M A^{-1} \ff_i$, and so
\[
  \gamma_1(\mu_0,\mu_1,\mu_2) = \det(A^{-\intercal} M A^{-1}).
\]
However, $\det A = | \det(\be_1,\be_2) / \det(\ff_1,\ff_2) | = 1$. Combining this with \eqref{eq:setKandh1} we obtain the claim.
\end{proof}

Essentially, the previous result shows that the operator $\MAop{\cdot}{h,S}{LBR}$ systematically overestimates the determinant of the Hessian for quadratic functions, and that we have equality if and only if the stencil $S$ contains a $M$--obtuse superbasis. For this reason, it is of interest to obtain conditions on the size of the stencil that guarantee that such a superbasis can be found. The following result is a restatement of \cite[Proposition 1.12]{BCM16}.

\begin{proposition}[stencil size estimate]
\label{prop:stencilsizeLBR}
The stencil
\[
  S = \left\{ \be \in \Z^2 : \gcd(\be) = 1, |\be|\leq 2\kappa \right\},
\]
contains a $M$--obtuse superbasis for every matrix $M \in \polS_+^2$ that satisfies 
\[
\|M \|_2 \| M^{-1} \|_2 \leq \kappa^2.
\]
\end{proposition}

Notice that the cardinality of the stencil stated in Proposition~\ref{prop:stencilsizeLBR} is quite large, approximately $\kappa^2$, and that if the solution degenerates, that is $\det D^2u(x_0) = 0$ at some point, then the stencil size must again grow unboundedly to maintain consistency.

\subsection{Discretization based on power diagrams}
\label{sub:Power}
In \cite{Mirebeau15} the following discretization of the \MA operator is proposed and analyzed. Let $S$ be a stencil such that $\Span S = \R^d$ and such that its elements have co-prime coordinates, that is, if $\be = (e_1,\ldots,e_d)^\intercal \in S$, then $\gcd(\be) = \gcd(e_1,\ldots,e_d) = 1$. We define
\begin{equation}
\label{eq:PowDiag}
  \MAop{w_h}{h,S}{PD}(x_h) = \left| \left\{ \bg \in \R^d : \forall \be \in S: 2 \bg \cdot \be \leq |\be|^2 \Del_{\be} w_h(x_h) \right\} \right|.
\end{equation}
Here, we denote the Lebesgue measure by $|\cdot|$. With this operator at hand, we define the problem: find $u_h \in X_h $ such that
\begin{subequations}
\label{eq:powdiag}
  \begin{align}
    \label{eq:powdiagin}
    \MAop{u_h}{h,S}{PD}(x_h) &= f(x_h), \quad \forall x_h \in \Omega_h, \\
  \label{eq:powdiagbc}
    u_h(x_h) &= g(x_h), \quad \forall x_h \in \p\Omega_h.
  \end{align}
\end{subequations}

Notice that the set entering the definition \eqref{eq:PowDiag} is a polytope. Efficient ways to compute the volume of a polytope are available. For instance, if the dimension is not too high (and recall that we are mostly interested in the cases $d=2$ or $d=3)$, one can first triangulate this polytope to then easily compute its volume.

Let us study the consistency of this scheme. To do so, we must introduce a definition.

\begin{definition}[Voronoi cells and facets]
\label{def:Voronoi}
Let $M\in \polS_+^d$. The Voronoi cell and facet are
\begin{align*}
  \Vor(M) &= \left\{ \bg \in \R^d: \forall \be \in \Z^d, \ \bg \cdot M \bg \leq (\bg - \be)\cdot M(\bg - \be) \right\}, \\
  \Vor(M,\be) &= \left\{ \bg \in \Vor(M): \bg \cdot M \bg = (\bg - \be)\cdot M(\bg - \be) \right\}.
\end{align*}
A $M$--Voronoi vector is an element $\be \in \Z^d \setminus\{0\}$ such that $\Vor(M,\be) \neq \emptyset$. It is a strict $M$--Voronoi vector if the facet $\Vor(M,\be)$ is $(d-1)$--dimensional.
\end{definition}

Now, the consistency of the operator defined in \eqref{eq:PowDiag} is as follows.

\begin{proposition}[consistency]
\label{prop:consistencyPowDiag}
Let $\varphi(x) = \tfrac12 x \cdot M x$ be a convex quadratic. Then we have that 
\[
  \MAop{\varphi}{h,S}{PD}(x) = \det M, \quad \forall x
\]
if and only if the stencil $S$ contains all the strict $M$--Voronoi vectors.
\end{proposition}
\begin{proof}
Let $\kappa = \sqrt{\| M \|_2 \| M^{-1} \|_2}$. We divide the proof in several steps.

\begin{enumerate}
  \item \label{Step1} Any point $\bg \in \Vor(M)$ must satisfy $|\bg| \leq \tfrac12 \kappa \sqrt{d}$. Any $M$--Voronoi vector $\be$ satisfies $|\be|\leq \kappa \sqrt{d}$ and has co-prime coordinates: \\
    Indeed, if $\bg \in \Vor(M)$, then let $\be_{\bg} \in \Z^d$ be obtained by rounding its coordinates to the nearest integer, so that $|\bg - \be_{\bg}| \leq \tfrac12 \sqrt{d}$. The estimate
    \[
     \frac1{\| M^{-1} \|_2} |\bg|^2 \leq \bg \cdot M \bg \leq (\bg - \be_{\bg})\cdot M(\bg - \be_{\bg}) \leq \| M \|_2 |\bg - \be_{\bg}|^2 \leq \frac{d}4 \| M \|_2
    \]
    yields the desired estimate. In addition, if $\be$ is a $M$--Voronoi vector, there is $\bg \in \Vor(M)$ for which $|\bg|= |\be - \bg|$ so that
    \[
      |\be| \leq 2 |\bg| \leq \kappa \sqrt{d}.
    \]
    Finally, to show that the coordinates must be co-prime consider $k \be \in \Z^d$ with $k\geq 2$ and notice that, for every $\bg \in \R^d$ we have
    \[
      (k\be - \bg)\cdot M(k\be - \bg) + (k-1) \bg \cdot M\bg = k (\be -\bg)\cdot M (\be-\bg) + (k^2-k) \be\cdot M \be.
    \]
    Consequently,
    \[
      (\be-\bg)\cdot M (\be - \bg) < \max\left\{ (k\be - \bg)\cdot M(k\be - \bg), \bg \cdot M\bg \right\},
    \]
    and $k\be$ cannot be a $M$--Voronoi vector.
    
  \item Let $E$ be the set of strict $M$--Voronoi vectors, then
  \[
    \Vor(M) \subset \left\{ \bg \in \R^d : \forall \be \in S: 2\bg \cdot M \be \leq \be \cdot M \be \right\},
  \]
  with equality if and only if $E \subset S$: \\
    Notice that $\bg \cdot M \bg \leq (\bg - \be) \cdot M (\bg - \be)$ is equivalent to saying that $2\bg \cdot M \be \leq  \be \cdot M \be$. This shows that $\Vor(M)$ is a convex polytope, defined by inequalities of this type where $\be$ runs over the set of strict $M$--Voronoi vectors. The bound established in the previous step shows that there can only be a finite number of them.

  \item $|\Vor(M)| = 1$: \\
    It follows from the observation that $\Vor(M)$ collects all elements $\bg \in \R^d$ that are closer to the origin (in the metric induced by the matrix $M$) than to any other point $\be \in \Z^d\setminus\{0\}$.
    
  \item Consistency: \\
    Recall that, for any $\be \in S$ we have that $|\be|^2\Del_{\be} \varphi(x) = \be \cdot M \be$. Consequently,
    \begin{align*}
      \MAop{\varphi}{h,S}{PD}(x) 
      &= \left| \left\{ \bg \in \R^d: \forall \be \in S: 2\bg \cdot \be \leq \be \cdot M \be \right\} \right|\\ 
      &= \left| M\left\{ \bg \in \R^d: \forall \be \in S: 2\bg \cdot M\be \leq \be \cdot M \be \right\} \right|.
    \end{align*}
    A combination of the second and third steps then yields
    \[
      \MAop{\varphi}{h,S}{PD}(x) \geq \det M |\Vor(M)| = \det M,
    \]
    with equality if $\Vor(M) = \{ \bg \in \R^d: \forall \be \in S: 2\bg \cdot \be \leq \be \cdot M \be \}$ with equality if $\Vor(M)$ contains all strict $M$--Voronoi vectors.
\end{enumerate}
This concludes the proof.
\end{proof}

Since the consistency of the operator $\MAop{\cdot}{h,S}{PD}$ requires the stencil to contain all strict Voronoi vectors, it is necessary to provide sufficient conditions for this to happen.

\begin{corollary}[stencil size estimate]
\label{cor:stencilsizePD}
Let $\kappa >0$ and define
\[
  S = \left\{ \be \in \Z^d: |\be| \leq \sqrt{d}\kappa, \ \gcd(\be) = 1\right\}.
\]
Let $\varphi(x) = \frac12 x\cdot M \cdot x$, then we have that
\[
  \MAop{\varphi}{h,S}{PD}(x) = \det M, \quad \forall x
\]
provided $\|M\|_2 \|M^{-1} \|_2 \leq \kappa^2$.
\end{corollary}
\begin{proof}
It immediately follows from the norm estimates given in Step~\ref{Step1} in the proof of Proposition~\ref{prop:consistencyPowDiag}.
\end{proof}

Let us now provide a convergence analysis of scheme \eqref{eq:powdiag}, which will follow from the framework provided in Section~\ref{sub:BarlesSouganidis}. To do so, we introduce the operator $F_{h,S}: \bar\Ome_h \times \R \times X_h \to \R$ via
\begin{equation}
\label{eq:PDasBS}
  F_{h,S}(x_h, t, w) = \begin{dcases}
                       \MAop{w}{h,S}{PD}(x_h) - f(x_h), & x_h \in \Omega_h, \\
                       g(x_h) - t, & x_h \in \p\Ome_h,
                     \end{dcases}
\end{equation}
and notice that \eqref{eq:powdiag} can be compactly written as
\[
  F_{h,S}(x_h,u_h(x_h), u_h) = 0, \ \forall x_h \in \bar\Omega_h.
\]
Let us also define the operator $F_S: \bar\Omega \times \R \times \polS_+^d \to \R$
\begin{equation}
\label{eq:PDcontasBS}
  F_S(x,t,M) = \begin{dcases}
                 |K(M)| - f(x), & x \in \Omega, \\
                 g(x) - t, & x \in \p\Ome,
               \end{dcases}
\end{equation}
where
\[
  K(M)= \{ \bv \in \R^d: \forall \be \in S,\ 2 \bv \cdot \be \leq \be \cdot M\be \} .
\]
Notice that, if $D^2u(x_0)$ exists for all $x_0 \in \Ome$ and its eigenvalues are properly bounded, see Corollary~\ref{cor:stencilsizePD}, we have that
\[
  \det D^2 u(x_0) - f(x_0)= F_S(x_0,u(x_0),D^2u(x_0)).  
\]
For this reason, we will consider the problem: find $u$ that is a viscosity solution of
\begin{equation}
\label{eq:PDcontinuouscompact}
  F_S(x,u(x),D^2u(x)) = 0, \quad x \in \bar\Ome.
\end{equation}

Following \cite[Section 2.3]{Mirebeau15} we will now show the convergence of scheme \eqref{eq:powdiag} via Theorem~\ref{thm:BarlesSouganidis}. To do so, we must show that scheme \eqref{eq:powdiag} is monotone, consistent, and stable in the sense of \eqref{eq:BarlesSougMonotone}, \eqref{eq:BarlesSougConsistent}, and \eqref{eq:BarlesSougStable}, respectively. We have shown consistency in Proposition~\ref{prop:consistencyPowDiag}. For stability, we refer the reader to \cite[Section 2.2]{Mirebeau15}, where stability is shown by proving global convergence of a damped Newton algorithm. We will focus then on the monotonicity of the scheme.

\begin{proposition}[monotonicity]
\label{prop:monotonePD}
The operator $F_{h,S}$, defined in \eqref{eq:PDasBS} is monotone in the sense of \eqref{eq:BarlesSougMonotone}.
\end{proposition}
\begin{proof}
Notice that, if $x_h \in \p\Omega_h$, then there is nothing to show. On the other hand, if $x_h \in \Omega_h$, then $\MAop{w}{h,S}{PD}(x_h)$ is an increasing function of the second differences $\Del_{\be} w_h(x_h)$. Indeed, increasing this difference makes the polytope larger. Notice also that $\Del_{\be} w_h(x_h)$ is a linear combination, with positive coefficients, of $w_h(x_h + \be h) - w_h(x_h)$ and $w_h(x_h + \be h) - w_h(x_h)$, with the obvious modification for points that are not interior with respect to $\be$. Thus, we can invoke \cite[Lemma 3.11]{NSWActa} to conclude the monotonicity.
\end{proof}

Next to be able to apply Theorem~\ref{thm:BarlesSouganidis} we must make sure that the operator $F_S$ satisfies a comparison principle. To establish this we begin with an auxiliary result.

\begin{lemma}[polytope comparison]
\label{lem:BMinnaplication}
Let $M_1,M_2 \in \polS_+^d$ and $x \in \Omega$. If $M_1 \leq M_2$ then, for every $t \in \R$ we have that $F_S(x,t,M_1) \leq F_S(x,t,M_2)$. In addition,
\begin{multline*}
  \left(F_S(x,t,M_1+M_2) + f(x) \right)^{1/d} \geq \\ \left( F_S(x,t,M_1)+ f(x) \right)^{1/d} + \left( F_S(x,t,M_2)+ f(x) \right)^{1/d}.
\end{multline*}
\end{lemma}
\begin{proof}
Notice that, since $x \in \Omega$ we have, independently of $t$,
\[
  F(x,t,M)+ f(x) = \left| K(M) \right|, \qquad K(M) = \left\{ \bv \in \R^d: \forall \be \in S,\ 2 \bv\cdot \be \leq \be \cdot M \be \right\}.
\]
Notice, in addition, that $M_1 \leq M_2$ implies that $\be \cdot M_1 \be \leq \be \cdot M_2 \be$ for every $\be \in \R^d$. Consequently, $M_1 \leq M_2$ implies $K(M_1) \subset K(M_2)$ from which the first statement follows.

Now, since $\be \cdot (M_1+M_2) \be = \be \cdot M_1 \be + \be \cdot M_2 \be$ we have that $K(M_1+M_2)$ contains $K(M_1) + K(M_2)$. The Brunn-Minkowski inequality given in Lemma~\ref{BM} allows us to conclude.
\end{proof}

Now we can establish a comparison principle for $F_S$.

\begin{proposition}[comparison]
\label{prop:comparePD}
Let $\overline{u} \in USC(\bar\Omega)$ and $\underline{u} \in LSC(\bar\Omega)$ be a sub and supersolution, respectively, of \eqref{eq:PDcontinuouscompact}. Then we have that $\overline{u} \leq \underline{u}$.
\end{proposition}
\begin{proof}
We begin by noticing that, since $F_{S,\star}(x,t,M) \leq F_S(x,t,M)$ for all $x \in \bar\Omega$ we obtain that, if $x_0 \in \partial\Omega$ we must have that
\[
  0 \leq F_{S,\star}(x_0,\varphi(x_0),D^2 \varphi(x_0) ) \leq F_S(x_0,\varphi(x_0),D^2 \varphi(x_0) ) = g(x_0) - \varphi(x_0),
\]
for every $\varphi$ sufficiently smooth that satisfies the conditions given in Definition~\ref{def:viscosol}. As a consequence, In this case, the condition defining a viscosity subsolution at boundary points reduces to $\overline{u} \leq g$ on $\p\Omega$. Similarly we can show that for a supersolution we must have $g \leq \underline{u}$ on $\p\Omega$. In conclusion, at the boundary $\p\Ome$ we have $\overline{u} \leq \underline{u}$.

By the semicontinuity assumption we can also define $\delta= \sup_{\bar\Omega}( \overline{u} - \underline{u}) \in \R$. Additionally, since $\Omega$ is bounded, there is $R>0$ such that $\Omega \subset B_R$. Assume now, for the sake of contradiction, that $\delta>0$.

Let us define, for $\eps>0$, the operator $F_{S,\eps}: \bar \Omega \times \R \times \polS^d \to \R$ by
\[
  F_{S,\eps}(x,t,M) = \begin{dcases}
                        F_S(x,t,M) - \eps (t - \underline{u}(x)), & x \in \Omega, \\
                        g(x) - t, & x \in \p\Omega,
                      \end{dcases}
\]
and notice that this operator satisfies all the conditions of the comparison principle given in \cite[Theorem 3.3]{CIL}. Moreover, since for all $x \in \bar\Omega$ we have that $F_{S,\eps}(x,\underline{u}(x),D^2 \underline{u}(x) ) = F_S(x,\underline{u}(x),D^2 \underline{u}(x) )$ we conclude that $\underline{u}$ is a supersolution for the operator $F_{S,\eps}$.

We now construct a subsolution. Define
\[
  v(x) = \frac{(\eps\delta)^{1/d}}2 \left( |x|^2 - R^2  \right), \qquad u_\eps(x) = \overline{u}(x) + v(x)
\]
and notice that $u_\eps \in USC(\bar\Omega)$ and, moreover, $u_\eps \leq \overline{u} \leq g \leq \underline{u}$ on $\p\Omega$. In addition, we have that, for $x \in \Omega$
\[
  D^2 v(x) = (\eps \delta)^{1/d} I, \qquad F_S(x,t,D^2v(x)) + f(x) = \eps\delta,
\]
see the proof of Proposition~\ref{prop:consistencyPowDiag}. Let now $x \in \Omega$ and, to shorten notation, denote
\[
    F_S[w] = F_S(x,w(x),D^2w(x))+f(x).
\]
If this is the case we have that, in the viscosity sense
\begin{align*}
  F_{S,\eps}(x,u_\eps(x),D^2 u_\eps(x) ) &= F_S[\overline{u}+v] - f(x) - \eps( \overline{u}(x) - \underline{u}(x) ) - \eps v(x) \\
    &\geq \left( F_S[\overline{u}]^{1/d} + F_S[v]^{1/d} \right)^d - f(x) - \eps( \overline{u}(x) - \underline{u}(x) ) \\
    &\geq \left( f(x)^{1/d} + F_S[v]^{1/d} \right)^d - f(x) - \eps( \overline{u}(x) - \underline{u}(x) )\\
    & \geq F_S[v]- \eps( \overline{u}(x) - \underline{u}(x) ) \\
    & = \eps \left( \delta - ( \overline{u}(x) - \underline{u}(x) ) \right) \geq 0.
\end{align*}
where we used Lemma~\ref{lem:BMinnaplication}, the fact that $v(x) \leq 0$ for all $x \in \bar\Omega$, that $\overline{u}$ is a subsolution for the operator $F_S$, the elementary identity
\[
  \left( x + y \right)^\theta \leq x^\theta + y^\theta, \quad \forall x,y \in \R_+, \quad \forall \theta \in (0,1],
\]
and the definition of $\delta$. In conclusion, $u_\eps$ is a subsolution  for the operator $F_{S,\eps}$. The comparison principle of \cite[Theorem 3.3]{CIL} then yields that
\[
  u_\eps(x) = \overline{u}(x) + v(x) \leq \underline{u}(x), \quad \forall x\in \bar\Omega
\]
or that
\[
  \underline{u}(x) - \overline{u}(x) \geq -\frac{(\eps\delta)^{1/d}}2 R^2, \quad \forall x \in \bar\Omega.
\]
Letting $\eps \downarrow 0$ we obtain $\overline{u}(x) \leq \underline{u}(x)$, contradicting that $\delta > 0$.

\end{proof}

As a consequence, we have convergence.

\begin{corollary}[convergence]
Let $\{u_h\}_{h>0} \subset X_h$ be the solutions to \eqref{eq:powdiag}. Then, as $h \downarrow 0$, we have that $u_h \to u$ locally uniformly, where $u$ is the (unique) viscosity solution of \eqref{eq:PDcontinuouscompact}.
\end{corollary}
\begin{proof}
Apply Theorem~\ref{thm:BarlesSouganidis}. It is only relevant to mention that owing to the comparison principle showed in Proposition~\ref{prop:comparePD}, $u$ must necessarily be unique.
\end{proof}


\subsection{Two scale methods}
\label{sub:twoscale}

We will now present and analyze the so--called {\it two scale method}, which can be understood as a generalization of the wide stencil schemes presented in Section~\ref{sub:Oberman} to unstructured meshes (see also \cite{Froese18MF}). Here and in what follows we will implicitly assume that $\Omega$ is uniformly convex. Additional assumptions will be explicitly stated. Next, for $h>0$, we introduce a quasiuniform (in the sense of \cite{CiarletBook}) simplicial triangulation $\mct$ of our domain $\Ome$. We denote by $\Nhi$ and $\Nhb$ the set of interior and boundary nodes, respectively, of $\mct$. We define $X_h$ to be the set of piecewise linear and continuous functions subject to this triangulation. The mesh size $h$ will constitute the fine scale of discretization. The large scale, denoted by $\delta$, will be the one at which second order differences will be evaluated. Notice that, since now we are dealing with continuous functions, these can be evaluated at any point. Indeed, given $x_h \in \Nhi$ and $\bw \in \R^d$ with $|\bw|=1$ we define, for $w_h \in X_h$
\begin{equation}
\label{eq:newseconddiff}
  \nabla^2_{\delta\bw} w_h(x_h) = \frac{ w_h(x_h+ \rho \delta \bw) - 2w_h(x_h) + w_h(x_h - \rho\delta\bw) }{\rho^2 \delta^2},
\end{equation}
where $\rho\in (0,1]$ is the largest number so that $x_h \pm \rho \delta \bw \in \bar\Ome$; compare with \eqref{eq:2nddiff} and \eqref{eq:2nddiffbdry}.
As a final discretization ingredient, as in the case of the wide stencil schemes of Section~\ref{sub:Oberman}, we need a directional discretization. That is if, as before, $\calV$ denotes the set of all orthonormal bases of $\R^d$ we must construct, for $\theta>0$, a set $\calV_\theta$ of collections of $d$ unit vectors such that if $\{\bw_i\}_{i=1}^d \in \calV$, then there is $\{\bw_i^\theta\}_{i=1}^d \in \calV_\theta$ such that
\begin{equation}
\label{eq:calVthetaisclose}
  \max_{i=1,\ldots,d} |\bw_i - \bw_i^\theta| \leq \theta.
\end{equation}
It is important to notice that the elements of $\calV_\theta$ are not required to be orthonormal collections of vectors.

Having defined all the discretization ingredients, which are parametrized by the triple $\eps = (h,\delta,\theta)$, following \cite{NochettoNtogkasZhangfirsttwoscale} we introduce the two scale discrete \MA operator by defining, for $w_h \in X_h$, and $x_h \in \Nhi$,
\begin{equation}
\label{eq:TSop}
  \begin{aligned}
  \MAop{w_h}{h,\delta,\theta}{2S}(x_h) = \min_{\{\bw_i\}_{i=1}^d \in \calV_\theta} &\left[\prod_{i=1}^d\left(\nabla^2_{\delta\bw_i} w_h(x_h) \right)^+ \right. \\ &-\left. \sum_{i=1}^d \left( \nabla^2_{\delta\bw_i} w_h(x_h) \right)^- \right],
  \end{aligned}
\end{equation}
compare with the scheme discussed in Remark~\ref{rem:variantWS}. With these ingredients at hand, the two scale method seeks a function $u_h^\eps \in X_h$ such that
\begin{subequations}
\label{eq:twoscale}
  \begin{align}
    \label{eq:twoscalein}
    \MAop{u_h^\eps}{h,\delta,\theta}{2S}(x_h) &= f(x_h), \quad \forall x_h \in \Nhi, \\
  \label{eq:towscalebc}
    u_h^\eps(x_h) &= g(x_h), \quad \forall x_h \in \Nhb.
  \end{align}
\end{subequations}

\begin{remark}[generalization]
\label{rem:WSis2S}
Starting from the Cartesian mesh $\Omega_h$ used to define the wide stencil schemes \eqref{eq:WSFD} it is possible to construct a simplicial triangulation of $\Omega$ without introducing new vertices: in two dimensions this is accomplished by subdividing each square by its diagonal, and a similar construction is possible in three dimensions. Once this is done, it can be seen that scheme \eqref{eq:twoscale} is, after little modifications, a generalization of the wide stencil scheme \eqref{eq:WSFD}.
\eremk\end{remark}

\begin{remark}[domain approximation]
\label{rem:Omegah}
Notice that, since the domain $\Omega$ is assumed to be uniformly convex, it is not possible to triangulate it exactly. If we 
denote $\bar\Omega_{\mct} = \cup_{T \in \mct} \bar T$, then we have $\bar\Omega_{\mct} \subsetneq \bar \Omega$. In our discussion we will ignore this fact. This is because we can either replace $\Omega$ by $\Omega_{\mct}$ in all the statements that we shall make, or we can consider all functions in $X_h$ as defined in $\Omega$ by extending them to $\Omega\setminus\Omega_{\mct}$ by a constant in the normal direction to faces. This is a standard construction and we shall not delve into it further.
\eremk\end{remark}

Let us now provide, following \cite{NochettoNtogkasZhangfirsttwoscale,NochettoNtogkasZhang,LiNochettoZhangMA}, an analysis of \eqref{eq:twoscale}. We will first introduce a discrete notion of convexity, based on the positivity of the second differences defined in \eqref{eq:newseconddiff}. The operator \eqref{eq:TSop} turns out to have a comparison principle, and acts in a particular way on discretely convex functions. This will allow us to establish existence, uniqueness, and stability of solutions to \eqref{eq:twoscale}. In addition, since the size large scale $\delta$ is reduced near the boundary, the consistency can only hold sufficiently far away from it. For this reason, appropriate barrier functions need to be constructed. All these ingredients will allow us to assert convergence of the method. Finally, using the comparison principle and suitable barriers, we will establish rates of convergence for classical solutions.

\subsubsection{Discrete convexity}
\label{subsub:dconvexity}

The second order differences defined in \eqref{eq:newseconddiff} and the set of directions $\calV_\theta$ give a discrete notion of convexity.

\begin{definition}[discrete convexity]
\label{def:discrconvex}
We say that the function $w_h \in X_h$ is discretely convex if
\[
  \nabla^2_{\delta\bw_j} w_h(x_h) \geq 0, \quad \forall x_h \in \Nhi, \quad\forall \bw_j \in \{\bw_i\}_{i=1}^d \in \calV_\theta.
\]
\end{definition}

It is well known that if a function is convex, then its second order differences are nonnegative. On the other hand, discrete convexity does not imply convexity. This is due, for instance, to the fact that convexity and interpolation are not easily compatible. In other words, if $w \in C(\bar\Ome)$ is convex, then its Lagrange interpolant $\calI_h w \in X_h$ satisfies $\calI_h w \geq w$ so that it is discretely convex, but $\calI_h w$ is not necessarily convex.

On the other hand, discrete convexity implies nonnegativity of the two scale discrete \MA operator; see \cite[Lemma 2.2]{NochettoNtogkasZhangfirsttwoscale}.

\begin{lemma}[discrete convexity]
\label{lem:discrconvex}
A function $w_h \in X_h$ is discretely convex if and only if
\[
  \MAop{w_h}{h,\delta,\theta}{2S}(x_h) \geq 0, \quad \forall x_h \in \Nhi.
\]
Moreover, for a discretely convex function we have that
\[
  \MAop{w_h}{h,\delta,\theta}{2S}(x_h) = \min_{\{\bw_i\}_{i=1}^d \in \calV_\theta} \prod_{i=1}^d \nabla^2_{\delta\bw_i} w_h(x_h).
\]
\end{lemma}

\subsubsection{A comparison principle}
\label{subsub:2smonotone}

Let us now show that the operator defined in \eqref{eq:TSop} is monotone and has a comparison principle. From this we will obtain uniqueness of solutions to \eqref{eq:twoscale}.

\begin{lemma}[monotonicity]
\label{lem:twoscalemonotone}
Let $v_h,w_h \in X_h$ be such that $v_h - w_h$ attains its maximum at the interior node $x_h \in \Nhi$. Then we have
\[
  \MAop{w_h}{h,\delta,\theta}{2S}(x_h) \geq \MAop{v_h}{h,\delta,\theta}{2S}(x_h).
\]
\end{lemma}
\begin{proof}
Since $x_h$ is the maximum, for suitable $\rho>0$ and any unit vector $\bw$ we have
\[
  v_h(x_h) - w_h(x_h) \geq v_h(x_h \pm \rho \delta \bw) - w_h(x_h \pm \rho \delta \bw),
\]
which implies that
\[
  \nabla^2_{\delta \bw} v_h(x_h) \leq \nabla^2_{\delta \bw} w_h(x_h).
\]
multiplying this inequality as $\bw$ runs over all elements of $\calV_\theta$ allows us to conclude.
\end{proof}

The previous result gives us a comparison principle for the operator \eqref{eq:TSop}.

\begin{proposition}[comparison]
\label{prop:tscompare}
Let $v_h,w_h \in X_h$ be such that $v_h \leq w_h$ on $\p\Ome$, and
\[
  \MAop{v_h}{h,\delta,\theta}{2S}(x_h) \geq \MAop{w_h}{h,\delta,\theta}{2S}(x_h), \quad \forall x_h \in \Nhi,
\]
then we must have that $v_h \leq w_h$ in $\bar\Ome$.
\end{proposition}
\begin{proof}
We consider two cases for the inequality between the operators:
\begin{enumerate}
  \item {\it The inequality is strict}. Let us assume, for the sake of contradiction, $v_h - w_h$ attains a maximum at an interior node. Lemma~\ref{lem:twoscalemonotone} then gives a contradiction.
  
  \item {\it The inequality is not strict}. Since $\Omega$ is bounded, there is $R>0$ such that the convex quadratic $q(x) = \tfrac12( |x|^2 -R)$ is nonpositive on $\bar\Omega$. Let $q_h = \calI_h q \in X_h$. This function is strictly convex and satisfies
  \[
    \nabla^2_{\delta \bw} q_h (x_h) \geq \nabla^2_{\delta \bw} q (x_h) = \frac{\partial^2 q(x_h)}{\partial \bw^2} = 1.
  \]
  We claim now that, for all $\alpha >0$ and $x_h \in \Nhi$, we have that
  \begin{equation}
  \label{eq:newforcompare}
    \MAop{v_h + \alpha q_h}{h,\delta,\theta}{2S}(x_h) \geq \MAop{v_h}{h,\delta,\theta}{2S}(x_h)  + \min \left\{ \frac{\alpha^d}{2^d} + \frac\alpha2 \right\}.
  \end{equation}
  Indeed, fix $\{\bw_i\} \in \calV_\theta$ and assume first that $\nabla^2_{\delta \bw_i}(v_h(x_h) + \tfrac\alpha2 q_h(x_h) ) \geq 0$ for all $i$. In this case
  \begin{align*}
    \prod_{i=1}^d \Big( \nabla^2_{\delta \bw_i}(v_h(x_h) + &\alpha q_h(x_h) \Big) 
    \geq    \prod_{i=1}^d \left( \nabla^2_{\delta \bw_i}(v_h(x_h) + \tfrac\alpha2 q_h(x_h) ) + \frac\alpha2\right) \\
    &\geq \min_{\{\bw_i\} \in \calV_\theta} \prod_{i=1}^d \left( \nabla^2_{\delta \bw_i}(v_h(x_h) + \tfrac\alpha2 q_h(x_h) ) \right) + \frac{\alpha^d}{2^d} \\
    &\geq \left( \prod_{i=1}^d \left( \nabla^2_{\delta \bw_i}v_h(x_h) \right)^+ - \sum_{i=1}^d \left( \nabla^2_{\delta \bw_i}v_h(x_h) \right)^- \right) + \frac{\alpha^d}{2^d}.
  \end{align*}
  On the other hand, if there is $i \in \{1,\ldots,d\}$ for which $\nabla^2_{\delta \bw_i}(v_h(x_h) + \alpha q_h(x_h) ) < 0$, then this implies that $\nabla^2_{\delta \bw_i}v_h(x_h) < 0$. Thus,
  \[
    \prod_{i=1}^d \left( \nabla^2_{\delta \bw_i}v_h(x_h) \right)^+ = 0,
  \]
  and
  \begin{align*}
  -\sum_{i=1}^d \Big(\nabla^2_{\delta \bw_i}(v_h(x_h) + &\alpha q_h(x_h) ) \Big)^- \geq 
  -\sum_{i=1}^d \left(\nabla^2_{\delta \bw_i}v_h(x_h)  \right)^- + \frac\alpha2 \\
  &=
  \left( \prod_{i=1}^d \left( \nabla^2_{\delta \bw_i}v_h(x_h)\right)^+ -\sum_{i=1}^d \left(\nabla^2_{\delta \bw_i}v_h(x_h) \right)^- \right) + \frac\alpha2.
  \end{align*}

  A combination of these two cases, since $\{\bw_i\}_{i=1}^d \in \calV_\theta$ was arbitrary, implies \eqref{eq:newforcompare}.
  
  Finally, since, $v_h + \alpha q_h \leq v_h \leq w_h$ on $\p\Ome$ and, on the basis of \eqref{eq:newforcompare}, we have
  \[ 
    \MAop{v_h + \alpha q_h}{h,\delta,\theta}{2S}(x_h) > \MAop{v_h}{h,\delta,\theta}{2S}(x_h) \geq \MAop{w_h}{h,\delta,\theta}{2S}(x_h), \ \ \forall x_h \in \Nhi,
  \]
  the previous step then implies that $v_h + \alpha q_h \leq w_h$. Letting $\alpha \downarrow 0$ we can conclude.
\end{enumerate}
\end{proof}

\begin{remark}[discrete interior barrier]
\label{rem:interiorbarrier}
Notice, that, in the course of the second case of the proof of this result we effectively constructed a discrete interior barrier. If $q(x) = \tfrac12(|x|^2-R)$ with $R>0$ sufficiently large, then we have that
\[
  \calI_h q \leq 0, \ \text{ on } \p \Omega, \qquad \MAop{\calI_h q_h}{h,\delta,\theta}{2S}(x_h) \geq 1, \ \forall x_h \in \Nhi.
\]
\eremk\end{remark}

As an immediate consequence, we also have uniqueness of solutions to \eqref{eq:twoscale}.

\begin{corollary}[uniqueness]
\label{cor:tsunique}
Scheme \eqref{eq:twoscale} cannot have more than one solution.
\end{corollary}

As a final application of the comparison principle, let us now show existence and uniform bounds on the solution to \eqref{eq:twoscale}.

\begin{theorem}[existence and stability]
For all $\eps = (h,\delta,\theta)>0$ scheme \eqref{eq:twoscale} has a solution $u_h^\eps \in X_h$. Moreover, this solution is stable in the sense that $\| u_h^\eps \|_{L^\infty(\Omega)}$ is bounded independently of $\eps$.
\end{theorem}
\begin{proof}
The existence proceeds via Perron's method. For this reason, we will only indicate how to construct a discrete subsolution, that is a function $u_h^0 \in X_h$ such that $u_h^0 = \calI_h g$ on $\partial\Omega$ and
\[
  \MAop{u_h^0}{h,\delta,\theta}{2S}(x_h) \geq f(x_h), \quad \forall x_h \in \Nhi.
\]
To construct this function, we define
\[
  s(x) = \sum_{x_h \in \Nhi} s_{x_h}(x), \qquad s_{x_h}(x) = \frac{\delta \rho_{x_h}}2 f(x_h)^{1/d} |x-x_h|,
\]
where $\rho_{x_h} \in (0,1]$ is the largest number such that, for all $\bw \in \R^d$ with $|\bw|=1$ we have $x_h\pm \rho_{x_h} \bw \in \bar \Omega$. Notice that $\nabla^2_{\delta\bw}s_{x_h}(y_h) \geq 0$ for all $y_h \in \Nhi$, and that
\[
  \nabla^2_{\delta\bw} s_{x_h}(x_h) = f(x_h)^{1/d}, \quad \forall \bw \in \R^d, \ |\bw|=1.
\]
Consequently, for $y_h \in\Nhi$
\[
  \nabla^2_{\delta\bw} \calI_h s(y_h) \geq \nabla^2_{\delta\bw}  s(y_h) \geq f(y_h)^{1/d} \geq 0,
\]
which, by Lemma~\ref{lem:discrconvex} implies
\[ 
  \MAop{\calI_h s}{h,\delta,\theta}{2S}(x_h) = \min_{\{\bw_i\}_{i=1}^d \in \calV_\theta} \prod_{i=1}^d \nabla^2_{\delta\bw_i} \calI_h s(x_h) \geq f(x_h), \quad \forall x_h \in \Nhi.
\]

Let now $w \in C(\bar\Omega)$ be the convex envelope of $(\calI_h(g-s))|_{\p\Omega}$, and set $w_h = \calI_h w$. By convexity of $w$ we have that
\[
  \MAop{w_h}{h,\delta,\theta}{2S}(x_h) \geq 0, \quad \forall x_h \in \Nhi.
\]
Thus, we define
\[
  u_h^0= w_h + \calI_h s.
\]
This function, by construction, is discretely convex and $u_h^0 = \calI_h g$ on $\p\Omega$. Since the second differences of $w_h$ are nonnegative, then we have that
\begin{align*}
  \MAop{u_h^0}{h,\delta,\theta}{2S}(x_h) 
  &= 
  \min_{\{\bw_i\}_{i=1}^d \in \calV_\theta} \prod_{i=1}^d \left[ \nabla^2_{\delta\bw_i} w_h (x_h) + \nabla^2_{\delta\bw_i} \calI_h s(x_h) \right] \\
  &\geq \min_{\{\bw_i\}_{i=1}^d \in \calV_\theta} \prod_{i=1}^d \nabla^2_{\delta\bw_i} \calI_h s(x_h)
  \geq f(x_h),
\end{align*}
and so $u_h^0$ is a discrete subsolution.

It remains to show the uniform boundedness. To achieve this we will show that every discrete subsolution is uniformly bounded. Let then $w_h \in X_h$ be a discrete subsolution and $b_h = \max_{x \in \p \Omega} g(x) \in X_h$. We have then that
\[
  \MAop{b_h}{h,\delta,\theta}{2S}(x_h) = 0 \leq f(x_h) \leq \MAop{w_h}{h,\delta,\theta}{2S}(x_h), \quad \forall x_h \in \Nhi.
\]
Since, in addition, we have that $b_h \geq w_h$ on $\p\Omega$, the comparison principle of Proposition~\ref{prop:tscompare} implies that
\[
  w_h \leq b_h.
\]
This is enough since Perron's method shows existence of a solution by constructing an increasing sequence of subsolutions. Thus, $u_h^0$ is a lower bound for the solution and, evidently, $\|u_h^0\|_{L^\infty(\Omega)}$ is independent of $\eps$.
\end{proof}

\subsubsection{Consistency and discrete barriers}
\label{subsub:twoscaleconsistent}

Let us now examine the consistency of the operator \eqref{eq:TSop}. As we have stated above, the operator can only be consistent at points sufficiently far away from the boundary. For this reason, we define the $\delta$-interior and $\delta$-boundary layer of $\Ome$ via
\[
  \Omega_\delta = \bigcup_{T \in \mct: \dist(T,\p\Omega)> \delta} T,
  \qquad
  (\p \Omega)_\delta = \bar\Omega \setminus \Omega_\delta.
\]
For an interior node $x_h \in \Nhi$ its interior patch is
\[
  \omega_{x_h} = \bigcup_{T \in \mct: \dist(x_h,T) < \rho \delta } \bar T,
\]
where, as before, $\rho \in (0,1]$ is the largest number such that, for any $\bw \in \R^d$ with $|\bw|=1$ we have $x_h \pm \rho \delta \bw \in \bar \Omega$.

The following result follows, essentially, by a Taylor expansion argument.

\begin{lemma}[consistency of second differences]
\label{lem:2nddiffconsistent}
Let $x_h \in \Nhi$ and assume that $\varphi \in C^{1,1}(\omega_{x_h})$, then for all $\bw \in \R^d$ with $|\bw|=1$ we have
\[
  |\nabla^2_{\delta\bw} \calI_h \varphi(x_h)| \leq C |\varphi|_{C^{1,1}(\omega_{x_h})}.
\]
If, in addition, we have that $x_h \in \Omega_\delta$ and $\varphi \in C^{k+2,\alpha}(\omega_{x_h})$ for $k = 0,1$, and $\alpha \in (0,1]$ then we also have that
\[ 
  \left| \nabla^2_{\delta\bw} \calI_h \varphi(x_h) - \frac{\p^2 \varphi(x_h)}{\p \bw^2} \right| \leq C
    \left(  |\varphi|_{C^{k+2,\alpha}(\omega_{x_h})} \delta^{k+\alpha} + \frac{h^2}{\delta^2} |\varphi|_{C^{1,1}(\omega_{x_h})} \right).
\]
Finally, if $\varphi$ is, in addition, convex then we have
\[ 
  \frac{\p^2 \varphi(x_h)}{\p \bw^2} - \nabla^2_{\delta\bw} \calI_h \varphi(x_h) \leq C
    |\varphi|_{C^{k+2,\alpha}(\omega_{x_h})} \delta^{k+\alpha} .
\]
\end{lemma}

The previous result can be applied to obtain interior consistency of \eqref{eq:TSop}. The following result was first obtained in \cite[Lemma 4.2]{NochettoNtogkasZhangfirsttwoscale} under the assumption that $\calV_\theta \subset \calV$. This assumption was later removed in \cite[Lemma 2.4]{LiNochettoZhangMA}.

\begin{theorem}[interior consistency]
\label{thm:twoscaleconsistencyin}
Let $x_h \in \Nhi$ and $\varphi \in C^{k+2,\alpha}(\omega_{x_h})$ with $k=0,1$ and $\alpha \in (0,1]$ be convex. In this setting the following estimates are valid:
\begin{enumerate}
  \item \label{it:detisclose} If $\theta \leq \tfrac1{4d}$ then, for any $\{\bw_i\}_{i=1}^d \in \calV_\theta$, we have
  \[
    \det D^2 \varphi(x_h) \leq \prod_{i=1}^d \frac{ \p^2 \varphi(x_h)}{\p \bw_i^2}\left( 1 + 16 \theta^2(d-1)^2 \right).
  \]
  
  \item \label{it:detisclose2} If $\{\bv_i\}_{i=1}^d \in \calV$ realizes the minimum in the variational characterization of the determinant given in Lemma~\ref{lem:detisvar}, then for any $\{\bv_i^\theta \}_{i=1}^d \in \calV_\theta$ that satisfies \eqref{eq:calVthetaisclose} we have
  \[
    \left| \frac{\p^2 \varphi(x_h)}{\p \bv_i^2} - \frac{\p^2 \varphi(x_h)}{\p {\bv_i^\theta}^2} \right| \leq C |\varphi|_{C^{1,1}(\omega_{x_h})} \theta^2.
  \]
  
  \item Finally, if $x_h \in \Nhi \cap \Omega_\delta$, then
  \[
    \left| \det D^2 \varphi(x_h) - \MAop{\calI_h \varphi}{h,\delta,\theta}{2S}(x_h) \right| \leq C_1 \delta^{k+\alpha} + C_2 \left( \frac{h^2}{\delta^2} + \theta^2 \right),
  \]
  where the constants $C_1$ and $C_2$ depend only on the smoothness of $\varphi$, the domain $\Omega$, and the dimension $d$.
\end{enumerate}
\end{theorem}
\begin{proof}
We prove each statement separately.
\begin{enumerate}
  \item Let $W_\theta = (\bw_1, \ldots, \bw_d)$. We have
  \[
    \det(W_\theta^\intercal W_\theta) \det D^2\varphi(x_h) = \det (W_\theta^\intercal D^2 \varphi(x_h) W_\theta) \leq 
    \prod_{i=1}^d \frac{ \p^2 \varphi(x_h)}{\p \bw_i^2},
  \]
  where, in the last step, we used that $W_\theta^\intercal D^2 \varphi(x_h) W_\theta$ is positive semidefinite and Hadamard's inequality. We now need to estimate the determinant of $W = W_\theta^\intercal W_\theta$ from below. Write
  \[
    W = \begin{pmatrix}
          W_0 & \bw \\ \bw^\intercal & 1
        \end{pmatrix}
    = \begin{pmatrix}
        I & 0 \\ \bw^\intercal W_0^{-1} & 1
      \end{pmatrix}
      \begin{pmatrix}
        W_0 & \bw \\ 0 & 1- \bw\cdot W_0 \bw
      \end{pmatrix}
  \]
  implying that $ \det W = (1- \bw\cdot W_0 \bw) \det W_0$,
  which holds if the submatrix $W_0$ is nonsingular. Notice, however, that $W_{i,i} = 1$ 
   and $|W_{i,j}| \leq 2 \theta$ as the columns of $W_\theta$ form an element of $\calV_\theta$. This implies, for $\theta \leq \tfrac1{4d}$, that $W_0 \geq \tfrac12 I$ and $|\bw| \leq 2\theta \sqrt{d-1}$. Thus, $W_0^{-1} \geq 2 I$ and
  \[
    |\bw \cdot W_0 \bw | \leq 8 \theta^2(d-1) \qquad \det W \geq (1-8 \theta^2(d-1)) \det W_0,
  \]
  which by repeating this process yields
  \[
    \det W \geq (1-8 \theta^2(d-1))^{d-1} \geq 1-8 \theta^2(d-1)^2,
  \]
  and using, again the bound on $\theta$
  \[
    \frac1{\det W} \leq 1+ 16 \theta^2(d-1)^2.
  \]
  
  \item We begin by noticing that, given the minimization assumption, $\{\bv_i\}_{i=1}^d$ must be the normalized eigenvectors of $D^2 \varphi(x_h)$. Set $\bv_i^\theta = \bv_i + \bw_i$ and write
  \[
    \frac{ \p^2 \varphi(x_h)}{\p {\bv_i^\theta}^2} = \bv_i^\theta \cdot D^2 \varphi(x_h) \bv_i^\theta = 
    \frac{ \p^2 \varphi(x_h)}{\p \bv_i^2} + 2 \bw_i\cdot D^2\varphi(x_h) \bv_i + \bw_i\cdot D^2\varphi(x_h) \bw_i.
  \]
  Since $\{\bv_i\}_{i=1}^d$ are eigenvectors, $|\bw_i|\leq \theta$, and $|\bv_i \cdot \bw_i | \leq \tfrac12 \theta^2$ we then have
  \[
    \left| \frac{ \p^2 \varphi(x_h)}{\p {\bv_i^\theta}^2} -\frac{ \p^2 \varphi(x_h)}{\p \bv_i^2} \right| \leq C \theta^2.
  \]

  \item By Lemma~\ref{lem:discrconvex}, since $\calI_h \varphi$ is discretely convex, we have that
  \[
    \MAop{\calI_h \varphi}{h,\delta,\theta}{2S}(x_h) = \min_{\{\bw_i\}_{i=1}^d \in \calV_\theta} \prod_{i=1}^d \nabla^2_{\delta\bw_i} \calI_h \varphi(x_h).
  \]
  Let $\{\bw_i\}_{i=1}^d \in \calV_\theta$ be the set that realizes the minimum in this expression. Using Lemma~\ref{lem:detisvar} we can write that
  \begin{align*}
    \det D^2 \varphi(x_h) - \MAop{\calI_h \varphi}{h,\delta,\theta}{2S}(x_h) &\leq 
    \prod_{i=1}^d \frac{ \p^2 \varphi(x_h)}{\p \bw_i^2} - \prod_{i=1}^d \nabla^2_{\delta\bw_i} \calI_h \varphi(x_h) \\ &\leq C \delta^{k+\alpha},
  \end{align*}
  where, in the last step, we used repeatedly Lemma~\ref{lem:2nddiffconsistent}.

  Let now $\{\bv_i\}_{i=1}^d \in \calV$ be the normalized eigenpairs of $D^2\varphi(x_h)$, and $\{\bv_i^\theta\}_{i=1}^d \in \calV_\theta$ the collection that realizes \eqref{eq:calVthetaisclose}. Then we have
\begin{align*}
    &\MAop{\calI_h \varphi}{h,\delta,\theta}{2S}(x_h) - \det D^2 \varphi(x_h)\\
    &\quad \leq 
    \left( \prod_{i=1}^d \nabla^2_{\delta\bv^\theta_i} \calI_h \varphi(x_h) - \prod_{i=1}^d \frac{ \p^2 \varphi(x_h)}{\p {\bv_i^\theta}^2} \right) +
    \left( \prod_{i=1}^d \frac{ \p^2 \varphi(x_h)}{\p {\bv_i^\theta}^2} - \prod_{i=1}^d \frac{ \p^2 \varphi(x_h)}{\p {\bv_i}^2}\right).
  \end{align*}
  The first term can be handled by repeatedly applying Lemma~\ref{lem:2nddiffconsistent}, while the second by applying the previous step.
\end{enumerate}
All the estimates have been proved and the interior consistency is thus obtained.
\end{proof}

As mentioned before, the operator is not consistent near the boundary. For this reason we will, instead, construct discrete barriers which will allow us to control the behavior of the solution near the boundary.

\begin{proposition}[discrete barrier I]
\label{prop:discrbdrybarrier}
Let $E>0$ be arbitrary and $x_h \in \Nhi$ be such that $\dist(x_h,\p\Omega) \leq \delta$. Then, there is $p_h \in X_h$ such that 
\[
  p_h \leq 0, \text{ on } \p\Omega, \quad \MAop{p_h}{h,\delta,\theta}{2S}(y_h) \geq E, \ \forall y_h \in \Nhi, \quad |p_h(x_h)| \leq C E^{1/d}\delta,
\]
where the constant $C$ depends only on the domain $\Omega$.
\end{proposition}
\begin{proof}
Without loss of generality, we can assume that $x_h = (0,\ldots,0,z)^\intercal$ with $z>0$ so that $0 \in \partial \Omega$ and $z = \dist(x_h,\p\Omega)$. The uniform convexity of $\Omega$ shows that there is $R>0$ such that, in this system of coordinates,
\[
  \Omega \subset \left\{ x \in \R^d: \sum_{i=1}^{d-1} x_i^2 - (x_d-R)^2 \leq R^2 \right\}.
\]
Let
\[
  p(x) = \frac{E^{1/d}}2 \left( \sum_{i=1}^{d-1} x_i^2 - (x_d-R)^2 - R^2 \right).
\]
We claim that $p_h = \calI_h p$ is the desired barrier. Indeed, by construction $p_h \leq 0$ on the boundary $\p\Omega$ and, since $z \leq \delta$ we have that $|p_h(x_h)| \leq C E^{1/d}\delta$. Finally, since $p_h$ is discretely convex, for any interior node $y_h$ we have
\[
  \MAop{p_h}{h,\delta,\theta}{2S}(y_h) \geq \MAop{p}{h,\delta,\theta}{2S}(y_h) = \prod_{i=1}^d E^{1/d} = E,
\]
as claimed.
\end{proof}

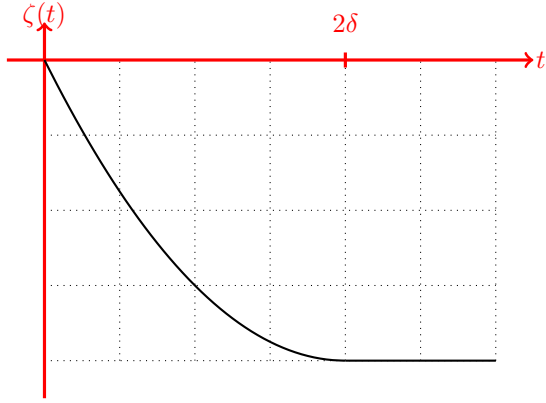
\begin{figure}
\begin{center}
  \begin{tikzpicture}
    \draw[step = 1, dotted] (1,1) grid (7,5);
    \draw[red, very thick, ->] (0.5,5) -- (7.5,5);
      \node[red] at (7.6,5) {$t$};
    \draw[red, very thick] (5,4.9) -- (5,5.1);
      \node[red] at (5,5.5) {$2\delta$};
    \draw[red, very thick, ->] (1,0.5) -- (1,5.5);
      \node[red] at (1,5.6) {$\zeta(t)$};
    \draw[thick, domain=1:5, smooth] plot (\x, {((\x-5)*(\x-5)/4.) + 1. } ) ;
    \draw[thick] (5,1) -- (7,1);
  \end{tikzpicture}
\end{center}
\caption{The function $\zeta$ used to define the discrete barrier of Proposition~\ref{prop:otherbarrier}.}
\label{fig:functionzeta}
\end{figure}

To obtain rates of convergence we shall also require another discrete barrier that was originally introduced in \cite[Section 6.2]{NochettoZhang16}. We define
\[
  \zeta : [0,\infty) \to (-\infty, 0], \qquad \zeta(t) = \begin{dcases}
                                                           (t-2\delta)^2 - (2\delta)^2, & t \in [0,2\delta], \\
                                                           -(2\delta)^2, & t \in (2\delta, \infty).
                                                         \end{dcases}
\]
The graph of this function is illustrated in Figure~\ref{fig:functionzeta}. With this function at hand, we define
\[
  b(x) = \zeta(\dist(x,\p\Omega)),
\]
and $b_h = \calI_h b$. The properties of this barrier are as follows.

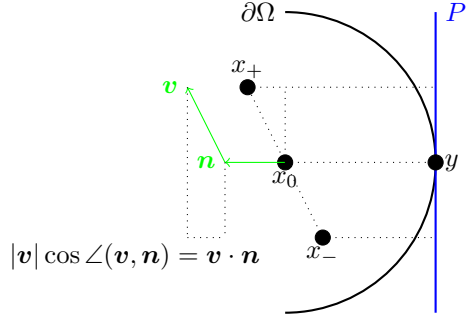
\begin{figure}
\begin{center}
  \begin{tikzpicture}
    \draw[thick] (0,0) arc[radius=2, start angle = -90, end angle = 90] node[left] {$\p\Omega$};
    \draw[blue, thick] (2,0)--(2,4) node[right] {$P$};
    
    \draw[fill] (2,2) circle[radius=0.1] node[right] {$y$};
    \draw[fill] (0,2) circle[radius=0.1] node[below] {$x_0$};
    \draw[green,->] (0,2) -- (-.8,2) node[left] {$\bn$};
    
    \draw[fill] (-0.5,3) circle[radius=0.1] node[above] {$x_+$};
    \draw[fill] (.5,1) circle[radius=0.1] node[below] {$x_-$};
    \draw[dotted] (.5,1) -- (-0.5,3);
    \draw[dotted] (0,2) -- (0,3);
  
    \draw[dotted] (-0.5,3) -- (2,3);
    \draw[dotted] (0,2) -- (2,2);
    \draw[dotted] (0.5,1) --(2,1);
    
    \draw[green, ->] (-0.8,2) -- (-1.3,3) node[left] {$\bv$};

    \draw[dotted] (-1.3,3) -- (-1.3,1);
    \draw[dotted] (-0.8,2) -- (-0.8,1);
    \draw[dotted] (-1.3,1) -- (-0.8,1);

    \node at (-2,0.75) {$|\bv|\cos \angle(\bv,\bn) = \bv \cdot \bn$};
    
  \end{tikzpicture}
\end{center}
\caption{The construction Proposition~\ref{prop:otherbarrier} that shows that the function $b$ is convex. The distance between $x_+$ and the supporting hyperplane $P$ equals the sum of the distance from $x_0$ to the boundary $\p\Omega$ and the inner product between $\bn$ and $\bv$.}
\label{fig:proofbarrier}
\end{figure}

\begin{proposition}[discrete barrier II]
\label{prop:otherbarrier}
For $\theta \leq \tfrac1{2\sqrt{d}}$ the barrier function $b_h$ satisfies:
\begin{enumerate}
  \item For all $x_h \in \Nhi$ and any $\bw \in \R^d$ with $|\bw|=1$,
  \[
    \nabla^2_{\delta \bw}b_h(x_h) \geq 0.
  \]

  \item \label{it:otherbarrierbdry} For all $x_h \in \Nhi \setminus \Omega_\delta$ and $\{\bw_i^\theta\}_{i=1}^d \in \calV_\theta$,
  \[
    \max_{i=1,\ldots,d} \nabla^2_{\delta \bw_i^\theta}b_h(x_h) \geq \frac1{2d}.
  \]

  \item For all $x \in \bar\Omega$
  \[
    -4\delta^2 \leq b_h(x) \leq 0.
  \]
\end{enumerate}
\end{proposition}
\begin{proof}
We consider each property separately.
\begin{enumerate}
  \item {Let $x_+,x_-\in \Omega$ with $x_+\neq x_-$.  The convexity of $\Omega$ ensures that $x_0 = \frac12(x_++x_-)\in \Omega$.
Denote by $y \in \partial\Omega$ the closest point to $x_0$. Since $\Omega$ is convex, there is a supporting hyperplane $P$ at $y$, whose normal is $\bn = \tfrac1{|x_0 - y|}(x_0 -y)$. Let now $\bv = \pm(x_+ - x_-)$, where the sign is chosen so that
$\bn\cdot \bv\ge 0$.}
 Consequently, see Figure~\ref{fig:proofbarrier},
  \[
    \dist(x_\pm,\p\Omega) \leq \dist(x_\pm, P) = \dist(x_0,\p\Omega) \pm \bn \cdot \bv.
  \]
  With this estimate, and using that $\zeta$ is nonincreasing, we can compute
  \begin{align*}
    b(x_+) + b(x_-) &\geq \zeta\left( \dist(x_0,\p\Omega) + \bv\cdot\bn \right) + \zeta\left( \dist(x_0,\p\Omega) - \bv\cdot\bn \right)\\
    & \geq 2\zeta\left( \dist(x_0,\p\Omega)  \right) = 2b(x_0),
  \end{align*}
where the second inequality follows directly from the definition of $\zeta$.
We then conclude (cf.~\cite[Pages 1--2]{MR0126722}) that the function $b$ is convex and the stated property of $b_h$ follows. 
  
  \item With the notation of the previous step, if we take a node $x_h \in \Nhi \setminus \Omega_\delta$, and $\bv \in \R^d$ with $|\bv| = 1$, then $\dist(x_h,\p\Omega) \pm \rho \delta \bw \cdot \bn \in [0,2\delta]$. Since $\zeta$ is nonincreasing and quadratic on that interval
  \[
    \nabla^2_{\delta \bv} b_h(x_h) \geq \nabla^2_{\delta\bv} b(x_h) \geq 2 \frac{ \rho^2 \delta^2 |\bv \cdot \bn|^2}{\rho^2\delta^2} = 2 |\bv \cdot \bn|^2.
  \]
  Now, if we let $\bv$ run over $\{\bw_i^\theta\}_{i=1}^d \in \calV_\theta$ we have obtained that
  \[
    \max_{i=1,\ldots,d} \nabla^2_{\delta \bw_i^\theta} b_h(x_h) \geq 2 \max_{i=1,\ldots,d} |\bw_i^\theta \cdot \bn|^2.
  \]
  Let now $\{\bw_i\}_{i=1}^d \in \calV$ be such that it satisfies \eqref{eq:calVthetaisclose}. Since $|\bn|=1$ we must have that
  \[
    \sum_{i=1}^d |\bn \cdot \bw_i|^2 = 1, \quad \implies \max_{i=1,\ldots,d} |\bn \cdot \bw_i| \geq \frac1{\sqrt{d}}.
  \]
  Therefore,
  \[
    |\bw_i^\theta \cdot \bn| \geq |\bw_i\cdot \bn| - |(\bw_i - \bw_i^\theta)\cdot \bn| \geq |\bw_i\cdot \bn| - \theta \geq \frac1{2\sqrt{d}},
  \]
  where we used that $\theta \leq \frac1{2\sqrt{d}}$. This implies the estimate.
  
  \item The last property follows directly from the definition of the function $\zeta$.

\end{enumerate}

\end{proof}

\subsubsection{Convergence}
\label{subsub:twoscaleconverge}

Let us now show convergence. We will do so by adapting the arguments developed in Section~\ref{sub:BarlesSouganidis} to take into account that test functions must be convex. We will rely on Proposition~\ref{prop:equiv}.

\begin{theorem}[convergence]
\label{thm:twoscaleconverge}
Let $\Omega$ be uniformly convex, $f \in C(\bar\Omega)$ such that $f \geq 0$, and $g \in C(\p\Omega)$. As $\eps = (h,\delta,\theta) \to 0$ with $h\delta^{-1} \to 0$ we have that the family $\{u_h^\eps\}_{\eps}$ of solutions of \eqref{eq:twoscale} converges uniformly to $u \in C(\bar\Omega)$, the solution of \eqref{eqn:MA}.
\end{theorem}
\begin{proof}
In a similar way to Theorem~\ref{thm:twoscaleconsistencyin} we have that, for all $x_0 \in \Omega$, $x_h \in \Nhi \cap \Omega_\delta$ and all $\varphi \in C^{2,\alpha}(\omega_{x_h})$, it holds that
\begin{equation}
\label{eq:tsnewcons}
  \begin{aligned}
    |\MAop{\varphi}{}{}(x_0) - \MAop{\calI_h \varphi}{h,\delta,\theta}{2S}(x_h)| &\leq C_1 (\delta^\alpha + |x_0-x_h|^\alpha) \\ &+ C_2 \left( \frac{h^2}{\delta^2} + \theta^2 \right).
  \end{aligned}
\end{equation}
Indeed, we only need to use that the operations $t \mapsto t^\pm$ are Lipschitz and with Lipschitz constant equal one.

We now extend the ideas of Theorem~\ref{thm:BarlesSouganidis}. As there we define
\[
  \overline{u}(x) = \limsup_{\eps \to 0, \frac{h}\delta \to 0, y \to x} u_h^\eps(y),
  \qquad
  \underline{u}(x) = \liminf_{\eps \to 0, \frac{h}\delta \to 0, y \to x} u_h^\eps(y)
\]
and we will show that $\overline{u}$ is a subsolution. For that we assume that $\overline{u}-\varphi$, with $\varphi \in C^{2,\alpha}(\bar\Omega)$ attains a maximum at $x_0 \in \Omega$. Let $\{x_h\}$ be the sequence of nodes such that $x_h \to x_0$ and $u_h^\eps - \calI_h\varphi$ attains a maximum at $x_h$. By the monotonicity result of Lemma~\ref{lem:twoscalemonotone} we obtain then that
\[
  \MAop{\calI_h \varphi}{h,\delta,\theta}{2S}(x_h) \geq \MAop{u_h^\eps}{h,\delta,\theta}{2S}(x_h) \geq f(x_h),
\]
the consistency, as expressed in \eqref{eq:tsnewcons}, implies by passing to the limit that
\[
  \MAop{\varphi}{}{}(x_0) \geq f(x_0).
\]

It remains to understand the boundary behavior of $\overline{u}$. We will show that the boundary condition is attained in a classical sense, that is $\overline{u} = g$. Let $x \in \p\Omega$ and $p_k$ be the continuous quadratic constructed during the proof of existence of the boundary barrier function in Proposition~\ref{prop:discrbdrybarrier} with constant $E =k$. As $k$ can be taken arbitrarily large, the sequence of points where $g \pm p_k$ attains a maximum (minimum) over $\p\Omega$, converges to $x$.

We now observe that the monotonicity of Lemma~\ref{lem:twoscalemonotone} implies that if $v_h \in X_h$ is such that $\MAop{v_h}{h,\delta,\theta}{2S}(x_h) > 0$ for all $x_h \in \Nhi$, then $v_h$ attains its maximum on $\p\Omega$. Since
\[
  \MAop{u_h^\eps + \calI_h p_k}{h,\delta,\theta}{2S}(x_h)>0, \quad \forall x_h \in \Nhi,
\]
we can apply this observation to $u_h^\eps + \calI_h p_k$ to obtain that, for $x \in \p\Omega$,
\begin{align*}
  \overline{u}(x) &\leq \limsup_{\eps \to 0, \frac{h}\delta \to 0, y \to x} \left( u_h^\eps(y) + \calI_h p_k(y) \right) - \liminf_{\eps \to 0, \frac{h}\delta \to 0, y \to x} \calI_h p_k(y) \\ &\leq 
  \limsup_{\eps \to 0, \frac{h}\delta \to 0, y \to x} \max_{z \in \p\Omega} \calI_h(g(z) + p_k(z)) - p_k(x) \leq g(x_k) + p_k(x_k) - p_k(x),
\end{align*}
where $x_k$ is the point where $g+p_k$ attains its maximum over $\p\Omega$. Letting $k \to \infty$ we conclude $\overline{u} \leq g$. Similarly $\underline{u} \geq g$.

Finally we invoke the comparison principle of Proposition~\ref{prop:detisvarhascompare} to conclude.
\end{proof}

\begin{remark}[convergence by regularization]
\label{rem:otherconvergencetwoscale}
It is interesting to note that by invoking the continuous dependence result given in Proposition~\ref{prop:contdependencevisco}, and the approximation result of Proposition~\ref{prop:smoothapproxvisco}, another proof of convergence can be developed. See \cite[Section 5.3]{NochettoNtogkasZhangfirsttwoscale} for details.
\eremk\end{remark}

\subsubsection{Rates of convergence}
\label{subsub:twoscalerates}

The ingredients used to assert the convergence of the two scale method \eqref{eq:twoscale} were employed in \cite{NochettoNtogkasZhang} to obtain rates of convergence. The techniques used in this reference were very similar to those that we will describe in Section~\ref{sec:OP} and so, to avoid repetition, we shall not elaborate on them here. This is further justified by that fact that, although \cite{NochettoNtogkasZhang} was the first work to provide rates of convergence for wide stencil--type methods, the rates of convergence obtained in this work were suboptimal.

Let us here, instead, present the results obtained in \cite{LiNochettoZhangMA}, where optimal rates of convergence have been obtained. The main tools in this are the comparison principle of Proposition~\ref{prop:tscompare} and the discrete barriers constructed in Section~\ref{subsub:twoscaleconsistent}.

We begin by noticing that we shall only assume
\[
  f \geq 0,
\]
so that the \MA equation \eqref{eqn:MA} may be degenerate. The main result about rates of convergence for classical solutions is the following.

\begin{theorem}[error estimate]
\label{thm:ratetwoscale}
Let $u \in C^{2,\alpha}(\bar\Omega)$, with $\alpha \in (0,1]$, solve \eqref{eqn:MA} and $u_h^\eps \in X_h$ solve \eqref{eq:twoscale}. If $\theta \leq \tfrac1{4d}$ then we have
\[
  \| u - u_h^\eps \|_{L^\infty(\Omega)} \leq C \left[ h^2\left( 1+\delta^{-2}\right)|u|_{C^{1,1}(\bar\Omega)} + \delta^\alpha|u|_{C^{2,\alpha}(\bar\Omega)}  \right],
\]
where the constant $C$ depends on the domain $\Omega$, the dimension $d$, and the shape regularity of the mesh $\mct$, but is independent of $h$, and the solution $u$.
\end{theorem}
\begin{proof}
Recall that a standard interpolation estimate yields
\[
  \| u - \calI_h u \|_{L^\infty(\Omega)} \leq C h^2 |u|_{C^{1,1}(\bar\Omega)},
\]
so that we only need to bound the difference $u_h - \calI_h u$. To do so, we will construct a suitable discrete subsolution $u_h^-$ and supersolution $u_h^+$ and use the comparison principle of Proposition~\ref{prop:tscompare}.

Let $u_h^- = \calI_h u + K_1 q_h \in X_h$, where $q_h$ is the interior barrier of Remark~\ref{rem:interiorbarrier} and $K_1>0$ is to be chosen later. Notice that, by construction
\[
  u_h^- \leq \calI_h u = \calI_h g, \quad \text{on } \p\Omega.
\]
Thus, to guarantee that this is a subsolution we must show that
\[
  \MAop{u_h^-}{h,\delta,\theta}{2S}(x_h) \geq f(x_h) = \det D^2 u(x_h), \quad \forall x_h \in \Nhi.
\]
However, since $u_h^-$ is discretely convex, showing this inequality reduces to showing that, for all $\{\bw_i\}_{i=1}^d \in \calV_\theta$ we have
\[
  \prod_{i=1}^d \nabla^2_{\delta \bw_i} u_h^-(x_h) \geq \det D^2u(x_h), \quad \forall x_h \in \Nhi,
\]
see Lemma~\ref{lem:discrconvex}. Using the convexity of $u$, we have, according to Lemma~\ref{lem:2nddiffconsistent}, that
\[
  \nabla^2_{\delta \bw_i} \calI_h u(x_h) \geq \frac{\p^2 u(x_h)}{\p \bw_i^2} - C|u|_{C^{2,\alpha}(\bar\Omega)} \delta^\alpha,
\]
so that, upon choosing
\[
  K_1 = C\left[ \delta^\alpha |u|_{C^{2,\alpha}(\bar\Omega)} + \left( \frac{h^2}{\delta^2} + \theta^2 \right)|u|_{C^{1,1}(\bar\Omega)} \right],
\]
where $C$ is sufficiently large, 
we have
\begin{align*}
  \nabla^2_{\delta \bw_i} u_h^-(x_h) &\geq \frac{\p^2 u(x_h)}{\p \bw_i^2} - C|u|_{C^{2,\alpha}(\bar\Omega)} \delta^\alpha + K_1 \geq \frac{\p^2 u(x_h)}{\p \bw_i^2} + C\theta^2 |u|_{C^{1,1}(\bar\Omega)} \\
    &\geq 
  \left( 1 + 16\theta^2(d-1)^2\right)\frac{\p^2 u(x_h)}{\p \bw_i^2} \geq \left( 1 + 16\theta^2(d-1)^2\right)^{1/d}\frac{\p^2 u(x_h)}{\p \bw_i^2}.
\end{align*}
Finally, since $\theta \leq \tfrac1{4d}$, we multiply this inequality over $i = 1,\ldots,d$ and invoke Theorem~\ref{thm:twoscaleconsistencyin} item~\ref{it:detisclose} to conclude that $u_h^-$ is a subsolution. The comparison principle of Proposition~\ref{prop:tscompare} then yields that
\begin{align*}
  u_h^\eps \geq u_h^- 
  &= \calI_h u + C\left( \delta^\alpha |u|_{C^{1,1}(\bar\Omega)} + \left(\frac{h^2}{\delta^2}+\theta^2\right) |u|_{C^{1,1}(\bar\Omega)} \right)q_h\\ 
  &\geq \calI_h u - C\left( \delta^\alpha |u|_{C^{1,1}(\bar\Omega)} + \theta^2 |u|_{C^{1,1}(\bar\Omega)} \right).
\end{align*}

We now define
\[
  u_h^+ = \calI_h u - K_1 q_h -K_2 b_h,
\]
where $q_h$ and $K_1$ are as before, $b_h$ is the barrier of Proposition~\ref{prop:otherbarrier} and $K_2>0$ is to be chosen. We show that $u_h^+$ is a supersolution. First of all, because of the choice of signs
\[
  u_h^+ \geq \calI_h u = \calI_h g, \quad \text{on } \p\Omega.
\]
Now, to show the inequality between operators we must consider in $\Omega_\delta$ and outside of it separately. Let $x_h \in \Nhi \cap \Omega_\delta$ and $\{\bv_i\}_{i=1}^d \in \calV$ such that
\[
  f(x_h) = \det D^2 u(x_h) = \prod_{i=1}^d \frac{\p^2 u(x_h)}{\p \bv_i^2}.
\]
Let now $\{\bv_i^\theta\}_{i=1}^d \in \calV_\theta$ satisfy \eqref{eq:calVthetaisclose}. The interior consistency of second differences of Lemma~\ref{lem:2nddiffconsistent}, together with the estimate of Theorem~\ref{thm:twoscaleconsistencyin} item~\ref{it:detisclose2} gives us that
\[
  \left| \nabla_{\delta \bv_i^\theta}^2 \calI_h u(x_h) - \frac{\p^2 u(x_h)}{\p \bv_i^2} \right| \leq C\left[ \delta^\alpha |u|_{C^{1,1}(\bar\Omega)} + \left( \frac{h^2}{\delta^2} + \theta^2 \right)|u|_{C^{1,1}(\bar\Omega)} \right],
\]
which, using that $\nabla_{\delta \bv_i^\theta}^2 q_h(x_h) \geq 1$, $\nabla_{\delta \bv_i^\theta}^2 b_h(x_h) \geq 0$, and the definition of $K_1$ immediately implies that
\[
  \nabla_{\delta \bv_i^\theta}^2 u_h^+(x_h) \leq \frac{\p^2 u(x_h)}{\p \bv_i^2}.
\]

Notice now that $u_h^+$ might not be discretely convex, so that $\nabla_{\delta \bv_i^\theta}^2 u_h^+(x_h)$ might be negative. To deal with this we define the function
\[
  G: \R^d \to \R, \qquad G(\bz) = \prod_{i=1}^d (\bz\cdot\be_i)^+  - \sum_{i=1}^d (\bz\cdot\be_i)^-,
\]
where $\{\be_i\}_{i=1}^d$ is the canonical basis of $\R^d$. Notice that this function is monotone in each coordinate of $\bz$. Moreover if, for $\{\bw_i\}_{i=1}^d \in \calV_\theta$ and $w_h \in X_h$, we define the vectors
\begin{align*}
  \bchi(w_h,\{\bw_i\}) 
  &= \left( \nabla_{\delta \bw_1}^2 w_h(x_h), \ldots, \nabla_{\delta \bw_d}^2 w_h(x_h) \right)^\intercal,\\
  \bgamma &= \left( \frac{\p^2 u(x_h)}{\p \bv_1^2}, \ldots, \frac{\p^2 u(x_h)}{\p \bv_d^2} \right)^\intercal. 
\end{align*}
Then we have that
\[
  \MAop{w_h}{h,\delta,\theta}{2S}(x_h) = \min_{\{\bw_i\}_{i=1}^d \in \calV_\theta} G(\bchi(w_h,\{\bw_i\})).
\]
Therefore
\[
  \MAop{u_h^+}{h,\delta,\theta}{2S}(x_h) \leq G(\bchi(u_h^+,\{\bv_i^\theta\})) \leq G(\bgamma) = \prod_{i=1}^d \frac{\p^2 u(x_h)}{\p \bv_i^2} = f(x_h).
\]

Consider now a node close to the boundary, that is $x_h \in \Nhi\setminus \Omega_\delta$, and let $\{\bw_i^\theta\}_{i=1}^d \in \calV_\theta$. Using Proposition~\ref{prop:otherbarrier} item \ref{it:otherbarrierbdry} we have that
\[
\max_{i=1,\ldots,d} \nabla_{\delta \bw_i^\theta}^2 b_h(x_h) \geq \frac1{2d}.
\]
Assume that this maximum is attained for index $k$. Using Lemma~\ref{lem:2nddiffconsistent} we can conclude that
\begin{align*}
  \nabla_{\delta \bw_k^\theta}^2 u_h^+(x_h) 
  &\leq \nabla_{\delta \bw_k^\theta}^2 \calI_h u(x_h) - K_2 \nabla_{\delta \bw_k^\theta}^2 b_h(x_h)\\
  & \leq \nabla_{\delta \bw_k^\theta}^2 \calI_h u(x_h) - \frac1{2d}K_2 \leq 
  C |u|_{C^{1,1}(\bar\Omega)} - \frac1{2d}K_2 \leq 0,
\end{align*}
where the last step holds upon choosing $K_2$ sufficiently large. This shows that
\[
\min_{i=1,\ldots,d}  \nabla_{\delta \bw_i^\theta}^2 u_h^+(x_h) = 0 \quad \implies \MAop{u_h^+}{h,\delta,\theta}{2S}(x_h) = 0 \leq f(x_h).
\]

We have shown that, for all $x_h \in\Nhi$, we have $\MAop{u_h^+}{h,\delta,\theta}{2S}(x_h) \leq f(x_h)$, so that $u_h^+$ is a supersolution. The discrete comparison principle of Proposition~\ref{prop:tscompare} then allows us to conclude that
\[
  u_h \leq u_h^+ = \calI_h u - K_1 q_h -K_2 b_h \leq \calI_h u + C_1 K_1 + C_2 \delta^2 K_2,
\]
where we used the lower bounds on $q_h$ and $b_h$. Recalling the choices of $K_1$ and $K_2$ allows us to conclude.
\end{proof}

Choosing relations between the discretization parameters $h$, $\delta$ and $\theta$ we can obtain explicit rates of convergence.

\begin{corollary}[rates of convergence]
\label{cor:explicitratetwoscale}
In the setting of Theorem~\ref{thm:ratetwoscale}, if $\delta = C_1 h^{\frac2{2+\alpha}}$ and $\theta = C_2h^{\frac2{2+\alpha}}$, we have
\[
  \| u - u_h^\eps \|_{L^\infty(\Omega)} \leq C h^{\frac{2\alpha}{2+\alpha}}.
\]
On the other hand, choosing $\delta = h^{2/3}$ and $\theta = h^{1/3}$, then we have
\[
  \| u - u_h^\eps \|_{L^\infty(\Omega)} \leq C h^{\frac{2\alpha}{3}}.
\]
In both estimates the hidden constant is independent of $h$.
\end{corollary}

Notice that both choices of relations between the coarse parameters and the mesh size $h$ in Corollary~\ref{cor:explicitratetwoscale} have its benefits and drawbacks. While the first choice yields a faster
rate of convergence, it requires knowledge of the regularity of $u$. On the other hand, the second choice yields a slower convergence rate, but does not require {\it a priori} knowledge of the smoothness of $u$.

\begin{remark}[error estimates under different assumptions]
\label{rem:extendratestwosacle}
The results of Theorem~\ref{thm:ratetwoscale} have been extended in \cite{LiNochettoZhangMA} in several directions:
\begin{enumerate}
  \item {\it Smoother solutions}: If $u \in C^{3,\alpha}(\bar\Omega)$ {\it mutatis mutandis} the proof of Theorem~\ref{thm:ratetwoscale} it follows a rate of convergence. The discretization parameters can be related to each other in such a way that the error is $\calO(h)$, and numerical experiments indicate that this is sharp. 
  
  \item {\it Estimates for solutions with Sobolev regularity}: Assuming that $u \in W^{s,p}(\Omega)$ with $s \leq 3$ and $s-d/p>2$, and that $D^2 u(x) \geq \lambda I$, it has been shown \cite[Theorem 5.7]{LiNochettoZhangMA} that we have 
  \[
    \| u - u_h^\eps \|_{L^\infty(\Omega)} \leq C \left( \frac{h^2}{\delta^2} + \theta^2 + \delta^2 + \frac{\delta^{s-2}}\lambda \right),
  \]
  where the constant depends on the smoothness of $u$. Once again, the discretization parameters can be optimized to obtain a rate $\calO(h^{2-4/s})$.
\end{enumerate}
\eremk\end{remark}

\subsection{Extensions, generalizations, and applications}
\label{sub:otherFDstuff}

We conclude the discussion on finite difference schemes and its variants by briefly describing some  connections, extensions, generalizations, and applications of the schemes discussed here.

\subsubsection{Hamilton Jacobi Bellman formulation and semi-Lagrangian schemes}
\label{subsub:FengJensen}

Let
\[
  \Lambda = \left\{ \blam \in \R^d: \blam_i \geq 0, \ i=1,\ldots,d, \ \sum_{i=1}^d \blam_i = 1\right\}.
\]
Define the function $h: \polS^d \times \R_+ \to \R$ by
\[
  h(M,t) = \mathop{\sup_{\{\bw_i\}_{i=1}^d \in \calV}}_{\blam \in \Lambda} \left[ -\frac1d \sum_{i=1}^d \blam_i \bw_i \cdot M \bw_i + t^{1/d} \left( \prod_{i=1}^d \blam_i \right)^{1/d} \right].
\]
The following result is from \cite{KrylovBook88}, see also \cite[Proposition 6.13]{NSWActa}.

\begin{proposition}[determinant]
\label{prop:Krylov}
For $M \in \polS^d$ and $\delta \in \R_+$ we have that
\[
  h(M,\delta) = 0,
\]
if and only if $M \in \polS^d_+$ and $\det M = \delta$.
\end{proposition}

This motivates to define the function $F_{HJB}: \bar\Omega \times \R \times \polS^d$ by
\[
  F_{HJB}(x,r,M) = \begin{dcases}
                     h(M,f(x)), & x \in \Omega, \\
                     g(x) - r, & x\in \p\Omega,
                   \end{dcases}
\]
and consider the problem: find $u \in C(\bar\Omega)$ that is a viscosity solution of
\begin{equation}
\label{eq:MAasHJB}
  F_{HJB}(x,u(x),D^2u(x)) = 0, \quad x \in \bar\Omega.
\end{equation}
It turns out that this problem has an intimate connection with \eqref{eqn:MA}, as shown in \cite[Theorems 3.3 and 3.5]{FengJensen}.

\begin{theorem}[equivalence]
\label{thm:MAasHJB}
Let $f\in C(\Omega)$ be nonnegative. The function $u \in C(\Omega) \cap B(\bar\Omega)$ is a viscosity solution of \eqref{eq:MAasHJB}, in the sense of Definition~\ref{def:viscosol}, if and only if it is a viscosity solution on the set of convex functions of \eqref{eqn:MA}, in the sense of Definition~\ref{def:viscosolMA}.
\end{theorem}

It is remarkable that the convexity assumption on the solution is not enforced in \eqref{eq:MAasHJB}, it is rather a consequence of the formulation. This motivated \cite{FengJensen} to use \eqref{eq:MAasHJB} for numerical purposes. They proposed a so-called semi-Lagrangian scheme which we now describe. Over a quasiuniform triangulation $\mct$ of size $h>0$ we introduce $X_h$ as the space of piecewise linear and continuous functions. On the basis of \eqref{eq:MAasHJB} we introduce over $X_h$ the operator
\[
  \MAop{w_h}{h,k}{SL}(x_h) = \mathop{\sup_{\{\bw_i\}_{i=1}^d \in \calV}}_{\blam \in \Lambda} \left[
    -\frac1d \sum_{i=1}^d \blam_i \nabla^2_{k \bw_i} w_h(x_h) + f(x_h)^{1/d} \left( \prod_{i=1}^d \blam_i \right)^{1/d}
  \right],
\]
where $x_h \in \Nhi$ and $k>0$ is a discretization parameter. The semi-Lagrangian scheme then seeks for $u_h \in X_h$ such that
\begin{subequations}
\label{eq:SL}
  \begin{align}
    \label{eq:SLin}
    \MAop{u_h}{h,k}{SL}(x_h) &= 0, \quad \forall x_h \in \Nhi, \\
  \label{eq:SLbc}
    u_h(x_h) &= g(x_h), \quad \forall x_h \in \Nhb.
  \end{align}
\end{subequations}
Reference \cite{FengJensen} showed existence and uniqueness of solutions to \eqref{eq:SL} as well as, provided $(h,k) \to 0$ with $\tfrac{h}k \to 0$, convergence to the viscosity solution of \eqref{eq:MAasHJB} and, as a consequence of Theorem~\ref{thm:MAasHJB}, to the viscosity solution of \eqref{eqn:MA} over the set of convex functions. Rates of convergence, however, were not provided.

Although rates of convergence for general semi-Lagrangian schemes were given in \cite[Corollary 7.3]{MR3042570} let us here explore a connection between the solutions of the scheme \eqref{eq:SL} and the two scale method of Section~\ref{sub:twoscale} as described in \cite[Section 6]{LiNochettoZhangMA}. For that one needs to notice, first, that the scheme given in \eqref{eq:SL} is not fully practical. This is because in the operator $\MAop{\cdot}{h,k}{SL}$ the supremum runs over all of $\calV$. We need to introduce a directional discretization by, as before, using $\calV_\theta$ whose elements satisfy \eqref{eq:calVthetaisclose}. With this we introduce the new operator
\[
  \MAop{w_h}{h,k,\theta}{SL}(x_h) = \mathop{\sup_{\{\bw_i\}_{i=1}^d \in \calV_\theta}}_{\blam \in \Lambda} \left[
    -\frac1d \sum_{i=1}^d \blam_i \nabla^2_{k \bw_i} w_h(x_h) + f(x_h)^{1/d} \left( \prod_{i=1}^d \blam_i \right)^{1/d}
  \right],
\]
and denote by $u_h^{(k,\theta)} \in X_h$ the solution to \eqref{eq:SL} but with this new operator. The following is a rather surprising fact. For a proof see \cite[Proposition 6.2]{LiNochettoZhangMA}.

\begin{proposition}[equivalence]
\label{prop:SlandTScoincide}
Let $u_h^\eps \in X_h$ denote the solution to the two scale method \eqref{eq:twoscale} and $u_h^{(k,\theta)} \in X_h$ the solution to the modified semi-Lagrangian scheme with the operator $\MAop{\cdot}{h,k,\theta}{SL}$. In this case, we have $u_h^\eps = u_h^{(k,\theta)}$.
\end{proposition}

From Proposition~\ref{prop:SlandTScoincide} and the results of Section~\ref{subsub:twoscalerates}, rates of convergence for \eqref{eq:SL} can be deduced.

\begin{remark}[nonconvex domains]
\label{rem:nonconvexOmega}
Notice that convexity of the solution is not a constraint in \eqref{eq:MAasHJB} but rather a consequence of it. This has motivated \cite{JensenNonConvex} to explore the possibility of using \eqref{eq:MAasHJB} as an extension of the \MA equation to nonconvex domains, or cases with nonconvex data.
\eremk\end{remark}

\subsubsection{Filtered two scale schemes}
\label{subsub:FilterTwoScale}

In \cite{NochettoNtogkasFilter} the ideas of two scale methods of Section~\ref{sub:twoscale} and those of filtered schemes of Section~\ref{sub:Filtered} were extended to construct a {\it filtered two scale scheme}. Let $\calT_{2h}^2$ be a quasiuniform triangulation of $\bar\Omega$ of size $2h>0$. The superscript in this triangulation indicates that we are doing a quadratic approximation of the boundary. This can be accomplished, for instance, by the use of isoparametric approximation of the boundary; see \cite[Section 10.4]{BrennerBook} and \cite[Section 4.3]{CiarletBook}. Over this mesh we construct $X_{2h}^2$, the space of piecewise quadratic and continuous functions. For $w_{2h} \in X_{2h}^2$ and $x_{2h} \in \Omega_{2h}^i$ we define
\begin{equation}
\label{eq:quadraticMAh}
  \begin{aligned}
    \MAop{w_{2h}}{2h,\tilde\delta,\tilde\theta}{2Sq}(x_{2h}) = \min_{\{\bw_i\}_{i=1}^d \in \calV_{\tilde\theta}} &\left[\prod_{i=1}^d\left(\tilde{\nabla}^2_{\tilde\delta\bw_i} w_{2h}(x_{2h}) \right)^+ \right. \\ &- \left. \sum_{i=1}^d \left( \tilde{\nabla}^2_{\tilde\delta\bw_i} w_{2h}(x_{2h}) \right)^- \right],
  \end{aligned}
\end{equation}
where $\Omega_{2h}^i$ denotes the set of internal degrees of freedom of $X_{2h}^2$, which includes now the vertices and edge midpoints of $\calT_{2h}^2$, and $\tilde{\nabla}^2_{\tilde\delta\bw}$ is a more accurate, say using five points, discretization of the second derivative in direction $\bw$ at scale $\tilde\delta$.

Following the ideas presented in Theorem~\ref{thm:twoscaleconsistencyin} we can show that operator \eqref{eq:quadraticMAh} is consistent with order $\calO( \tilde{\delta}^{k+\alpha} + \tfrac{h^3}{\tilde{\delta}^2} + \tilde{\theta}^2)$; see \cite[Lemma 5.8]{NochettoNtogkasFilter}. However, this scheme is {\it not} monotone. It will, instead serve as the two scale analogue of the accurate scheme \eqref{eq:MAFD}.

By refining in a conforming way once $\calT_{2h}^2$ we obtain $\mct$, over which we can apply the two scale scheme of Section~\ref{sub:twoscale}. Notice that there is a bijection between $\Omega_{2h}^i$ and $\Nhi$ so that the elements of $X_{2h}^2$ and $X_h$ can be compared by looking at their nodal values. In light of this observation we alleviate the notation and carry out the rest of the discussion using the scale $h$.

We combine \eqref{eq:quadraticMAh} and \eqref{eq:TSop} into a {\it filtered} two scale operator: for $w_h \in X_h$ and $x_h \in \Nhi$
\begin{align*}
  \MAop{w_h}{h,\delta,\theta,\tilde\delta,\tilde\theta}{F}(x_h) &= \MAop{w_h}{h,\delta,\theta}{2S}(x_h) \\ &+ \tau \tilde{S}\left(
  \frac{\MAop{w_h}{2h,\tilde\delta,\tilde \theta}{2Sq}(x_h) - \MAop{w_h}{h,\delta,\theta}{2S}(x_h)}\tau\right),
\end{align*}
where $\tilde{S}(t) = \min\{S(t), 0\}$ and the function $S$ is defined in \eqref{eq:exoffilter}. As explained in \cite[Section 2]{NochettoNtogkasFilter} the choice of filter function ensures discrete convexity in the case that the right hand side degenerates, that is if $f(x_h) = 0$, for some $x_h \in \Nhi$.

With these ingredients the filtered two scale scheme seeks for $u_h^F \in X_h$ such that
\begin{subequations}
\label{eq:Filter2S}
  \begin{align}
    \label{eq:Filter2Sin}
      \MAop{u_h^F}{h,\delta,\theta,\tilde\delta,\tilde\theta}{F}(x_h) &= f(x_h), \quad \forall x_h \in \Nhi, \\
    \label{eq:Fiter2Sbc}
      u_h^F(x_h) &= g(x_h), \quad \forall x_h \in \Nhb.
  \end{align}
\end{subequations}
The theory of almost monotone schemes of Corollary~\ref{cor:BSnearmonotone} was combined with the convergence results of Section~\ref{subsub:twoscaleconverge} in \cite[Section 6]{NochettoNtogkasFilter} to assert the convergence of any solution to \eqref{eq:Filter2S}.

\subsubsection{Approximation of convex envelopes}
\label{sub:convexenvelope}

Let us describe the results obtained in \cite{LiNochettoCE} regarding the approximation of the convex envelope of a function, which was introduced in Definition~\ref{def:convexenvelope}. Let $f \in C(\bar\Omega)$. As shown in \cite{MR3621814} the convex envelope $u = \Gamma f$ of $f$ can be characterized as the viscosity solution of the problem
\begin{subequations}
\label{eq:CE}
  \begin{alignat}{2}
    \label{eq:CEin}
      \mathrm{CE}[u](x) &= 0, &&x \in \Omega, \\
    \label{eq:CEbc}
      u(x)&= f(x), &&x \in \p\Omega,
  \end{alignat}
\end{subequations}
where the operator $\mathrm{CE}[\cdot]$ is given by
\begin{equation}
\label{eq:CEop}
  \mathrm{CE}[w](x) = \min\left\{ f(x) - u(x), \min \sigma \left( D^2w(x) \right) \right\}.
\end{equation}

The intuition behind \eqref{eq:CE} is clear. First, we have that $u(x) \leq f(x)$ for every $x \in \bar\Omega$. In addition, if we define the {\it contact set}
\[
  \calC(f) = \left\{ x \in \bar \Omega: u(x) = f(x) \right\},
\]
we obtain, upon denoting $\lambda_1(w) = \min \sigma \left( D^2w(x) \right)$, that for $x \in \calC(f)$ we must have $\lambda_1(u) \geq 0$. On the other hand, if $x \notin \calC(f)$, then we must have $\lambda_1(u) = 0$. In conclusion, $u$ must be convex.

We remark, however, that problem \eqref{eq:CE} is very degenerate. Indeed, it can be shown, see for instance \cite[Lemma 3.1]{LiNochettoCE}, that if $\dist(x,\calC(f)) \geq d \delta$ and $\bp \in \p u(x)$ there is $\bv \in \R^d$ with $|\bv|=1$ such that
\[
  x_\pm = x \pm \delta \bv, \quad u(x_\pm) = u(x) + \delta \bp\cdot \bv, \quad \nab_{\delta\bv}^2 u(x) = 0, \quad \bp \in \p u (x_\pm).
\]
In other words, if we are sufficiently far away from the contact set $\calC(f)$, then the graph of $u$ is flat in at least one direction. As a consequence, in general, the convex envelope cannot be arbitrarily smooth, regardless of the smoothness of the domain $\Omega$ and data $f$. Indeed, \cite{MR3324933} shows that if $\Omega$ is strictly convex with $\p\Omega \in C^{3,1}$,  and $f \in C^{3,1}(\bar\Omega)$, then $u \in C^{1,1}(\bar\Omega)$, and that this is optimal. This very low regularity is one of the main obstacles in the analysis of numerical schemes for \eqref{eq:CE}.

Formulation \eqref{eq:CE} was already used for numerical purposes in \cite{Oberman08ConEnv} via wide stencil schemes like those presented in Section~\ref{sub:Oberman}. Let us present here, instead, the two scale methods of \cite{LiNochettoCE}. We will follow the notation of Section~\ref{sub:twoscale}. In addition, if $\calS$ denotes the unit ball in $\R^d$ we introduce, in full analogy to \eqref{eq:calVthetaisclose}, a discretization $\calS_\theta$ of $\calS$ such that, for every $\bw \in \calS$ there is $\bw^\theta \in \calS_\theta$ that satisfies
\[
  |\bw - \bw^\theta| \leq \theta.
\]

Over the space of piecewise linear functions $X_h$ subordinated to the triangulation $\mct$ we define
\begin{equation}
\label{eq:CEh}
  \mathrm{CE}_{h,\delta,\theta}[w_h](x_h) = \min\left\{ f(x_h) - w_h(x_h), 
    \min_{\bw \in \calS_\theta} \nab_{\delta \bw}^2 w_h(x_h) \right\}
\end{equation}
where $w_h \in X_h$ and $x_h \in\Nhi$. With the aid of this operator we define the discrete convex envelope of a function $f$ as the function $u_h^\eps \in X_h$ that solves
\begin{subequations}
\label{eq:CEdiscr}
  \begin{align}
    \label{eq:CEdiscrin}
      \mathrm{CE}_{h,\delta,\theta}[u_h^\eps](x_h) &= 0, & x_h \in \Nhi, \\
    \label{eq:CEdiscrbc}
      u_h^\eps(x_h)&= f(x_h), & x_h \in \Nhb.
  \end{align}
\end{subequations}

The analysis of scheme \eqref{eq:CEdiscr} to a large extent follows that of two scale methods presented in Section~\ref{sub:twoscale}. Namely, owing to discrete convexity we can show that the scheme has a comparison principle, from which uniqueness of solutions follows. The existence of solutions is obtained via a discrete Perron method, and the stability by noticing that $u_h^- = \calI_h u$ and $u_h^+ = \calI_h f$ are discrete sub and supersolutions, respectively.

The considerations given above show that scheme \eqref{eq:CEdiscr} is monotone and stable. In addition, assuming smoothness of the arguments, one can show its consistency with similar arguments to those of Section~\ref{subsub:twoscaleconsistent}. Upon realizing that the operator \eqref{eq:CEop} has a comparison principle in the sense of Definition~\ref{def:comparisonVisco}, this is enough to appeal to the theory of Section~\ref{sub:BarlesSouganidis} and conclude that the scheme is convergent as $\eps = (h,\delta,\theta) \to 0$, provided $\tfrac{h}\delta \to 0$.

The derivation of rates of convergence, however, requires special attention. This is due to the fact that, as stated above, the best regularity we can expect is $u \in C^{1,1}(\bar\Omega)$, and this is not enough to exploit the consistency estimates that were used for convergence (which are applied to smooth test functions). To overcome this, one must take advantage of the flatness of the solution outside the contact set. To describe these results we must introduce some notation. Set, for $x_h \in \Nhi$,
\[
  \delta_{x_h} = \min \{ \delta, \dist(x_h,\p\Omega) \}, \quad B_{x_h} = \bigcup_{T \in \mct: \dist(x_h,T) < \delta_{x_h}} T,
\]
and
\[
  W_{x_h} = \left\{ x \in \bar\Omega : |x - x_h| \leq d \delta \right\}.
\]
The following is \cite[Proposition 3.3]{LiNochettoCE}.

\begin{proposition}[consistency]
Let $\Omega$ be strictly convex and $u$, the solution of \eqref{eq:CE} satisfy $u \in C^{k,\alpha}(\bar\Omega)$ with $k = 0,1$ and $\alpha \in (0,1]$. For $x_h \in \Nhi$ we have:
\begin{enumerate}
  \item If $\dist(x_h, \calC(f) ) \geq d \delta$, then 
  \[
    \min_{\bw \in \calS_\theta} \nab_{\delta \bw }^2 \calI_h u(x_h) \leq C 
    \left( \frac{(\delta\theta)^{k+\alpha} + h^{k+\alpha}}{\delta^2} \right)|u|_{C^{k,\alpha}(B_{x_h})}.
  \]

  \item If $\dist(x_h, \calC(f)) < d \delta$ but $\dist(x_h,\p\Omega) \geq d \delta$, then we have
  \[
    f(x_h) - u(x_h) \leq C_k \delta^{k+\alpha},
  \]
  where $C_k$ depends on $|u|_{C^{0,\alpha}(W_{x_h})} + |f|_{C^{0,\alpha}(W_{x_h})}$ for $k=0$ and on $|f|_{C^{1,\alpha}(W_{x_h})}$ for $k=1$.
  
  \item If $0 < \dist(x_h, \p\Omega) < d \delta$, then for all $\bw \in \calS$ we have
  \[
    \nab_{\delta \bw }^2 \calI_h u(x_h) \leq C \delta_i^{k+\alpha - 2 }|u|_{C^{k,\alpha}(B_{x_h})},
  \]
  and the previous item also holds provided $k=0$.
\end{enumerate}
\end{proposition}

To take advantage of this result two new discrete barriers were constructed. One  handles the first case, i.e., points sufficiently far away from the contact set. The other barrier handles points near the boundary, that is the third case in the previous result. Without going into details, we present here the main error estimate, and refer the reader to \cite[Theorem 3.7]{LiNochettoCE}.

\begin{theorem}[convergence rate]
Let $\Omega$ be strictly convex, $u$ be the viscosity solution of \eqref{eq:CE}, and $u_h^\eps$ the solution of \eqref{eq:CEdiscr}. If $u \in C^{k,\alpha}(\bar\Omega)$, with $k = 0,1$ and $\alpha \in (0,1]$, then 
\[
  \| u - u_h^\eps \|_{L^\infty(\Omega)} = \calO\left( h^{\frac{(k+\alpha)^2}{k+\alpha+2}} \right),
\]
provided, $\delta = \calO\left(h^{\frac{k+\alpha}{k+\alpha+2}}\right)$ and $\theta = \calO\left(h^{\frac2{k+\alpha+2}}\right)$. In particular, if $k=\alpha=1$, i.e., $u \in C^{1,1}(\bar\Omega)$, we obtain
\[
  \delta = \calO( h^{1/2}), \quad \theta = \calO(h^{1/2}) \quad \implies \quad 
  \| u - u_h^\eps \|_{L^\infty(\Omega)} = \calO( h).
\]
Similarly, if $k=0$ and $\alpha=1$, i.e, $u\in C^{0,1}(\bar\Omega)$, we get
\[
  \delta = \calO( h^{1/3}), \quad \theta = \calO(h^{2/3}) \quad \implies \quad 
  \| u - u_h^\eps \|_{L^\infty(\Omega)} = \calO( h^{1/3}).
\]
\end{theorem}


\subsubsection{The Gauss curvature problem}
\label{subsub:Gausscurvature}

As an application of the wide stencil finite difference schemes that were presented in Section~\ref{sub:Oberman} let us here, following \cite{FroeseGauss}, describe a discretization of the prescribed Gaussian curvature problem \eqref{eq:Gausscurvature}. To do so, we must begin by defining what is a solution of this problem. In a similar manner to the notion of Alexandrov solutions to \MA problem, introduced in Definition~\ref{def:AlexSolution}, we have

\begin{definition}[generalized solution]
\label{def:AlexSolGauss}
A convex function $u: \bar\Omega \to \R$ is a generalized solution to \eqref{eq:Gausscurvature} if the following two conditions hold:
\begin{enumerate}
  \item It is a generalized solution of \eqref{eq:Gausscurvaturein}. This means that, for all Borel sets $D \subset \Omega$, we have
  \[
    \int_{\p u (D)} \frac1{\left( 1+ |\bp|^2 \right)^{(d+2)/2}} \diff \bp = \int_D \calK(x) \diff x.
  \]

  \item It satisfies
  \[
    \limsup_{y \to x} u(y) \leq g(x), \quad \forall x \in \p \Omega
  \]
  and, if $v$ is any other generalized solution of \eqref{eq:Gausscurvaturein}, then $v \leq u$ in $\Omega$.
\end{enumerate}
\end{definition}

Under the assumptions of uniform convexity of $\Omega$; continuity of $g$; continuity, boundedness, and nonnegativity of $\calK$; and the compatibility condition
\[
  \int_{\R^d} \frac1{\left( 1+ |\bp|^2 \right)^{(d+2)/2}} \diff \bp > \int_\Omega \calK(x) \diff x;
\]
it can be shown that problem \eqref{eq:Gausscurvature} has a unique generalized solution; see \cite{MR1305147}.

It is also possible to extend the notion of viscosity solution presented in Definition~\ref{def:viscosol}, by allowing the operators in Definition~\ref{def:FLelliptic} to also depend on a variable $\bp \in \R^d$. In doing that, we note that the operator $F_{G\calK,c}: \bar\Omega \times \R \times \R^d \times \polS^d \to \R$ defined by
\[
  F_{G\calK,c}(x,r,\bp,M) = \begin{dcases}
                   \det M - \calK(x)\left( 1+ |\bp|^2 \right)^{(d+2)/2}, & x \in \Omega, \\
                   g(x) - r, & x \in \p\Omega,
                 \end{dcases}
\]
is, as the \MA operator $F_{MA}$ defined in \eqref{eq:MaasFNLelliptic}, only elliptic when $M \in \polS^d_+$, which implies that to have a reasonable notion of viscosity solution, we must require sub and supersolutions to be convex, and restrict the test functions to be convex, as in Definition~\ref{def:viscosolMA}. As we have seen throughout our discussion, the convexity constraint is rather difficult to impose explicitly during discretization.

Reference \cite{FroeseGauss} proposed to consider the following formulation of \eqref{eq:Gausscurvature}. If for $M \in \polS^d$ we set $\sigma(M) = \{ \lambda_1(M), \ldots, \lambda_d(M)\}$, where the eigenvalues are counted with multiplicity and arranged in nondecreasing order, then the operator
\[
  F_{G\calK}(x,r,\bp,M) = \begin{dcases}
                   F^{\mathrm{in}}_{G\calK}(x,\bp,M), & x \in \Omega, \\
                   g(x) - r, & x\in \p\Omega,
                 \end{dcases}
\]
with
\[
  F^{\mathrm{in}}_{G\calK}(x,\bp,M) = \min\left\{ \lambda_1(M),\ \prod_{i=1}^d \lambda_i(M)^+ - \calK(x)\left( 1+ |\bp|^2 \right)^{(d+2)/2} \right\}
\]
is elliptic in the sense of Definition~\ref{def:FLelliptic} and, at least formally, it is clear that if
\[
  F_{G\calK}(x,u(x),\nab u(x),D^2 u(x)) = 0, \ x \in \Omega,
\]
then we must have, that either, $\lambda_1(D^2u(x)) >0$, so that $u$ is convex, and
\[
  \det D^2u(x) = \calK(x)\left( 1+ |\nab u(x)|^2 \right)^{(d+2)/2},
\]
or $\lambda_1(D^2u(x)) = 0$ and, thus
\[
  0 = \det D^2 u(x) \geq \calK(x)\left( 1+ |\nab u(x)|^2 \right)^{(d+2)/2} \geq 0.
\]
In either case, the convexity of the solution is recovered.

With these constructions we have two options to define viscosity solutions to \eqref{eq:Gausscurvature}. The first one is, like in Definition~\ref{def:viscosolMA}, to require that it is a viscosity solution, in the set of convex functions, of the problem
\begin{equation}
\label{eq:convexviscosolGauss}
F_{G\calK,c}(x,u(x),\nab u(x),D^2u(x)) = 0, \quad \forall x \in \bar\Omega.
\end{equation}
The second, as in Definition~\ref{def:viscosol}, to require that it is a viscosity solution of
\begin{equation}
\label{eq:viscosolGauss}
F_{G\calK}(x,u(x),\nab u(x),D^2u(x)) = 0, \quad \forall x \in \bar\Omega.
\end{equation}

In full analogy to Proposition~\ref{prop:equiv} it is shown in \cite[Section 3]{FroeseGauss} that viscosity subsolutions to problem \eqref{eq:viscosolGauss} are convex and that a function is a viscosity solution of \eqref{eq:convexviscosolGauss} over the set of convex function if and only if it is a viscosity solution of \eqref{eq:viscosolGauss}. In addition it is shown that, under certain assumptions on $\calK$, this notion of solution, at least in the interior of the domain $\Omega$, coincides with that of Definition~\ref{def:AlexSolGauss}.

It is important to note that incorporating the boundary conditions into the definition of the operator is {\it essential} in this problem, as they may not be realized in a classical sense. The following is \cite[Example 1]{FroeseGauss}. 

\begin{example}[nonclassical boundary conditions]
\label{ex:Gaussbcsareamess}
Let $d = 1$, $\Omega = (0,1)$, and $\calK \equiv 1$. We set the boundary conditions $u(0) = -1$ and $u(1) = 1$. Then it is possible to show that
\[
  u(x) = - \sqrt{1-x^2}
\]
is a viscosity solution of \eqref{eq:viscosolGauss}. It is a classical solution over $[0,1)$ so it remains to understand what happens at $x=1$.

Note that $u'(x)$ grows unboundedly as $x \uparrow 1$ so that it is not possible to find a smooth $\varphi$ such that $u_\star - \varphi$ has a local minimum at $x_0$, in other words, the graph of $u$ cannot be touched from below at $x=1$. This makes $u$ automatically a supersolution.

To show that $u$ is also a subsolution we note that $u(1) = 0 < 1$ so that, if $\varphi$ touches the graph of $u$ from above at $x=1$, we must have $\varphi(1) = u(1) = 0$, and 
\[
  \left(F_{G\calK}\right)_\star(1,u(1),\varphi'(1),\varphi''(1)) \geq 1-u(1) = 1 >0.
\]
\end{example}

The behavior of Example~\ref{ex:Gaussbcsareamess} was characterized in \cite[Corollary 24]{FroeseGauss}. Namely, if $u$ is a viscosity solution of \eqref{eq:convexviscosolGauss} then at every $x \in \p\Omega$ we either have that $u_\star(x) = u^\star(x) =g(x)$, or $u_\star(x) \leq u^\star(x) \leq g(x)$ and $\p u_\star(x) = \emptyset$. The second option here corresponds to the right endpoint in Example~\ref{ex:Gaussbcsareamess}.

Existence of solutions to \eqref{eq:viscosolGauss} was shown using a variant of Perron's method. 
The usual argument to show uniqueness is obtained via a comparison principle of Definition~\ref{def:comparisonVisco}. This problem, however, does not have a comparison principle, as \cite[Example 3]{FroeseGauss} shows.

\begin{example}[lack of comparison]
\label{ex:Gaussdoesnotcompare}
In the setting of Example~\ref{ex:Gaussbcsareamess} we have that $u(x) = -\sqrt{1-x^2}$ is a viscosity solution, so that necessarily it is a supersolution. Let
\[
  v(x) = \begin{dcases}
           u(x), &x \in [0,1), \\ 1, & x = 1,
         \end{dcases}
\]
we see that $v\in USC([0,1])$ and, as in Example~\ref{ex:Gaussbcsareamess}, if $\varphi$ touches from above the graph of $u$ at $x=1$, then $\varphi(1) = v(1) = 1$ and
\[
  \left(F_{G\calK}\right)_\star(1,v(1),\varphi'(1),\varphi''(1)) \geq 1-v(1) = 0,
\]
showing that $v$ is a subsolution. Note, however, that $u(1) \leq v(1)$ and this problem does not have a comparison principle.
\end{example}

The previous result, combined with the behavior of solutions at the boundary shows that, in fact, a comparison principle takes place, but only in the interior of the domain; see \cite[Theorem 7]{FroeseGauss}.

\begin{theorem}[interior comparison]
\label{thm:intcompareGauss}
Let $\underline{u} \in USC(\bar\Omega)$ be a subsolution of \eqref{eq:viscosolGauss} and $\overline{u} \in LSC(\bar\Omega)$ a supersolution. Then we have $\underline{u} \leq \overline{u}$ in $\Omega$.
\end{theorem}

This weakened comparison principle is sufficient to guarantee uniqueness.

Having shown the existence and uniqueness of solutions to \eqref{eq:viscosolGauss}, it is possible now to construct numerical schemes. This is carried using variants of the wide stencil finite difference schemes of Section~\ref{sub:Oberman}. With the notation introduced there we define, for $w_h \in X_h$,
\begin{align*}
  \mathrm{GK}_{h,\theta}[w_h](x_h) = 
    &\min\left\{ \min_{ \{\bnu_i\}_{i=1}^d \in \calG_\theta } \Del_{\bnu_i} w_h(x_h),\ \MAop{w_h}{h,\theta}{WS}(x_h) \right. \\ &- \left. \calK(x_h)\left( 1 + |\nab_h w_h(x_h)|^2 \right)^{(d+2)/2}  \right\},
\end{align*}
where $\MAop{\cdot}{h,\theta}{WS}$ was defined in \eqref{eq:MAhtheta} and the vector $\nab_h w_h(x_h)$ is such that
\begin{equation}
\label{eq:monotonegrad}
  \begin{aligned}
    \nab_h w_h(x_h)\cdot \be_i = &\max\left\{ \frac{w_h(x_h) - w_h(x_h - h \be_i)}h,  \right. \\ & \left. \frac{w_h(x_h) - w_h(x_h + h \be_i)}h ,  0 \right\},
  \end{aligned}
\end{equation}
and $\{\be_i\}_{i=1}^d$ is the canonical basis of $\R^d$. With this operator, the finite difference approximation of \eqref{eq:viscosolGauss} is to find $u_h \in X_h$ such that
\begin{subequations}
\label{eq:FDGausscurvature}
  \begin{align}
  \label{eq:FDGausscurvaturein}
    \mathrm{GK}_{h,\theta}[u_h](x_h) &= 0, \quad x_h \in \Nhi, \\
  \label{eq:FDGausscurvaturebc}
    u(x_h ) &= g(x_h), \quad x_h \in \Nhb.
  \end{align}
\end{subequations}

In \cite[Section 6]{FroeseGauss} it is shown that scheme \eqref{eq:FDGausscurvature} is monotone, in the sense of \eqref{eq:BarlesSougMonotone}, stable, in the sense of \eqref{eq:BarlesSougStable}, and consistent, in the sense of \eqref{eq:BarlesSougConsistent}. Notice, however, that as Example~\ref{ex:Gaussdoesnotcompare} shows, problem \eqref{eq:viscosolGauss} does not have a comparison principle. As a consequence, Theorem~\ref{thm:BarlesSouganidis} cannot be applied. For this reason, the framework of Section~\ref{sub:BarlesSouganidis} was extended in \cite[Theorem 9]{FroeseGauss} to cases where problem \eqref{eq:BVPBarlesSouganidisStyle} only has an interior comparison principle like that of Theorem~\ref{thm:intcompareGauss} and there exist classical sub and supersolutions. The conclusion is the locally uniform convergence of $u_h$ to $u$.

\subsubsection{Transport boundary conditions}
\label{subsub:TransportBCs}

Let us conclude the discussion of wide stencil finite difference schemes by describing how these methods can be used to tackle the optimal transportation problem. Since {this will be one of the main topics of another Chapter in this volume}, we shall be brief.

We recall that, given $\Omega, \calO \subset \R^d$, which we assume bounded, with $\calO$ convex, and measures
$\rho_\Omega : \Omega \to \R$ and $\rho_\calO : \calO \to \R$, the optimal transportation problem (with quadratic cost) seeks for a map $T: \Omega \to \calO$ with $T_\sharp \rho_\Omega = \rho_\calO$ that minimizes
\[
  \frac12 \int_\Omega |x - T(x)|^2 \diff \rho_\Omega(x).
\]
We recall that $T_\sharp \mu$ denotes the pushforward of the measure $\mu$ under the mapping $T$. Assuming that the measures are absolutely continuous with respect to Lebesgue measure, with densities $f_\Omega,f_{\calO}$, this condition can be written as
\[
\int_E f_{\calO} (x)\diff x = \int_{T^{-1}(E)} f_\Omega(x)\diff x,
\]
and so by a change of variables, $\det(\nab T(x))f_{\calO}(T(x)) = f_\Omega(x)$.
Finally, we recall that since the cost is quadratic, it can be shown that $T$ is given by the gradient map of a convex potential $u : \Omega \to \R$. This allows us to, at least at the formal level, rewrite the optimal transportation problem as a \MA problem: find $u : \bar\Omega \to \R$ convex, such that
\begin{equation}
\label{eq:transporteq}
  \det D^2 u(x) = F(x,\nab u(x)), \quad x \in \Omega.
\end{equation}
where we set $F(x,\bp) = \rho_\Omega(x)/\rho_\calO(\bp)$. This problem is supplemented by the so--called {\it transport} or {\it second} boundary condition
\[
  \nab u(\bar\Omega)= \bar\calO.
\]
Notice that this, more than a boundary condition, is a set of constraints. It can be shown also that this condition can be replaced by
\begin{equation}
\label{eq:transportbc}
  \nab u(\p\Omega)= \p\calO.
\end{equation}
Thus, we want to construct numerical schemes to approximate the solution of \eqref{eq:transporteq}---\eqref{eq:transportbc}.

It is clear that the main issue is the discretization of the boundary condition \eqref{eq:transportbc}. If the boundary of the domain $\calO$ is given as the zero level set of some function $\Phi:\R^d \to \R$, then it is clear that \eqref{eq:transportbc} can be equivalently written as
\[
  \Phi(\nab u(x)) = 0, \quad \forall x \in \p\Omega.
\]
While we would be tempted to discretize this condition directly, the function $\Phi$ can be highly nonlinear and nonsmooth, which will make the design of monotone and consistent numerical schemes a daunting task. However, this can be achieved very easily if the domains are rectangles, say $\Omega = (0,1)^2 =\calO$. In this case, it is shown in \cite[Section 3.2]{FroeseTransportBC} that each side must be mapped to itself. If we consider the left side of the square, that is,
\[
  \{ (x_1,x_2) \in \R^2: x_1 = 0, x_2 \in [0,1]\},
\]
then the function that describes this is given by $\Phi(y_1,y_2) =y_1$. Thus, on this side we can write
\[
  \frac{\p u(0,x_2)}{\p x_1} = 0.
\]
Similarly, in the right, bottom and top sides, respectively, we can write
\[
  \frac{\p u(1,x_2)}{\p x_1} = 1, \quad \frac{\p u(x_1,0)}{\p x_2} = 0, \quad \frac{\p u(x_1,1)}{\p x_2} = 1.
\]
It is remarkable that on all sides the derivative that appears is actually the normal derivative. This motivated \cite{FroeseTransportBC} to replace the boundary condition \eqref{eq:transportbc} by a Neumann--type boundary condition
\[
  \frac{\p u(x)}{\p \bn} = \phi(x)
\]
for some unknown function $\phi$.

Obviously, the correct choice of function $\phi$ is $\phi(x) = \nab u(x) \cdot \bn(x)$, which motivates the introduction of the following iterative scheme: Given $u_0$, an initial guess, then
\begin{itemize}
  \item \texttt{For} $k \geq 0$ 
  \begin{itemize}
    \item Define, for $x \in \p \Omega$,
    \begin{equation}
    \label{eq:defofProj}
      \bp_k(x) = \PRoj_{\p \calO}(\nab u_k(x)),
    \end{equation}
    where by $\PRoj_S(\bw)$ we denoted a projection of the vector $\bw$ onto the set $S$.
    
    \item Find $u_{k+1}: \bar\Omega \to \R$ convex,  and $c_{k+1} \in \R$ that satisfy
    \begin{subequations}
    \label{eq:OTiterative}
      \begin{align}
      \label{eq:OTiterativecompat}
        \int_\Omega u_{k+1}(x) \diff x &= 0, \\
      \label{eq:OTiterativeeq}
        \det D^2 u_{k+1}(x) &= c_{k+1}F(x,\nab u_{k+1}(x)), \quad x \in \Omega, \\
      \label{eq:OTiterativebc}
        \frac{\p u_{k+1}(x)}{\p \bn} &= \bp_k(x), \quad x \in \p\Omega.
      \end{align}
    \end{subequations}
    
    \item \texttt{Set} $k \leftarrow k+1$
  \end{itemize}
  \item \texttt{EndFor}
\end{itemize}

\begin{remark}[iterative scheme]
\label{rem:OTiterative}
The iterative scheme \eqref{eq:defofProj}--\eqref{eq:OTiterative} deserves several observations.
\begin{enumerate}
  \item The introduction of the projection $\bp_k$ in step \eqref{eq:defofProj} is due to the fact that there is no reason to expect that $\nab u_k(\p \Omega) \subset \p \calO$. Thus, we settle for the closest point on the target boundary.
  
  \item Problem \eqref{eq:OTiterative} is a Neumann problem for an elliptic equation so that the solution, if it exists, is unique only up to a constant. Condition \eqref{eq:OTiterativecompat} forces uniqueness, while the introduction of the number $c_{k+1}$ in \eqref{eq:OTiterativeeq} relaxes the equation so that the necessary conditions for existence are fulfilled.
  
  \item The initialization of this scheme can done by choosing $\bp_0 = M x \cdot \bn$, where $\bn$ is the unit outer normal to $\p \Omega$ and $M>0$ is so large that the image of the mapping $\bar\Omega \ni x \mapsto Mx \in \R^d$ contains $\bar\calO$.
\end{enumerate}
\eremk\end{remark}

We are then going to discretize \eqref{eq:defofProj}--\eqref{eq:OTiterative}. Notice that now the boundary conditions \eqref{eq:OTiterativebc} are rather standard and can be approximated by, for instance, introducing a layer of ghost nodes near the boundary and computing centered differences.

It remains to discretize \eqref{eq:OTiterativeeq}. Setting $v_h = u_{h,k+1}$, the first alternative, proposed in \cite{FroeseTransportBC}, is to use
\begin{equation}
\label{eq:easydiscrOT}
  \MAop{v_h}{h,\theta}{WS}(x_h) = F(x_h, \nab_h v_h), \quad x_h \in \Nhi,
\end{equation}
where $\MAop{\cdot}{h,\theta}{WS}(x_h)$ is the operator defined in \eqref{eq:MAhtheta} or Remark~\ref{rem:variantWS}, and $\nab_h v_h$ is defined as in \eqref{eq:monotonegrad}. Another option, also from \cite{FroeseTransportBC}, is to take advantage of the directional difference that are already being computed to approximate the \MA operator. Notice that if $\{\bnu_i\}_{i=1}^d \in \calV$, we have that 
\[
  \nab w = \left( \frac{\p w}{\p x_1}, \dots, \frac{\p w}{\p x_d}\right)^\intercal = 
    \left( \sum_{i=1}^d \frac{\p w}{\p \bnu_i} \bnu_i \cdot \be_1, \ldots, \sum_{i=1}^d \frac{\p w}{\p \bnu_i} \bnu_i \cdot \be_d\right)^\intercal,
\]
where, as usual, $\{\be_i\}_{i=1}^d$ is the canonical basis of $\R^d$. This allows us to write that
\begin{multline*}
  \det D^2 w(x)- F(x,\nab w(x)) = \MAop{w}{}{}(x) - F(x,\nab w(x))\\
   =
    \min_{\{\bw_i\}_{i=1}^d \in \calV} \left[ \prod_{i=1}^d\left(\frac{\p^2 w(x)}{\p \bw_i^2} \right)^+ - \sum_{i=1}^d \left( \frac{\p^2 w(x)}{\p \bw_i^2} \right)^- \right] \\
  = \min_{\{\bw_i\}_{i=1}^d \in \calV} \left[ \prod_{i=1}^d\left(\frac{\p^2 w(x)}{\p \bw_i^2} \right)^+ - \sum_{i=1}^d \left( \frac{\p^2 w(x)}{\p \bw_i^2} \right)^- - F(x,\nab w(x) )\right] \\
  = \min_{\{\bw_i\}_{i=1}^d \in \calV} \Bigg[ \prod_{i=1}^d\left(\frac{\p^2 w(x)}{\p \bw_i^2} \right)^+ - \sum_{i=1}^d \left( \frac{\p^2 w(x)}{\p \bw_i^2} \right)^-\\ 
  - F\left(x,\left( \sum_{i=1}^d \frac{\p w(x)}{\p \bnu_i} \frac{\bnu_i \cdot \be_1}{|\bnu_i|}, \ldots, \sum_{i=1}^d \frac{\p w(x)}{\p \bnu_i} \frac{\bnu_i \cdot \be_d}{|\bnu_i|}\right)^\intercal \right)\Bigg] \\
  = \min_{\{\bw_i\}_{i=1}^d \in \calV} \mathrm{OT}_{\{\bw_i\}_{i=1}^d}[w](x).
\end{multline*}
In conclusion, an approximation of \eqref{eq:OTiterativeeq} is obtained by setting
\[
  \mathrm{OT}_{h,\theta}[v_h](x_h) = 0, \quad \forall x_h \in \Nhi,
\]
where
\[
  \mathrm{OT}_{h,\theta}[w_h](x_h) = \min_{\{\bnu_i\}_{i=1}^d \in \calG\theta} \mathrm{OT}_{\{\bnu_i\}_{i=1}^d}[w_h](x_h).
\]

Reference \cite{BenamouFroeseObermanTBC} considered a different treatment of the boundary condition \eqref{eq:transportbc}. Since for all $x \in \p\Omega$ we must have that $\nab u(x) \in \p \calO$, then we must have
\begin{equation}
\label{eq:HJBinbdy}
  H(\nab u(x)) = 0, \qquad H(y) = \begin{dcases}
                                    \dist(y,\p \calO), & y \in \calO, \\
                                    -\dist(y,\p\calO), & y \notin \calO,
                                  \end{dcases}
\end{equation}
where $H$ is nothing but the signed distance function to $\p\calO$. Notice that \eqref{eq:HJBinbdy} is a sort of Hamilton Jacobi equation posed on $\p\Omega$. Exploiting the convexity of $\calO$, the authors of \cite{BenamouFroeseObermanTBC} were able to rewrite the function $H$ as the supremum over linear expressions on $y$ (the supporting hyperplanes of $\calO$ at $y$)
\[
  H(y) = \sup_{\bn \in \R^d: |\bn| = 1} \left \{ y \cdot \bn - H^\star(\bn): \bn \cdot \bn_x > 0 \right\},
\]
where $\bn_x$ is the normal to $\p\Omega$ at $x$ and $H^\star$ is the support function of $\calO$, that is,
\[
  H^\star(\bn) = \sup_{\bz \in \p\calO} \bz \cdot \bn.
\]
This function can be precomputed or evaluated rather cheaply in the discrete setting. The reformulation of the function $H$ can be approximated by replacing the supremum by one over a finite set of directions, and, finally, the gradient appearing in \eqref{eq:HJBinbdy} can be discretized as in \eqref{eq:monotonegrad}. This gives a discretization of \eqref{eq:transportbc}. Finally, the discretization of \eqref{eq:transporteq} is proposed to be carried similarly to \eqref{eq:easydiscrOT}.


\section{Discretizations based on geometric considerations}
\label{sec:OP}

\epigraph{\tiny{\it In fact, geometrical representations, graphs and diagrams of all sorts, are used in all sciences, not only in physics, chemistry, and the natural sciences, but also in economics, and even in psychology. Using some suitable geometrical representation, we try to express everything in the language of figures, to reduce all sorts of problems to problems of geometry.}}{\tiny{G.~P\'olya \cite{MR3289212}}}

In this section we will describe the so-called Oliker-Prussner method, which is a discrete analogue of the notion of solution in the Alexandrov sense. We recall that Alexandrov solutions to the \MA equation were introduced in Definition~\ref{def:AlexSolution}. They make a connection between the \MA equation and the measure of the subdifferential of its solution. This, very geometric, notion enables us to define solutions that are not smooth, say not $C^2(\dm)$. The Oliker-Prussner method, in turn, will allow us to approximate these solutions.

\subsection{Description of the scheme}
\label{subsec:OP}
To be able to present the Oliker-Prussner method, we must begin by introducing some useful notions.

\subsubsection{Nodal set and domain partition}
To discretize the domain $\dm$ and its boundary $\bdry$, we introduce a translation invariant nodal set 
and an open, disjoint partition
of the domain. For a parameter $h>0$, we define the interior nodal set as
\begin{align}\label{translationinvariant}
\Omega_h = \left\{ x_h = h\sum_{j=1}^d z^j \tilde \be_j: z^j \in \Z \right\} \cap \Omega,
\end{align}
where $\{\tilde \be_j\}_{j=1}^d$ is a basis of $\bbR^d$ with $|e_j|\le 1$ for all $1\le j\le d$. To discretize the boundary $\bdry$, we set the boundary nodal set $\partial \Omega_h$ as a collection of points on the boundary and require that their spacing is of order $h$, namely, $\bdry \subset \cup_{x_h \in \partial \Omega_h} B_{h/2}(x_h)$.
We set the nodal set $\bar\Omega_h = \Omega_h \cup \partial \Omega_h$. 
We remark that this is a generalization of the finite difference discretizations introduced in Section~\ref{sub:Oberman}. Indeed, in that case the vectors $\{\tilde \be_j\}_{j=1}^d$ were the canonical basis of $\R^d$.

We define an open, disjoint partition {$\{ \omega_{x_h} \}_{x_h\in \bar\Omega_h}$ of the domain where, for $x_h\in \bar\Omega_h$,
\begin{align}\label{partition_MA}
  \omega_{x_h} = \left\{ x_h + \sum_{j=1}^d h^j \tilde \be_j : \ h^j \in \R, \ |h^j| \leq \frac h2 \right\} \cap \Omega.
\end{align}
Note that, by construction, the partition is translation invariant, that is, $\omega_{y_h} = y_h - x_h + \omega_{x_h}$ for all $x_h,y_h\in \Omega_h$ with $\omega_{y_h}, \omega_{x_h} \subset \Omega$.

\subsubsection{Nodal functions, their subdifferentials, and convex envelopes }
On the nodal set $\Omega_h$ constructed above, we define a nodal function $u_h$ to approximate the solution of the \MA PDE.
First, to mimic the convexity constraint for the PDE, we require the notion of convexity for nodal functions (compare to Definition \ref{def:discrconvex}).

\begin{definition}[nodal convexity]
\label{def:nodalConvex}
Let $w_h$ be a (nodal) function that maps the set of nodes $\bar{\Omega}_h$ to $\R$. 
We say that $w_h$ is {\it convex} if, for any node $x_h \in \bar{\Omega}_h$, there exist an affine function $L$, that is, $L(x) = \bp \cdot (x - x_h) + c$ for some $\bp \in \R^d$ and $c\in \R$,  such that
\begin{align}
  \label{def:convexnodalfunction}
    L(y_h) \leq w_h(y_h) \quad \forall y_h \in \bar{\Omega}_h
    \quad \mbox{and} \quad
    L(x_h) = w_h(x_h).
  \end{align}
\end{definition}
We define the subdifferential of a convex nodal function $w_h$ at a fixed node $x_h \in \bar{\Omega}_h$ as the set
\begin{align}\label{def:subdifferential}
  \p w_h(x_h) := \{ \bp \in \R^d:\ \bp \cdot (y_h - x_h) + w_h(x_h) \leq w_h(y_h)\ \ \forall y_h \in \bar{\Omega}_h\}.
\end{align}
In other words, this is the collection of slopes of affine functions that satisfy the condition that defines convexity for a nodal function. Note that nodal functions are only defined on $\bar\Omega_h$. To extend a nodal functionto the domain $\dm$, we introduce its convex envelope.

\begin{definition}[convex envelope of a nodal function]\label{convexenvelope}
Let $w_h$ be a nodal function defined on $\bar{\Omega}_h$. The {\it convex envelope} of $w_h$ is the piecewise linear function 
\[
  \Gamma (w_h)(x) = \sup_{L} \left\{ L(x): \quad \mbox{$L$ affine function and } L(x_h) \leq w_h(x_h)\ \ \forall x_h\in \bar\Omega_h \right\}
\]
for any $x \in \dm$.
\end{definition}

We note that, by definition, $\Gamma(w_h)(x_h) \leq w_h(x_h)$ for any node $x_h \in \bar{\Omega}_h$, and equality holds for all interior nodes
if $w_h$ is convex.
Indeed, if $w_h$ is convex, by \eqref{def:convexnodalfunction}, for any node $x_h \in \bar{\Omega}_h$, there exists an affine function $L(x)$ satisfying
\[
L(y_h) \leq w_h(y_h) \quad \forall y_h \in \bar{\Omega}_h
\quad \mbox{and} \quad
L(x_h) = w_h(x_h).
\]
Since $L(x) \leq \Gamma (w_h)(x)$ for any $x \in \dm$ by Definition~\ref{convexenvelope}, we deduce that $w_h(x_h) = L(x_h) \leq \Gamma(w_h)(x_h) $. Combining this inequality with the inequality in the other direction, we have $w_h(x_h) = \Gamma (w_h)(x_h)$ for all interior nodes. 
Thus, $\Gamma (w_h)$ is a natural extension to $\bar\Omega$ of the convex nodal function $w_h$. 
With an abuse of notation, we still use $w_h$ to denote the convex envelope of this nodal function. 

The convex envelope of a nodal function $w_h$ induces a triangulation of the domain $\dm$ and a piecewise linear function over this triangulation. However, this triangulation is not known {\it a priori}.  Here we give an example to illustrate this property.

\begin{example}[convex envelope and triangulation]
\label{ex:FEMCE1}
Define the nodal set $\bar\Omega_h = \{z_1,\ldots,z_5\}$ with
$z_1 = (1,0),\ z_2 = (0,1),\ z_3 = (-1,0),\
z_4 = (0,-1)$, and $z_5 = (0,0)$.
Consider the nodal functions satisfying
\begin{alignat*}{2}
&w_1(z_1) = w_1(z_3)=1,\quad w_2(z_2)=w_2(z_4) = 1,\\
&w_3(z_1) = w_3(z_2)=w_3(z_3)=w_3(z_4)=1,
\end{alignat*}
and $w_j(z_i)=0$ otherwise.
The convex envelopes are $\Gamma (w_1) = |x_1|$, $\Gamma (w_2) = |x_2|$,
and $\Gamma(w_3) = |x_1|+|x_2|$.
The convex envelopes are subordinate to the meshes depicted in Figure~\ref{meshFig}.
\eremk\end{example}
%
\begin{figure}
\begin{center}
\begin{tikzpicture}[scale = 1.25]
\draw[-,thick](1,0)--(0,1)--(-1,0)--(0,-1)--(1,0);
\draw[-,thick](0,-1)--(0,1);
\node[inner sep = 0pt,minimum size=6pt,fill=black!100,circle] (n2) at (1,0)  {};
\node[inner sep = 0pt,minimum size=6pt,fill=black!100,circle] (n2) at (0,1)  {};
\node[inner sep = 0pt,minimum size=6pt,fill=black!100,circle] (n2) at (-1,0)  {};
\node[inner sep = 0pt,minimum size=6pt,fill=black!100,circle] (n2) at (0,-1)  {};
\node[inner sep = 0pt,minimum size=6pt,fill=black!100,circle] (n2) at (0,0)  {};
\node at (1.2,0) {1};
\node at (0,1.2) {0};
\node at (-1.2,0) {1};
\node at (0,-1.2) {0};
\node at (0.2,0) {0};
\end{tikzpicture}\quad
\begin{tikzpicture}[scale = 1.15]
\draw[-,thick](1,0)--(0,1)--(-1,0)--(0,-1)--(1,0);
\draw[-,thick](-1,0)--(1,0);
\node[inner sep = 0pt,minimum size=6pt,fill=black!100,circle] (n2) at (1,0)  {};
\node[inner sep = 0pt,minimum size=6pt,fill=black!100,circle] (n2) at (0,1)  {};
\node[inner sep = 0pt,minimum size=6pt,fill=black!100,circle] (n2) at (-1,0)  {};
\node[inner sep = 0pt,minimum size=6pt,fill=black!100,circle] (n2) at (0,-1)  {};
\node[inner sep = 0pt,minimum size=6pt,fill=black!100,circle] (n2) at (0,0)  {};
\node at (1.2,0) {0};
\node at (0,1.2) {1};
\node at (-1.2,0) {0};
\node at (0,-1.2) {1};
\node at (0,0.2) {0};
\end{tikzpicture}\quad
\begin{tikzpicture}[scale = 1.25]
\draw[-,thick](1,0)--(0,1)--(-1,0)--(0,-1)--(1,0);
\draw[-,thick](-1,0)--(1,0);
\draw[-,thick](0,-1)--(0,1);
\node[inner sep = 0pt,minimum size=6pt,fill=black!100,circle] (n2) at (1,0)  {};
\node[inner sep = 0pt,minimum size=6pt,fill=black!100,circle] (n2) at (0,1)  {};
\node[inner sep = 0pt,minimum size=6pt,fill=black!100,circle] (n2) at (-1,0)  {};
\node[inner sep = 0pt,minimum size=6pt,fill=black!100,circle] (n2) at (0,-1)  {};
\node[inner sep = 0pt,minimum size=6pt,fill=black!100,circle] (n2) at (0,0)  {};
\node at (1.2,0) {1};
\node at (0,1.2) {1};
\node at (-1.2,0) {1};
\node at (0,-1.2) {1};
\node at (0.2,0.2) {0};
\end{tikzpicture}
\end{center}
\caption{\label{meshFig} Meshes corresponding
to convex envelopes $\Gamma(w_1) = |x_1|$ (left),
$\Gamma(w_2)= |x_2|$ (middle),
and $\Gamma(w_3) = |x_1|+|x_2|$ (right).}
\end{figure}
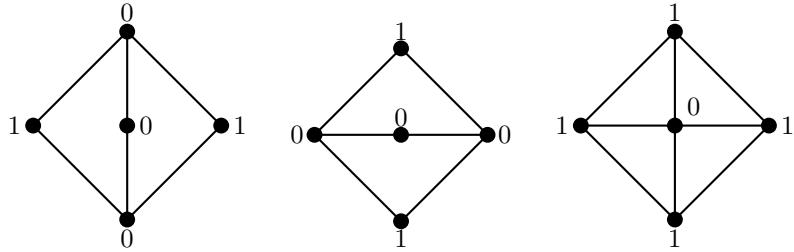

The above example shows that $\Gamma(w_h)$ is a piecewise linear function
that induces a mesh $\mct$ that depends on the values of $w_h$.
The example depicted in Figure \ref{fig:anisotropy} shows that,
if $w_h$ is the nodal interpolant of a function $w$, 
and if the Hessian $D^2 w$ is degenerate (or nearly degenerate), the induced mesh may be anisotropic.

\begin{figure}
\begin{center}
\begin{tikzpicture}[scale = 1.15]
\draw[-,thick](1,0)--(-1,1)--(-2,1)--(-1,0)--(1,-1)--(2,-1)--(1,0);
\draw[-,thick](1,0)--(0,0)--(-1, 0);
\draw[-,thick](-1,1)--(0,0)--(1, -1);
\draw[-,thick](-2,1)--(0,0)--(2, -1);
\node[inner sep = 0pt,minimum size=6pt,fill=black!100,circle] (n2) at (1,0)  {};
\node[inner sep = 0pt,minimum size=6pt,fill=black!100,circle] (n2) at (-1,1)  {};
\node[inner sep = 0pt,minimum size=6pt,fill=black!100,circle] (n2) at (-2,1)  {};
\node[inner sep = 0pt,minimum size=6pt,fill=black!100,circle] (n2) at (-1,0)  {};
\node[inner sep = 0pt,minimum size=6pt,fill=black!100,circle] (n2) at (1,-1)  {};
\node[inner sep = 0pt,minimum size=6pt,fill=black!100,circle] (n2) at (2,-1)  {};
\node at (1,0) [above right] {1};
\node at (-1,0) [below left] {1};
\node at (-1,1) [above right] {1};
\node at (1,-1) [below left] {1};
\node at (0,0) [above] {0};
\node at (-2,1) [above left] {0};
\node at (2,-1) [below right] {0};
\end{tikzpicture}
\caption{\label{fig:anisotropy} Mesh induced by the nodal interpolant  of $w(x) = (x \cdot e)^2$ where $e = (1,2)^\intercal$. Its convex envelope 
equals $|x \cdot e|$ in the star of $(0,0)$.}
\end{center}
\end{figure}
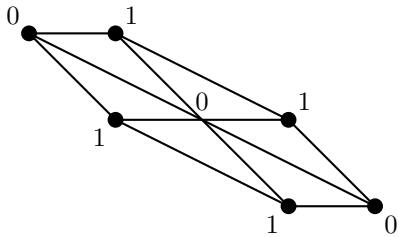


\subsubsection{The Oliker-Prussner method}
Now we are ready to introduce the Oliker-Prussner method \cite{Oliker88, NochettoZhangMA}.
We seek a convex nodal function $u_h$ satisfying the boundary condition
$u_h(x_h) = g(x_h)$ for all $x_h\in \p\Omega_h$ and
\begin{equation}
\label{OP}
|\p u_h(x_h)| =  \int_{\omega_{x_h}} f(x) \diff x,  \quad \forall x_h \in \Omega_h, 
\end{equation}
Note that, since  the partition $\{\omega_{x_h}\}_{x_h\in \Omega_h}$ is non-overlapping,
for all Borel sets $D\subset \Omega$, we have
\begin{align*}
  |\p u_h(D)| = \sum_{x_h\in D} f_{x_h}, \qquad \mbox{where } f_{x_h} = \int_{\omega_{x_h}} f(x) \diff x.
\end{align*}
Thus, the scheme is obtained
by replacing $f$ in \eqref{Alek} by a 
family of Dirac measures supported at the nodes in $\Omega_h$, 
and by replacing $g$
by its nodal interpolant on the boundary.
To implement the method, we need to derive a formula to compute the subdifferential of a nodal function $u_h$. 
This is a nontrivial task because it is non local. In fact, it involves computing the convex envelope of $u_h$. 
The following observation is useful in the characterization of the subdifferential. For a proof, see \cite{NochettoZhang16}.

\begin{lemma}[characterization of subdifferential]\label{char_subdifferential}
Let $w_h$ be a convex nodal function and $\Th$ be the mesh
induced by its convex envelope $\Gamma(w_h)$.
Then the subdifferential of $w_h$ at $x_h\in\Omega_h$ is the convex hull
of the constant gradients
$
\gradv \Gamma(w_h)|_T
$
for all $T \in \Th$ which contain $x_h$.
\end{lemma}

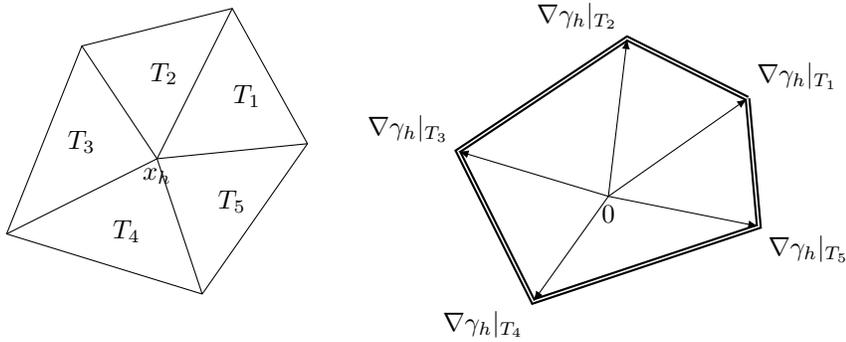
\begin{figure}[h!]
\begin{tikzpicture}
  \coordinate [label=below:$x_h$] (z1) at (1,1);
  \coordinate  (z2) at (3,1.2);
  \coordinate  (z3) at (2,3);
  \coordinate  (z4) at (0,2.5);
  \coordinate  (z5) at (-1,0);
  \coordinate  (z6) at (1.6,-0.8);
  \coordinate [label=below left:$T_1$] (k1) at (2.5, 2.1);
  \coordinate [label=below left:$T_2$] (k3) at (1.4,2.4);
  \coordinate [label=below left:$T_3$] (k4) at (0.3,1.5);
  \coordinate [label=below left:$T_4$] (k5) at (0.9,0.3);
  \coordinate [label=below left:$T_5$] (k2) at (2.3, 0.7);
  \coordinate [label=below :$0$] (o) at (7,0.5);
  \coordinate [label=above right:$\gradv \gamma_h|_{T_1}$] (g1) at (8.85,1.8);
  \coordinate [label=above left :$\gradv \gamma_h|_{T_2}$] (g2) at (7.25,2.6);
  \coordinate [label=above left:$\gradv \gamma_h|_{T_3}$] (g3) at  (5,1.1);
  \coordinate [label=below left:$\gradv \gamma_h|_{T_4}$] (g4) at  (6,-0.9);
  \coordinate [label=below right:$\gradv \gamma_h|_{T_5}$] (g5) at  (9,0.1);  
  \draw  (z1) -- (z2)  ;
  \draw  (z1) -- (z3)  ;
  \draw  (z1) -- (z4)  ;
  \draw  (z1) -- (z5)  ;
  \draw  (z1) -- (z6)  ;
  \draw  (z6)- - (z5) -- (z4) -- (z3) -- (z2) -- (z6);
  \draw [double, thick](g1) -- (g2) -- (g3) -- (g4) -- (g5) -- (g1);
  \draw [-latex] (o) -- (g1);
  \draw [-latex] (o) -- (g2);
  \draw [-latex] (o) -- (g3);
  \draw [-latex] (o) -- (g4);
  \draw [-latex] (o) -- (g5);
\end{tikzpicture}
\caption{Star centered at node $x_h$ corresponding to the mesh $\Th$
induced by the convex envelope $\gamma_h=\Gamma(w_h)$ 
and subdifferential $\partial w_h(x_h)$ of the convex
nodal function $w_h$ at node $x_h$. The latter is the convex hull of the constant
element gradients $\nabla \gamma_h|_{T_j}$ for $1\le j\le 5$.}
\label{fig:subdifferential}
\end{figure}

Figure \ref{fig:subdifferential} depicts the subdifferential
$\partial w_h(x_h)$ of a
convex nodal function $w_h$ at node $x_h$ for $d=2$.

\subsection{Stability, continuous dependence on data, and discrete maximum principle}
The Alexandrov estimate, which establishes the stability and continuous dependence of the \MA equation, is a cornerstone in the nonlinear PDE theory. In this subsection, we introduce a discrete version of the Alexandrov estimate suitable for nodal functions. We refer the reader to \cite{NochettoZhang16} for a complete proof.

\begin{lemma}[discrete Alexandrov estimate]\label{L:Alexandroff}
Let $w_h$ be a nodal function with $w_h(x_h) \geq 0$ at all $x_h \in \partial \Omega_h$. Then
\begin{align}\label{alex}
  \sup_{\dm_h} w_h^- \leq C \left( \sum_{x_h\in \mathcal{C}^-_h(w_h) } \abs {\partial w_h (x_h)} \right)^{1/d},
\end{align}
where $C = C(d, \Omega)$ is proportional to the diameter of $\Omega$
and $\mathcal{C}^-_h(w_h)$ is the (lower) contact set:
\begin{align}\label{contactset}
  \mathcal{C}^-_h(w_h) := \{ x_h \in \Nh, \quad \Gamma (w_h)(x_h) = w_h(x_h) \}.
\end{align}
\end{lemma}

The Alexandrov estimate establishes a lower bound for a nodal function in terms of the measure of the subdifferential at the (lower) contact set. 
Similarly, one can obtain an upper bound for a nodal function by the measure of the superdifferential at the (upper) contact set. 

Applying the discrete Alexandrov estimate, we are ready to compare two arbitrary nodal functions in terms of their subdifferentials. This is instrumental for the error analysis.

\begin{proposition}[stability of numerical solution]\label{stability}
Let $v_h$ and $w_h$ be two nodal functions with $v_h\geq w_h$ on $\p\Omega_h$.
Then
\[
  \sup_{\Omega_h} (v_h - w_h)^- \leq C \left( \sum_{x_h \in \mathcal{C}_h^-(v_h - w_h)} \Big( \abs{ \partial v_h (x_h)}^{1/d} - \abs{ \partial w_h (x_h)}^{1/d} \Big )^d  \right)^{1/d},
\]
where $C = C(d, \dm)$ is proportional to the diameter of $\dm$.
\end{proposition}
\begin{proof}
Let $v_h,w_h$ be two nodal functions. We introduce the convex envelope $\Gamma (v_h - w_h)$ as in Definition~\ref{convexenvelope}, and the nodal contact set $\mathcal{C}^-_h (v_h- w_h)$
defined in \eqref{contactset}. The discrete Alexandrov estimate of Lemma \ref{L:Alexandroff} yields 
  \begin{align}\label{stability-alexandroff}
    \sup_{\dm_h} ( v_h - w_h)^- \le C \left( \sum_{x_h \in \mathcal{C}^-_h(v_h - w_h)} |\partial  \Gamma (v_h - w_h)(x_h) | \right)^{1/d},
  \end{align}
whence we only need to estimate $|\partial  \Gamma (v_h - w_h)(x_h) |$
for all $x_h \in \mathcal{C}^-_h(v_h - w_h)$.
For these nodes, we easily see that
\[
\partial\Gamma(v_h - w_h)(x_h) \subset \partial(v_h - w_h)(x_h).
\]
We claim that 
\begin{align}\label{additioninequality}
\partial w_h(x_h) + \partial \Gamma (v_h - w_h)(x_h) \subset \partial v_h(x_h)
\quad\forall \, x_h \in \mathcal{C}^-_h(v_h - w_h).
\end{align}
Fix $x_h\in \calC^-_h(v_h-w_h)$, and let $\bp \in \partial w_h(x_h)$ and $\bq \in \partial \Gamma (v_h - w_h)(x_h)$, respectively, that is, by definition of the subdifferential \eqref{def:subdifferential},
\[
\bp \cdot (y_h - x_h) \leq w_h(y_h) - w_h(x_h)
\]
and
\[
  \bq \cdot (y_h - x_h) \leq   \Gamma (v_h - w_h)(y_h) -  \Gamma (v_h - w_h)(x_h) 
\]
for all nodes $y_h \in \Omega_h$. 
Adding both inequalities, we get
\[
(\bp+\bq) \cdot (y_h - x_h) \leq w_h(y_h) + \Gamma (v_h - w_h)(y_h) - \big( w_h(x_h) + \Gamma (v_h - w_h)(x_h) \big)
\]
Since $x_h$ is in the contact set $\mathcal{C}^-_h(v_h - w_h)$, we have $\Gamma (v_h - w_h)(x_h) = (v_h - w_h)(x_h)$. 
For all other nodes $y_h \in \Omega_h$, we have $\Gamma (v_h - w_h)(y_h) \leq (v_h - w_h)(y_h)$.
Hence, we deduce
\begin{align*}
(\bp+\bq) \cdot (y_h - x_h) \leq &\; w_h(y_h) + (v_h - w_h)(y_h) - \big( w_h(x_h) + (v_h - w_h)(x_h) \big) 
\\
= &\; v_h(y_h) - v_h(x_h). 
\end{align*}
This inequality implies $(\bp + \bq) \in \partial v_h(x_h)$ and proves the claim.

The Brunn-Minkowski inequality of Lemma~\ref{BM} applied to \eqref{additioninequality} yields
\begin{align*}
  |\partial w_h(x_h)|^{1/d} &+ |\partial \Gamma (v_h - w_h)(x_h)|^{1/d} 
  \\
  & \leq
  |\partial w_h(x_h) + \partial \Gamma (v_h - w_h)(x_h) |^{1/d} 
  \leq 
  | \partial v_h(x_h)|^{1/d} ,
\end{align*}
whence
\begin{align*}
  |\partial \Gamma (v_h - w_h)(x_h)| \leq &\; \left(| \partial v_h(x_h)|^{1/d} - |\partial w_h(x_h)|^{1/d}   \right)^{d} \quad \forall x_h\in \calC^-_h(v_h-w_h).
\end{align*}
This inequality gives us the desired estimate for $ |\partial \Gamma
(v_h - w_h)(x_h)|$. In view of \eqref{stability-alexandroff}, adding
over all $x_h \in \mathcal{C}^-_h (v_h - w_h)$ concludes the proof.
\end{proof}

A direct consequence of this stability result is the maximum principle
for nodal functions.

\begin{corollary}[discrete maximum principle]\label{MP}
Let $v_h$ and $w_h$ be two nodal functions over the nodal set $\bar\Omega_h$.
If $v_h(x_h) \geq w_h(x_h)$ at all $x_h \in \partial\Omega_h$ and 
$|\partial v_h(x_h)| \leq |\partial w_h(x_h)|$ at all $x_h \in \Omega_h$, then 
\[
w_h(x_h) \leq v_h(x_h) \quad \forall x_h \in \Omega_h.
\]
\end{corollary}

\begin{proof}
Since $v_h(y_h) \geq w_h(y_h)$ for all $y_h \in \p\Nh$, then for any node $x_h\in \mathcal{C}^-_h (v_h - w_h)$, we have
\[
\partial w_h(x_h) \subset \partial v_h(x_h).
\]
Combining this with the assumption that $|\partial v_h(x_h)| \leq |\partial w_h(x_h)|$ for all $x_h \in\Nh$, we
deduce $|\partial v_h(x_h)| = |\partial w_h(x_h)|$ for all
$x_h\in \mathcal{C}^-_h (v_h - w_h)$. Consequently, the stability of Proposition~\ref{stability} implies
\[
\sup_{\Omega_h} (v_h - w_h)^- = 0,
\]
whence $v_h - w_h \geq 0$. This completes the proof.
\end{proof}

Proposition~\ref{stability}  yields a lower bound on the difference between two nodal functions in terms of the difference of the measure of their subdifferentials. 
Similarly, to derive an upper bound, one may consider the functions $-w_h$ and $-v_h$ and derive 
\[
  \sup_{\Omega_h} (w_h - v_h)^- \leq C \left( \sum_{x_h \in \mathcal{C}_h^-(w_h - v_h)} \Big( \abs{ \partial w_h (x_h)}^{1/d} - \abs{ \partial v_h (x_h)}^{1/d} \Big )^d  \right)^{1/d}.
\]
Combining both bounds, we can derive a bound on $\norm{v_h - w_h}_{L^\infty(\Nh)}$ in terms of $|\p v_h(x_i)|$ and $|\p w_h(x_i)|$. 
In particular, the uniqueness of the solution of the Oliker-Prussner method follows immediately from Proposition \ref{stability}.

Finally, we notice that Proposition \ref{stability} is instrumental to derive error estimates. Define the nodal interpolation of a function $w$ as the nodal function $N_h w$ such that 
\begin{align}\label{def:nodalinterpolation}
  N_h w(x_h) = w(x_h) \quad \forall x_h\in \bar\Omega_h.
\end{align}
Setting $w_h = u_h$ and $v_h = N_h u$ In Proposition~\ref{stability}, where $u_h$ and $u$ solve \eqref{OP} and \eqref{eqn:MA}, respectively, we can derive an estimate for $\|{u_h - N_h u}\|_{L^\infty(\Omega)}$. 
It remains to estimate the discrepancy of the subdifferentials of the two nodal functions. 
While $|\p u_h(x_h)|=f_{x_h}$ is known by definition of the scheme \eqref{OP}, 
the measure of the subdifferential $|\p N_h u (x_h)|$ remains unknown. 
Therefore, the goal of our next step is to estimate the quantity
$|\p N_h u (x_h)|^{1/d} - f_{x_h}^{1/d}$ which will then 
be applied in Proposition \ref{stability} to derive a pointwise estimate.

%
\subsection{Consistency}\label{S:consistency}

In general, this method \eqref{OP} is consistent in the sense
that the right hand side of the \eqref{OP} can be written
equivalently as $\sum_{x_h \in \Nh} f_{x_h} \delta_{x_h}$ and this converges to $f$
in measure. 
However, such a concept of convergence is too weak to derive rates of convergence. 
Fortunately, we realize that if internal nodes are translation invariant,
then a reasonable notion of operator consistency holds for 
any convex quadratic polynomial; see Lemma \ref{L:consistency} below.
Such property is shown in \cite{BCM16,Mirebeau15,NochettoZhangMA}
for Cartesian nodes, see also section~\ref{sub:LBR}. In contrast, we give here an alternative
proof of consistency based on the geometric interpretation of
subdifferentials of convex quadratic polynomials
in the interior of the domain, extend the
results to $C^{2,\alpha}$ functions, and further investigate the
consistency error in the region close to the boundary.  To achieve this we, First, we require 
a definition.

\begin{definition}[adjacent set]\label{def:adjacentSet}
Given a convex nodal function $w_h$
and a node $x_h\in \Omega_h$, 
the {\it adjacent set} of $x_h$, denoted by $A_{x_h}(w_h)$,
is the collection of nodes $y_h\in \bar\Omega_h$ closest
to $x_h$ such that there exists
a supporting hyperplane $L$ of $w_h$ 
and $L(y_h)=w_h(y_h)$.
Thus, the set $A_{x_h}(w_h)$ is the collection
of nodes in the star associated with $x_h$ in the mesh $\mct$
induced by $\Gamma(w_h)$.
\end{definition}

\begin{lemma}[size of adjacent sets]\label{estimate_adjacent_set}
Let the nodal set $\Nh$ be translation invariant, and let $p$ be a $C^2$ convex function defined in $\bar\Omega$. If $ \lambda I
\leq D^2 p \leq \Lambda I$ in $\dm$ for some constants $\lambda, \Lambda>0$
and $p_h:=N_h p$ is the nodal function associated with $p$ defined in
\eqref{def:nodalinterpolation}, then 
the adjacent set of nodes $A_{x_h}(p_h)$ 
satisfies
\[
  A_{x_h}(p_h) \subset B_{Rh}(x_h) 
\]
where $R = \tfrac{\Lambda}{\lambda} d$,
and $B_{Rh}(x_h)$ is the ball centered at $x_h$ with radius $Rh$.
\end{lemma}
\begin{proof}
Let $z_h\in A_{x_h}(p_h)$ be 
such that 
\[
  |z_h - x_h| = \max \{ |y_h - x_h|: \, y_h \in A_{x_h}(p_h) \} .
\]
Without loss of generality, we may assume that  $p(x_h)=0$
and $\nab p(x_h) = 0$.
Let $\omega$ be the convex hull of the nodal set $\{ x_{\pm j} := x_h \pm h \tilde \be_j, \; j = 1, \ldots, d\}$ where $\{ \tilde \be_j\}_{j=1}^d$ is the basis defined in \eqref{translationinvariant}.
If $z_h \in \omega$, then the assertion is trivial because  $R \geq 1$. 

If $z_h\notin\omega$, then there is a constant $\tilde R\ge1$ such that $\tilde R^{-1} z_h \in \omega$, which implies that $|z_h| \leq \tilde R h$ and $|z_h| \geq \tilde R d^{-1/2} h$. Because $\omega$ is convex, we may write
\[
  \tilde R^{-1} z_h = \mathop{\sum_{j=1}^d}_{\sigma\in \{+,-\}} \alpha_{\sigma j} x_{\sigma j},
\qquad
\alpha_{\sigma j} \geq 0,
\qquad
  \mathop{\sum_{=1}^d}_{\sigma\in \{+,-\}} \alpha_{\sigma j} =1.
\]
We next note that $p(x_{\pm j}) \leq \tfrac 12 \Lambda h^2$ for all $j = 1, \ldots, d$ because $D^2 p \leq \Lambda I$, $|x_{\pm j} - x_h| \leq h$, $p(x_h)=0$ and $\nab p(x_h)=0$. Since $z_h
\in A_{x_h}(p_h)$, there exists a supporting hyperplane $L$ at $x_h$
such that
\[
L(z_h) = p_h(z_h),
\qquad
L(x_{\pm j}) \leq p_h(x_{\pm j}) \leq \frac 12 \Lambda h^2.
\]
Exploiting that $L$ is linear and $p_h(x_h) =0$ yields
\[
 p_h(z_h) =  L(z_h) = L\left( \tilde R \mathop{\sum_{j=1}^d}_{\sigma\in \{+,-\}} \alpha_{\sigma j} x_{\sigma j}\right)
 = \tilde R \mathop{\sum_{j=1}^d}_{\sigma\in \{+,-\}} \alpha_{\sigma j} L(x_{\sigma j}) \leq \frac 12 \Lambda h^2 \tilde R.
\]
On the other hand, since $D^2 p \geq \lambda I$ and $|z_h| \geq \tilde R h d^{-1/2} $, we have
\[
p_h(z_h) = p(z_h) \geq \frac {\lambda}2 |z_h|^2 \geq \frac {\lambda}{2} \tilde R^2 d^{-1} h^2.
\]
Combining the last two inequalities implies
\[
  \tilde R \leq R = \frac{\Lambda}{\lambda}d.  
\]
This completes the proof.
\end{proof}

The previous result shows that for any node $x_h$ with  $\dist(x_h, \bdry) > Rh$, all nodes in its adjacent set are contained in $\Omega_h$.
We apply this observation to establish the following consistency result.  

\begin{lemma}[properties of convex interpolation]\label{properties}
Let $p$ be a convex quadratic polynomial such that $\lambda I \leq D^2p \leq \Lambda I$, and let $p_h=N_h p$ be the nodal function 
defined by \eqref{def:nodalinterpolation}. Then the following properties hold:
\begin{enumerate}
  \item For all $x_h\in \Omega_h$ we have $\partial p_h(x_h) \neq \emptyset$.

  \item If the nodal set $\Nh$ is translation invariant and $\dist(x_h, \bdry_h) \geq R h$, with $R =  \tfrac{\Lambda}{\lambda} d$, under a uniform refinement from $\Nh$ to $\dm_{h/2}$, we have 
  \[
    \abs {\partial p_h(x_h)} = 2^d \abs{\partial p_{h/2}(x_h)}.
  \]
  
  \item If the nodal set $\Nh$ is translation invariant, $\dist(x_h, \bdry_hh) \geq Rh$, and $\dist(y_h, \bdry_h) \geq R h$, then
  \[
    \abs{\partial p_h(x_h)} = \abs{\partial p_h (y_h) }.
  \]
\end{enumerate}
\end{lemma}
\begin{proof}
To prove the first claim, we only need to note that if $\ell$ is the tangent plane of $p$ at $x_h$, then $\ell$ is a supporting plane of $p_h$ at $x_h$. Thus $\nab \ell \in \p p_h(x_h)$.

To prove the second claim, we may assume that $p(x_h) =0$, and $\gradv p(x_h) = 0$. 
Note that for homogeneous quadratic polynomials, we have
\[
  p(x) = 4p\left(\frac x 2\right).
\]
A simple calculation yields 
\[
  \partial p_h(x_h) = 2 \partial p_{ h/ 2} (x_h)  
\]
and therefore $\abs {\partial p_h(x_h)} = 2^d \abs{\partial p_{ h/2} (x_h)}$.

To prove the third claim. 
We consider the function 
\[
p^{x_h} (x) = p(x) - \gradv p(x_h) \cdot (x - x_h) - p(x_h),
\]
obtained by subtracting the tangent plane of $p$ at $x_h$. 
Since adding an affine function does not change the measure of the subdifferential, we have $|\partial p_h(x_h)| = |\partial p^{x_h}_h(x_h)| $. Further note that by subtracting the tangent plane at a node $y_h$, we obtain the same function up to a parallel translation, that is,
\[
p^{x_h}(x - x_h) = p^{y_h}(x - y_h). 
\]
Since the mesh is translation invariant, we have that if $L$ is a supporting plane of $p_h^{x_h}$ at $x_h$, then by a parallel translation it is also a supporting plane of $p_h^{y_h}$ at $y_h$. Hence, we have $|\partial p^{x_h}_h(x_h)| = |\partial p^{y_h}_h(y_h)|$. Since   $|\partial p_h(x_h)| = |\partial p^{x_h}_h(x_h)| $ for all nodes $x_h$, we conclude that $|\partial p_h(x_h)| = |\partial p_h(y_h)|$.
\end{proof}

Now we are ready to prove the consistency, for a proof see \cite[Lemma 5.3]{NochettoZhangMA}.

\begin{lemma}[consistency I]\label{L:consistency}
Let $p$ be a convex quadratic polynomial such that $\lambda I \leq
D^2p \leq \Lambda I$, and let $p_h:=N_h p$ be the corresponding convex nodal
function defined in \eqref{def:nodalinterpolation}.
Let $\Omega_h$ be translation invariant.
Then
\[
  | \partial p_h(x_h) | = \int_{\omega_{x_h}} \det D^2 p(x)\diff x 
\]
for any node $x_h \in \Omega_h$ such that $\dist(x_h,\bdry) \geq Rh$
with $ R= \tfrac{\Lambda}{\lambda}d$. 
\end{lemma}
\begin{proof}
Let $\phi$ be any continuous function with compact support in $\dm$. We consider a sequence of nested refinements $\Omega_{h_n}$ with $h_n = 2^{-n} H$, for a fixed $H>0$.

By Lemma~\ref{weakconvergence} we immediately obtain, as $n \to \infty$, that
\[
  \sum_{y_{h_n} \in \dm_{h_n}} \phi(y_{h_n}) \abs{\partial p_{h_n}(y_{h_n})} \rightarrow  \int_{\dm} \phi \det D^2 p(x)\diff x  = \det D^2 p(x) \int_{\dm} \phi(x)\diff x .
\]
Thus, we only need to prove that as $n \to \infty$
\[
    \sum_{y_{h_n} \in \dm} \phi (y_{h_n}) \abs{\partial p_{h_n}(y_{h_n})}  \rightarrow \frac{ \abs{\partial p_H(x_H)} }{\abs{\omega_{x_H}}} \int_{\dm} \phi(x)\diff x.
\]
In view of second and third result in Lemma~\ref{properties}, we have
\[
  \abs{\partial p_{h_n}(y_{h_n})} = \abs{\partial p_{h_n}(x_H)} =2^{-nd} \abs{\partial p_H(x_H)}.
\]
The refinement strategy implies that $\abs{\omega_{y_{h_n}}} = 2^{-nd}\abs{\omega_{x_H}}$. Thus, we infer that
\begin{align*}
 \sum_{y_{h_n} \in \dm_{h_n}} \phi(y_{h_n}) \abs{\partial p_{h_n}(y_{h_n})} 
 &=  \frac {\abs{\partial p_H(x_H)}}{\abs{\omega_{x_H}}} \sum_{y_{h_n} \in \Omega_{h_n}}  \phi(y_{h_n}) \abs{\omega_{y_{h_n}} }\\
&  \rightarrow
  \frac{\abs{\partial p_H(x_H)} }{\abs{\omega_{x_H}}} \int_{\dm} \phi(x)\diff x.
\end{align*}
This completes the proof.
\end{proof}

Moreover, for convex cubic polynomials, we have the following consistency error estimate. {This result, to our knowledge, has not appeared elsewhere}.

\begin{lemma}[consistency II]
\label{lem:C^4}
  Let $x_h \in \Omega_h$ and $q$ be a convex cubic polynomial such that $ \lambda I \leq D^2 q \leq \Lambda I$ in the ball $B_{Rh} := \overline{B_{Rh}(x_h)} \subset \dm$, with $R = \tfrac\Lambda\lambda d$. Then 
\begin{align*}
  \left| |\p N_h q(x_h)| - \int_{\omega_{x_h}} \det D^2 q (x)\,dx \right| \leq C h^{d+2} |q|_{C^{3}(B_{Rh})}..
\end{align*}
\end{lemma}
\begin{proof}
Without loss of generality, we may assume that $x_h = 0$ and $q(0) = 0$ and $\gradv q(0) = 0$.
We decompose the cubic polynomial $q(x)$ as 
\[
q(x) = p(x) + h r(x),
\]
where $p(x)$ is a quadratic polynomial such that $D^2 p = D^2 q(0)$ and $r(x)$ is a homogeneous cubic polynomial.
Since, by Lemma \ref{estimate_adjacent_set}, the adjacent set $A_{x_h}(q)$ of the node $x_h =0$ is contained in a ball of radius $Rh$ we deduce that 
\[
  |p(z_h)| \leq C_q R^2 h^2, \quad \quad  |r(z_h)| \leq C_r R^2 h^2 \quad \forall z_h \in A_{x_h}(q),
\]
where $C_q$ and $C_r$ depends on $D^2 p$ and $D^3 r$, respectively.
We set 
\[
q_t(x) = p(x) + t r(x) \quad t \in [-h, h],
\]
and note that $\lambda I \leq D^2 q_t(0) \leq \Lambda I$ for all $t$.
Therefore, the adjacent set of $q_t$ at $0$ remains in the ball $B_{Rh}$. 

We set the measure of its subdifferential of $q_t$ at $x_h$ as a function of $t$
\[
  m(t) = |\p N_h q_t(x_h)| = |\p N_hq_t(0)|,
\]
and note that we aim to show that
\begin{align*}
  \left| m(h)-\int_{\omega_{x_h}} \det D^2 q(x) \diff x\right|\le C h^{d+2} |q|_{C^3(B_{Rh})}.
\end{align*}
Now we proceed to prove the lemma in the following steps.

\begin{enumerate}
  \item We aim to show that $m(t)$ is a polynomial of degree $d$  
  \begin{align}\label{mtpolynomial}
    m(t) = \sum_{k=0}^d C_k t^k. 
  \end{align}
  and the coefficients $C_k$ satisfy $|C_k| \leq C h^d$ where $C$ depends on $|D^2 p|$, $|D^3 r|$, and the dimension $d$. 
  By the characterization of the subdifferential, given in Lemma \ref{char_subdifferential}, the subdifferential of $N_h q_t$ at $0$ is the convex hull of the piecewise gradient of its convex envelope $\gradv \Gamma(N_h q_t)|_T$ for all $T\in \mct$ that have $x_h$ as a vertex; see Figure \ref{fig:subdifferential}. We label these simplices as $T_1, \cdots, T_N$ and, to simplify notation,   we set the piecewise gradient of $\Gamma (N_h p)$ and $\Gamma (N_h r)$ at $T_i$ as
  \[
    \bv_i = \gradv \Gamma (N_h p)|_{T_i}, 
    \qquad
    \bw_i = \gradv \Gamma (N_h r)|_{T_i}, \qquad i = 1, \ldots, N.
  \]
  Hence, we have 
  \[
    \bu_i:= \gradv \Gamma (N_h q_t)|_{K_i} = \bv_i + t \bw_i.
  \]
  To compute the measure of the convex hull of $\{\bu_i \}$, we may divide the convex hull into a set of disjoint simplices $\{S_i, \; i= 1, \cdots, N\}$ and label the vertices of $S_i$ as  $\{{\bm 0}, \bu_{i_1}, \cdots, \bu_{i_d}  \}$. 
  Thus, we obtain 
  \[
    m(t) = \sum_{i=1}^N |S_i|  \quad \mbox{where $|S_i|$ denotes the signed volume of $S_i$. }
  \]
  and so, by the volume formula of simplices, we get
  \begin{align}\label{mtvolume}
    m(t) = \frac 1{d!} \sum_{i=1}^{N} \det 
    \begin{pmatrix}
      1 & {\bm 0}^\intercal\\
      1 & \bu_{i_1}^\intercal\\
      \vdots & \vdots\\
      1 & \bu_{i_d}^\intercal
    \end{pmatrix}.\
    \qquad
    \bu_{i_j} = \bv_{i_j} + t \bw_{i_j}.
    \end{align}
  Now, it is clear that
  \[
    |S_i| = \sum_{k = 0}^d C_k^i t^k
  \]
  is a polynomial of $t$ with degree at most $d$. Thus, $m(t)$ must be a polynomial with degree at most $d$ as well. Furthermore, by the volume formula of simplices \eqref{mtvolume}, the coefficients $|C_k^i| \leq C h^d$ because both $|\vf v_{i_j}| \leq C h$ and $|\vf w_{i_j}| \leq C h$. Finally, the number $N$ of simplices $S_i$ is finite and bounded by the number of vertices in the adjacent set $A$. 

  \item We show that $m'(0) = 0$.  To do so, it suffices to show that the function $m$ is even, that is $m(t) = m(-t)$ for all $-h \leq t \leq h$. Note that if $\bv \in \p N_h (p + t r)(0)$, then $-\bv \in \p N_h (p - t r)(0)$ for any $t \in (0, h]$. Indeed, since the subdifferential set is determined by the function values on the adjacent set which is contained in the ball $B_{Rh}(0)$, if 
  $
  \bv \cdot y_h \leq (p + t r)(y_h) 
  $
  for all $y_h \in B_{Rh}(0)$,
  then 
  \[
    \bv \cdot (-y_h) \leq (p + t r)(-y_h) \quad \forall y_h \in B_{Rh}(0).
  \]
  Hence, $-\bv \in \p N_h (p - t r)(0)$ because $p(y_h) = p(-y_h)$.
  Thanks to this symmetry property, we deduce that $ |\p N_h (p - t r)(0)| =  |\p N_h (p + t r)(0)|$, i.e., $m(t) = m(-t)$.
  
  \item We show that
  \[
    |m(h) - m(0)| \leq C h^{d+2}. 
  \]
  Combining the previous two steps we get that 
  \[
    m(t) = m(0) + C_2 t^2 + \dots + C_d t^d
  \]
  because $C_1 = m'(0) = 0$.  Since $|C_j| \leq C h^d$ for $j = 2, \dots, d$, we deduce that 
  \[
    |m(t) - m(0)| \leq C h^{d+2} \quad \forall t \in [0, h].
  \]

  \item It remains to show that 
  \[
    \left| \int_{\omega_{x_h}} \det D^2q(x) \diff x - m(0) \right| \leq C h^{d+2}.
  \]
  By  the consistency for quadratics given in Lemma~\ref{L:consistency}, we have $m(0) = \int_{\omega_{x_h}} \det D^2p(x) \diff x$. Therefore, it is sufficient to show that
  \[
    \left| \int_{\omega_{x_h}} \left( \det D^2 q(x) - \det D^2 p(x) \right)\diff x \right| 
    \leq C h^{d+2}.
  \]
  A Taylor expansion of  $\det D^2q  = \det D^2 (p + h r)$ reveals that
  \[
    \left| \det D^2q(x) - \det D^2 p (x) -  h \cof D^2 p(x) : D^2 r(x) \right| \leq C h^2.
  \]
  where the constant $C$ depends on $D^2p$ and $D^3 r$. This implies that 
  \begin{align*}
   & \left| \int_{\omega_{x_h}} \left( \det D^2 q(x) - \det D^2 p(x) \right)\diff x \right| \\
   &\qquad\qquad \leq 
    h \left| \int_{\omega_{x_h}} \cof D^2 p(x) : D^2 r(x) \diff x \right|
    + C h^2 \abs{ \omega_{x_h} }.
  \end{align*}
  Noting that $\cof  D^2 p :D^2 r $ is an odd function and $\omega_{x_h}$ is symmetric respect to the origin, we obtain
  \[
     \int_{\omega_{x_h}} \cof D^2 p(x) :D^2 r(x)\diff x = 0
  \]
  and 
  \[
    \left| \int_{\omega_{x_h}} \left( \det D^2 q(x) - \det D^2 p(x) \right)\diff x \right| \leq C h^{d+2}.
  \]
\end{enumerate}
This completes the proof.
\end{proof}

Now for any function $w$ that can be approximated locally by a quadratic polynomial 
such that $w(x) = p(x) + \mathcal{O}(h^{2+\alpha})$ in $B_{Rh}(x_h)$ or by a cubic polynomial such that $w(x) = q(x) + O(h^{3+\alpha})$ in $B_{Rh}(x_h)$, we show that the consistency error of the Oliker-Prussner method is of order $\mathcal{O}(h^{\alpha})$ and $\mathcal{O}(h^{1+\alpha})$, respectively.

\begin{proposition}[interior consistency]
\label{P:Holder-interior}
Let $\Omega_h$ be a translation invariant set of nodes, and $x_h\in \Omega_h$ be such that $\dist(x_h, \bdry_h) \geq Rh$ with $R =\tfrac{\Lambda}{\lambda}d$. If $w \in C^{2+k, \alpha}(\overline{B_{Rh}})$, with $k\in \{0, 1\}$, and $\alpha \in (0, 1]$) is a convex function with $\lambda I \leq D^2 w \leq \Lambda I$, then we have
\[
\left| \abs{\partial N_h w (x_h)} - \int_{\omega_{x_h}} \det D^2 w(x)\diff x \right| \leq C h^{k + \alpha} |w|_{C^{2+k,\alpha}(\overline{B_{Rh}})}\abs{\omega_{x_h}},
\]
where $C=C(d,\lambda,\Lambda)$.
\end{proposition}
\begin{proof}
We divide the proof into two cases $k =0$ and $k =1$.
\begin{enumerate}
  \item {\bf Case $k=0$:} We only need to show the inequality
  \begin{align*}
    | \partial N_h w (x_h) | 
    \leq \int_{\omega_{x_h}} \det D^2 w(x)\diff x + Ch^{\alpha} |w|_{C^{2,\alpha}(\overline{B_{Rh}})} \abs{\omega_{x_h}},
    \end{align*}
  because the reverse inequality can be derived similarly. Since $w \in C^{2, \alpha}(\overline{B_{Rh}})$, we estimate $w$ by a quadratic polynomial $p$ so that 
  \[
    w(x) \leq  p(x) \quad \forall x \in B_{Rh}(x_h),
  \]
  where $p(x_h) = w(x_h)$, $\gradv p(x_h) = \gradv w(x_h)$ and 
  \[
    D^2 p = D^2 w(x_h) + C h^\alpha |w|_{C^{2,\alpha}(\overline{B_{Rh}})} I
  \]
  for a fixed, and sufficiently large, constant $C$. Let $p_h=N_h p$, and note that
  \[
    |\partial N_h w(x_h)| \leq |\partial p_h(x_h)|.
  \]
  It remains to show that
  \[
    |\partial p_h(x_h)| \leq \int_{\omega_{x_h}} \det D^2 w(x)\diff x + Ch^{\alpha} |w|_{C^{2,\alpha}(\overline{B_{Rh}})}\abs{\omega_{x_h}}.
  \]
  Since $(\lambda + Ch^{\alpha})I \leq D^2 p \leq (\Lambda + Ch^{\alpha})I$ and
  \[
    \frac{\Lambda + Ch^{\alpha}} {\lambda + Ch^{\alpha}} \leq  \frac {\Lambda}{\lambda} \quad \text{because $\Lambda \geq \lambda$},
  \]
  invoking the consistency of Lemma~\ref{L:consistency} we obtain
  \[
    | \partial p_h(x_h) | = \int_{\omega_{x_h}} \det D^2 p(x)\diff x 
  \]
  provided that dist$(x_h,\partial\dm_h) \ge Rh$. Recalling that $w\in C^{2,\alpha}(\overline{B_{Rh}})$, we can write $D^2 p = D^2 w(x) + E(x)$ for all $x\in \overline{B_{Rh}}$, where $|E(x)| \le C h^\alpha |w|_{C^{2,\alpha}(\overline{B_{Rh}})}$. A Taylor expansion yields
  \[
    | \partial p_h(x_h) | \leq \int_{\omega_{x_h}} \det D^2 w(x)\diff x + Ch^{\alpha} |w|_{C^{2,\alpha}(\overline{B_{Rh}})} \abs{\omega_{x_h}}.
  \]
  
  \item {\bf Case $k=1$:} If $w \in C^{3, \alpha}(\overline{B}_{Rh})$, we approximate $w$ by a cubic polynomial $q$ so that 
  \[
    w(x) \leq  q(x) \quad \forall x \in B_{Rh}(x_h),
  \]
  where $q(x_h) = w(x_h)$, $ \gradv q(x_h) = \gradv w(x_h)$,
  \[
    D^2 q(x_h) = D^2 w(x_h)  + C h^{1+\alpha} |w|_{C^{3,\alpha}(\overline{B_{Rh}})},
  \]
and $D^3 q = D^3 w(x_h)$ with universal constant $C$. The rest of the proof is similar to the previous case. 
\end{enumerate}
Combing both cases, we conclude the proof of the estimate.
\end{proof}

\subsection{Pointwise error estimate}
We are now ready to show a pointwise error estimate for the method \eqref{OP} under suitable regularity assumptions on the solution $u$. We aim to apply the stability of the numerical scheme shown in Proposition~\ref{stability} to derive a lower bound of the difference $v_h - u_h$, for a suitable convex piecewise linear function $v_h$.

Assume that the convex solution $u$ of the \MA equation \eqref{eqn:MA} is $C^{k,\alpha}$ near the boundary of the domain $\dm$ where $k\in \{2,3\}$ and $\alpha\in (0,1]$.
We first extend the solution to a larger convex domain 
\[
  \dm_{4Rh} = \{ x \in \mathbb R^d, \dist(x, \dm) \leq 4Rh \}
\]
such that, for sufficiently small $h$, the extended function, which we still denote as $u$, remains $C^{k,\alpha}$-continuous in the extended region and satisfies 
\begin{align}\label{bound_of_u}
  \frac \lambda 2 I  \leq D^2 u(x) \leq 2 \Lambda I \quad \mbox{for any $x \in \dm_{4Rh}$.}
\end{align}
Next, we extend the translation invariant interior nodal set $\dm_h$ to the extended domain $\dm_{4 Rh}$ and, by an abuse of notation, we still denote the set as $\Omega_h$, that is, 
\[
  \Omega_h = \left\{ x_h = \sum_{j=1}^d z^j \tilde \be_j: \quad z^j \in \mathbb Z \right\}\cap \Omega_{4Rh}.
\]
We construct the piecewise linear function $v_h = \Gamma (N_h u)$ by taking the convex envelope of the nodal interpolation of the solution $u$ on $\Omega_h $ in the extended domain and then restrict the piecewise linear function $v_h$ to the domain $\Omega$. Thus, this procedure yields a piecewise linear function $v_h$ defined on the domain $\dm$.

We claim that the piecewise linear function $v_h$ satisfies the following two conditions which are useful in the error estimate.
First, the adjacent set size estimate of Lemma \ref{estimate_adjacent_set} and the bound of $D^2u$ given in \eqref{bound_of_u} imply that for any interior node $x_h \in \Omega_h \cap \dm$, its adjacent set $A_{x_h}(v_h)$ is contained in the extended 
domain $\dm_{4Rh}$.
Second, we notice that $|v_h(x) - u(x)| \leq C h^2$ on the boundary $\bdry$ where the constant $C$ depends on $\norm{u}_{C^2(\dm)}$. 
This is simply due to the fact that the diameter of any patch of a node $z \in \Omega_h \cap \dm$ is bounded by $4Rh$ 
and interpolation theory of piecewise linear function.

Now we are ready to derive the main error estimate. 

\begin{theorem}[error estimate]
\label{thm:OPMainEst}
Let $u$  be the solution of the \MA equation \eqref{eqn:MA}, $0 \leq \lambda I \leq D^2 u \leq \Lambda I$ and $u \in C^{2+k, \alpha}(\bar{\dm})$ with $k\in \{0,1\}$
and $\alpha\in (0,1]$.
Let $\Nh$ be a translation invariant nodal set satisfying \eqref{translationinvariant}, and let $u_h$ be the solution of discrete \MA equation \eqref{OP} defined on $\Nh$. Then we have 
\[
  \norm{u - u_h}_{L^\infty{(\dm)}} \leq C h^{k + \alpha},
\]
where the constant $C$ depends only on $\norm{u}_{C^{2+k, \alpha}(\dm)}$, $\lambda$, $\Lambda$, $\textrm{diam}(\dm)$, and space dimension $d$.
\end{theorem}
\begin{proof}
Let $v_h$ be the interpolation of the extension of the solution $u$ defined above. 
Since $|v_h - u_h| \leq C h^2$ on the boundary $\bdry$, we have $v_h + C h^2 \geq u_h$. 
By the stability of the numerical solution, Lemma \eqref{stability}, we obtain 
\[
  \sup_{\Omega_h} (v_h + Ch^2 - u_h)^- \leq C  \left( \sum_{x_i \in \mathcal{C}^-_h(v_h - u_h)}  ( |\p v_h (x_i)|^{1/d} - |\p u_h (x_i)|^{1/d})^d \right)^{1/d}.
\]
Invoking the consistency error estimate, Proposition \ref{P:Holder-interior}, we immediately obtain 
\[
  \sup_{\Omega_h} (v_h + Ch^2 - u_h)^- \leq C h^{k + \alpha}.
\]
By a simple algebraic manipulation, the estimate yields a lower bound for the error $v_h - u_h \geq - C h^2 - C h^{k+\alpha}$. Similarly, an estimate for the upper bound follows by considering the function $u_h + Ch^2 - v_h$. Combining both estimates, we get the desired result.
\end{proof}

\subsection{$W^{2,p}$ error estimate}
The results and arguments of the previous section have recently been
extended to the derivation of $W^{2,p}$ error estimates
of the Oliker-Prussner scheme \cite{NeilanZhangSINUM}.  Here, the discrete
$W^{2,p}$ norm is taken to be the
sum of weighted second-order differences:
\begin{align*}
\|v\|_{W^{2,p}_f} = \left(\sum_{x_h\in \Omega_h} f_{x_h} |\Delta_{ \be} v(x_h)|^p\right)^{1/p}.
\end{align*}

The starting point is a simple observation that 
the contact set of a nodal function contains
information of its second order difference.
In particular, if $u_h$ is the solution to \eqref{OP}
and $v_h$ is some approximation to $u$,
then we can define the perturbed error
\begin{equation}\label{eqn:wepsilonDefOP}
w^\epsilon_h = v_h-(1-\epsilon)u_h
\end{equation}
with parameter $\epsilon\in (0,1)$. Now,
by using the identity $\Delta_{\be} w^\epsilon_h(x_h)\ge \Delta_{\be} \Gamma w_h^\epsilon(x_h)\ge 0$
for $x_h\in \calC_h^-(w_h^\epsilon)$, we have, after some algebraic manipulations,
\begin{align*}
\Delta_e (u_h-v_h)(x_h)\le \frac{\epsilon}{1-\epsilon}\Delta_e v_h(x_h)\qquad \forall x_h\in \calC_h^-(w_h^\epsilon).
\end{align*}
The right-hand side of this expression is uniformly bounded for appropriate $v_h$ if $u$ is sufficiently smooth,
and therefore we find that the error $\Delta_e (u_h-v_h)(x_h)$ is controlled on the contact set $\calC_h^-(w_h^\epsilon)$. 
However, noting that $w_h^\epsilon$ is not necessary convex, we must estimate 
$\Delta_e(u_h-v_h)(x_h)$ on the complement set 
\begin{equation}\label{eqn:EepsilonDefOP}
E^\epsilon:=\Omega_h\backslash \calC_h^-(w_h^\epsilon).
\end{equation}
This is done by estimating its cardinality in terms of the consistency of the method.
%
\begin{lemma}[size of complement set]
\label{lem:SizeOfCoomplement}
Let $u_h$ and $v_h$ be convex nodal functions with
$u_h = v_h$ on $\p \Omega_h$ and $u_h\le v_h$
on $\Omega_h$.  Set 
\[
|\p u_h(x_h)| = f_{x_h}\quad \text{and}\quad |\p v_h(x_h)| = g_{x_h}\qquad x_h\in \Omega_h.
\]
Then there exists a constant $C>0$ depending only on $f$ such that
\begin{align*}
\sum_{x_h\in E^\epsilon} 
f_{x_h}\le C \frac{(1-\epsilon)}{\epsilon} \|f^{1/d}-g^{1/d}\|_{\ell^d(\calC_h^-(w^\epsilon_h))},
\end{align*}
where $w^\epsilon_h$ and $E^\epsilon$ are defined by \eqref{eqn:wepsilonDefOP}
and \eqref{eqn:EepsilonDefOP}, respectively.
\end{lemma}

The last ingredient to develop $W^{2,p}$ estimates
is a simple result of the discrete $L^1$ norm of a nodal
function in terms of its level sets.  Roughly speaking
this result gives a relation between Riemann and Lebesgue 
sums; see \cite[Lemma 5.1]{NeilanZhangSINUM}
\begin{lemma}\label{lem:PrelimQ}
Let $s_h$ be a  nodal function with $|s_h(x_h)|\le M$ for some $M>0$.
Then, for any $\sigma>0$,
\begin{align*}
  \sum_{x_h\in \Omega_h} f_{x_h}|s_h(x_h)|\le \sigma \sum_{k=0}^{{M}} \sum_{x_h\in A_k} f_{x_h},
\end{align*}
where
\[
A_k = \{x_h\in \Omega_h:\ |s_h(x_h)|\ge k \sigma\}.
\]
\end{lemma}

\begin{theorem}[$W^{2,p}$ error estimate]
\label{thm:OPW2pResult}
Suppose that the conditions of Theorem 
\ref{thm:OPW2pResult} are satisfied with $k+\alpha = 2$.
Then there holds
\begin{align*}
\|u-u_h\|_{W^{2,p}_f} \le \left\{
\begin{array}{ll}
Ch^{1/p} & p\in (d,\infty)\\
C |\log h|^{1/d} h^{1/d} & p\in (1,d].
\end{array}
\right.
\end{align*}
\end{theorem}

We now give a sketch of the main ideas to prove Theorem 
\ref{thm:OPW2pResult} and refer the reader to \cite{NeilanZhangSINUM} for details.  
To communicate
the main ideas, we make the simplifying assumption
that the consistency estimate in Proposition \ref{P:Holder-interior}
holds up to the boundary.  We also assume 
homogeneous boundary conditions, i.e., $g=0$ in \eqref{eqn:MA2}.  
These assumptions, which
do not hold in general, allow us to derive better rates
of convergence than those stated in Theorem 
\ref{thm:OPW2pResult}.

As a first step we set $v_h = (1-C h^{2})^{1/d} N_h u$,
where $C>0$ is sufficiently large such that (cf.~Proposition \ref{P:Holder-interior})
\begin{align*}
g_{x_h} = |\p v_h(x_h)| = (1-C h^{2})|\p N_h u(x_h)|\le f_{x_h}.
\end{align*}
Therefore by the comparison principle in Corollary \ref{MP}, we have $v_h\ge u_h$
on $\Omega_h$.  We also have $|f_{x_h}-g_{x_h}|\le C h^{2+d}$.

To deduce the estimate, it suffices bound
\begin{align*}
\sum_{x_h\in \Omega_h} f_{x_h} \big(\Delta_{\be} (u_h-v_h)(x_h)\big)^+.
\end{align*}
Bounding the negative part of the error can be obtained by similar arguments.

For parameter $\epsilon_k$ with $\epsilon_k/(1-\epsilon_k) = C k^{1/p} h^{2}$,
we define
\[
A_k = \left\{x_h\in \Omega_h:\ \Delta_{\be}(u_h-v_h)(x_h)\ge \frac {\epsilon_k} {1-\epsilon_k} {\Delta_{\be} v_h(x_h)} \right\},
\]
and note that $ A_k \subset E^{\epsilon_k}$.  
Let $s_h(x_h) = \left|\big(\Delta_{\be} (u_h-v_h)\big)^+\right|^p$, and note that $|s_h(x_h)| \leq C h^{-2p}$ because $u_h$ and $v_h$ are bounded.
Applying Lemma
\ref{lem:PrelimQ}, with $\sigma = C h^{2p}$,
we have
\begin{align*}
\sum_{x_h\in \Omega_h} f_{x_h} \left| \big(\Delta_{\be} (u_h-v_h)(x_h)\big)^+\right|^p\le C h^{2p} \left(1+ \sum_{k=1}^{Ch^{-2p}} \sum_{x_h\in A_k} f_{x_h}\right).
\end{align*}
On the other had, using Lemma \ref{lem:SizeOfCoomplement} and the consistency of the scheme yields, for $h$ sufficiently small,
\begin{align*}
\sum_{x_h\in A_k} f_{x_h} \le \sum_{x_h\in E^{\epsilon_k}} f_{x_h}
&\le C \frac{1-\epsilon_k}{\epsilon_k} \|f^{1/d} -g^{1/d}\|_{\ell^d(\calC_h^-(w_h^{\epsilon_k}))}\\
&\le C h^2 \frac{1-\epsilon_k}{\epsilon_k} = C k^{-1/p}.
\end{align*}
Thus, we find that
\begin{align*}
\sum_{x_h\in \Omega_h} f_{x_h} \left| \big(\Delta_{\be} (u_h-v_h)(x_h)\big)^+\right|^p
&\le C h^{2p}\left(1+ \sum_{k=1}^{C h^{-2p}} \frac{1}{k^{1/p}}\right)\\
&\le
C \left\{
\begin{array}{ll}
h^2 |\log h| & \text{if $p=1$},\\
h^2 & \text{if $p>1$}.
\end{array}
\right.
\end{align*}

In certain settings,
Theorems \ref{thm:OPMainEst} and \ref{thm:OPW2pResult} immediately
give us $W^{1,p}$ error estimates as well.  To make this precise, we assume
that the basis $\{\tilde \be_j\}_{j=1}^d = \{\be_j\}_{j=1}^d$ defined in \eqref{translationinvariant}
is the canonical one.  We then define the backward difference operator
\begin{equation*}
D^-_{\be} v(x_h) = \frac{v(x_h)-v(x_h-\be h)}{h},
\end{equation*}
and the discrete norms/semi-norms, for $p\in (1,\infty)$,
\begin{align*}
\|v\|_{L^p_h(\Omega_h)} &= \left(h^d \sum_{x_h\in \Omega_h} |v(x_h)|^p\right)^{1/p},\\
\|v\|_{W^{1,p}_h(\Omega_h)} &=\left(\|v\|_{L^p_h(\Omega_h)}^p+h^d\sum_{j=1}^d \|D^-_{\be_j} v\|_{L^p_h(\Omega_h)}^p\right)^{1/p},\\
\|v\|_{W^{2,p}_h(\Omega_h)} & = \left(\|v\|_{W^{1,p}_h(\Omega_h)}^p + 
h^d \sum_{j=1}^d \big(\|\Delta_{\be_j} v\|_{L^p_h(\Omega_h)}^p + \mathop{\sum_{i=1}^d}_{j\neq i} \|D^-_{\be_i} D^-_{\be_j} v\|_{L^p_h(\Omega_h)}^p\big) \right)^{1/p}.
\end{align*}
We then have \cite[Lemmas 2.60--2.61]{SuliBook}
\begin{align*}
\|v\|_{W^{1,p}_h(\Omega_h)} \le C\|v\|_{L^p_h(\Omega_h)}^{1/2} \|v\|_{W^{2,p}_h(\Omega_h)}^{1/2}.
\end{align*}
Therefore noting that $\|v\|_{L^p_h(\Omega_h)}\le C\|v\|_{L^\infty(\Omega)}$, 
{and,
\begin{align*}
D_{\be_i}^- D_{\be_j}^-v(x)
& = \frac12 \left(\Delta_{\be_i} v(x-h\be_i)+\Delta_{\be_j} v(x-h\be_j) -\Delta_{\tilde \be_{i,j}} v(x-h(\be_i+\be_j))\right)
\end{align*}
with $\tilde \be_{i,j} = \be_i-\be_j$,} 
we have the following, by Theorems \ref{thm:OPMainEst} and \ref{thm:OPW2pResult}.
\begin{corollary}
Suppose that 
the conditions in \ref{thm:OPMainEst} are satisfied with $k+\alpha=2$,
and assume that $f\ge f_0>0$ in $\Omega$.  Then there holds
\begin{align*}
\|u-u_h\|_{W^{1,p}_h(\Omega_h)}
& \le \left\{
\begin{array}{ll}
C h^{1+\frac1{2p}} & p\in (d,\infty),\\
C|\log h|^{\frac{1}{2d}} h^{1+\frac{1}{2d}} & p\in (1,d].
\end{array}
\right.
\end{align*}
\end{corollary}

\begin{remark}[extensions]
In this section we 
showed that the stability estimate 
given in Proposition \ref{stability}
provides a powerful tool
to develop error estimates for the \MA
equation, 
as it allows us to derive $L^\infty$
and $W^2_p$ error estimates
when the solution enjoys regularity $u\in C^{2+k,\alpha}(\bar\Omega)$.
Thanks to this stability estimate, it also possible 
to extend these estimates if 
the solution is of lower regularity and/or degenerate. 
The key observation is that the stability estimate measures the consistency error in the $\ell^d$-norm. 
If the solution is rough in a region of small measure and smooth elsewhere, so that the consistency error is small in $\ell^d$-norm, then by the stability estimate, we may still derive a rate of convergence for the low regularity case.
This is explored in  \cite[Theorem 6.3]{NochettoZhangMA} to prove a rate of convergence for solutions in $C^{1,1}(\Omega)$, but not in $C^2(\Omega)$.
\eremk\end{remark}

\section{Finite Element Methods}\label{sec:FEM}

\epigraph{\tiny{\it It will be found that most classical mathematical approximation procedures as well as the various direct approximations used in engineering fall into this category. It is thus difficult to determine the origins of the finite element method and the precise moment of its invention.}}{\tiny{O. Zienkewicz \cite{ZTZ}}}

In this section, we summarize recent developments of finite element methods for the \MA problem with Dirichlet boundary conditions \eqref{eqn:MA}. For simplicity, throughout this section, we assume that boundary conditions in \eqref{eqn:MA} are homogeneous, i.e., $g=0$. The extension to non-homogenous boundary conditions is straightforward.

The main difficulty to construct (and analyze) finite element schemes for fully nonlinear problems is that the PDEs are non-variational. Recall that a finite element method is typically derived by
\begin{enumerate}
  \item[(i)] multiplying the PDE by a test function; 
  \item[(ii)] integrating the resulting product over the domain; 
  \item[(iii)] performing integration by parts to arrive at a variational formulation; 
  \item[(iv)] posing the variational formulation on a finite dimensional space, usually consisting of piecewise polynomials.  
\end{enumerate}
Note that the third step usually requires some structure conditions of the PDE, e.g., that the PDE is in divergence-form, which is not present for fully nonlinear problems. Another obvious difficult to construct convergent finite element schemes is that the notion of viscosity solutions, {given in Definition~\ref{def:viscosol}, and Alexandrov solutions, as in Definition~\ref{def:AlexSolution}} for the \MA equation are non-variational, and it is unclear how this solution concept can be adopted 
within a finite element framework.

{We must remark, however, that the \MA operator \eqref{eqn:MA1} does possess a divergence-form.  Using well-known algebraic identities and the divergence-free property of cofactor matrices, there holds  $\det D^2 u  = \frac1{d}\nabla \cdot \big(\cof  D^2 u \nab u)$. Note however that variational formulations based on this identity would still involve second-order derivatives,
and therefore, at this time, it is unclear whether numerical methods based on this approach are advantageous.}

Nonetheless, assuming some regularity of the solution, 
well-defined finite element methods can be formulated and analyzed for fully nonlinear PDEs.
One approach is to omit the third step of the four--step process described above.
For example, multiplying the \MA equation \eqref{eqn:MA1} 
by a function $v$ and integrating over $\Omega$ yields the identity
\begin{equation}\label{eqn:NaiveC1}
\int_\Omega \big(f -\det D^2 u \big)v\diff x = 0.
\end{equation}
A simple calculation involving H\"older's inequality and Sobolev embeddings show 
that 
expression \eqref{eqn:NaiveC1} is well-defined 
provided $u,v\in W^{2,d}(\Omega)$.
Finite element methods can then be constructed 
based on the identity \eqref{eqn:NaiveC1}.  Namely, an
obvious finite element method based on the identity
\eqref{eqn:NaiveC1} seeks $u_h\in V_h$ satisfying
\begin{equation}\label{eqn:NaiveC1FEM}
\int_\Omega \big(f -\det D^2 u_h \big)v_h\diff x = 0\qquad \forall v_h\in X_h,
\end{equation}
where $X_h$ is a finite dimensional space
consisting of piecewise polynomials with respect
to a partition of $\Omega$ that vanish on the boundary.
While this method may be convergent (cf.~ \cite{Bohmer08,Neilan14,Awanou15F,Awanou14H,DavydovSaeed13,Awanou17ABC}), the appearance
of global second-order derivatives in the method necessitates the use of $C^1$
finite element spaces which
can be 
arduous to implement
and  are not found in most finite element software packages.
In addition, $C^1$ finite element generally require high-degree polynomial
bases, resulting in a relatively large algebraic system.

Because of the many disadvantages of the finite element 
method \eqref{eqn:NaiveC1FEM} several
finite element methods with simpler spaces
have been developed.  These include
$C^0$ penalty methods, discontinuous Galerkin (DG) methods,
mixed finite element methods, and methods based on high-order regularizations.
We now discuss these methods in the subsequent sections.

\subsection{Continuous finite element methods}\label{sec:C0FEMs}

Here we summarize
finite element methods 
presented in \cite{BGNS11,BrennerNeilan12,Neilan13} for the \MA equation
which employ spaces consisting
of continuous, piecewise polynomials, i.e., the Lagrange finite element space.
These are arguably the simplest finite element spaces and are
available on virtually all finite element software programs and libraries.
In addition, we provide a slightly new and improved convergence analysis
based on recent results for finite element methods for linear non-divergence form PDEs \cite{FengNeilan16}.
To describe these methods and their accompanying analysis, we require some notation.

As before, we assume that $\Omega\subset \mathbb{R}^d\ (d=2,3)$ 
is a bounded, convex domain.
Let $\calT_h$ denote a shape-regular and simplicial 
triangulation of $\Omega$.  
We denote the sets of interior and boundary
$(d-1)$--dimensional faces of $\calT_h$ by $\calF_h^I$ and $\calF_h^B$, restrictively.
The jump of a vector valued function $\bv$
across an interior face ${F}= \p T_+\cap \p T_-\in \calF_h^I$
is given by
\begin{align}\label{eqn:JumpDef}
\jump{\bv} = \frac12\Big(\bv_+\otimes \bn_+ +\bn_+\otimes \bv_+ + \bv_-\otimes \bn_- + \bn_-\otimes \bv_-\Big),
\end{align}
where $\bn_{\pm}$ is the outward unit normal of $\p T_{\pm}$,
and $\bv_{\pm} = \bv|_{T_\pm}$.  We also define
the average of $B$ (a scalar, vector, or matrix-valued function) across ${F}$ as
\begin{align}\label{eqn:AvgDef}
\avg{B} = \frac12 \big(B_++B_-\big).
\end{align}
If ${F} = \p T_+\cap \p \Omega\in \calF_h^B$, then we define
\begin{align}\label{eqn:BJumpAvgDef}
\jump{\bv} &= \frac12\big(\bv_+\otimes \bn_++\bn_+\otimes \bv_+\big),\qquad
\avg{B}  = B_+.
\end{align}

For an integer $r\ge 2$, 
the Lagrange finite element
space with homogeneous boundary conditions is given by
\begin{align*}
V_h = \left\{v_h\in W^{1,\infty}_0(\Omega):\ v_h|_T\in \bbP_r(T)\ \forall T\in \calT_h \right\},
\end{align*}
where $\bbP_r(T)$ is the space of polynomials with degree less than or equal to $r$
with domain $T$.
In addition, for a number $p\in (1,\infty)$ and integer $m$, we define 
\begin{align*}
W^{m,p}(\mct) = \prod_{T\in \calT_h} W^{m,p}(T),\qquad
V_p = W^{1,p}_0(\Omega)\cap W^{2,p}(\mct),
\end{align*}
and note that $V_h \subset V_p$ for all $p\in (1,\infty)$.  We also set
$H^m(\mct) = W^{m,2}(\mct)$. 

Because of the non-inclusion $V_h\not\subset W^{2,	d}(\Omega)$,
the finite element formulation \eqref{eqn:NaiveC1FEM}
is not well-defined if $X_h$ is taken to be
the Lagrange finite element space.
A na\"ive
approach to bypass this issue is to redefine
this formulation so that integration is done piecewise over the mesh, i.e.,
to consider
\begin{align}
\label{eqn:NaiveC0FEM}
\sum_{T\in \calT_h} \int_T \big(f -\det D^2 u_h \big)v_h\diff x\qquad \forall v_h\in V_h.
\end{align}
While this method is well-defined (i.e., all quantities are defined and bounded),
it is easy to see that the scheme is ill-posed.  For example,
if $w_h\in V_h$ is strictly piecewise linear, then $\det D^2 w_h =0$
on each $T\in \mct$, and consequently, uniqueness (and stability) is dramatically lost.

The arguments given in \cite{BGNS11} offer an alternative explanation
on why the formulation \eqref{eqn:NaiveC0FEM} leads to an ill-posed
problem.  Namely, the main point  in \cite{BGNS11} is that the linearization 
of the discrete problem \eqref{eqn:NaiveC0FEM} is not consistent
with respect to the linearization of the continuous
problem \eqref{eqn:MA1}.
  Instead, to ensure consistency and stability, finite element methods
for the \MA problem
should be designed such that the discrete linearization at the solution $u$
is a coercive operator over the finite element space.  We now
explain how to construct methods with stable linearizations.
To do so, 
we first 
assume that the exact solution to the \MA equation 
satisfies $u\in C^{k,\alpha}(\bar\Omega)$ with $k+\alpha > 2$
and is strictly convex.  

Define 
\[
{\polF}[u] = f-\det D^2 u
\]
 to be the \MA operator,
and let $L$ be the linearization of $F$ at the solution $u$, i.e.,
\begin{align}\label{eqn:LinearizationOperator}
Lw = \lim_{t\to 0} \frac{{\polF}[u+tw]-{\polF}[u]}{t} = -\cof D^2 u:D^2 w,
\end{align}
where $\cof D^2 u$ denotes the cofactor matrix of $D^2 u$,
and `$:$' denotes the Frobenius inner product.
The assumptions on $u$
imply that matrix $\cof D^2 u $ is positive definite on $\bar\Omega$ and uniformly continuous.

\begin{figure}
\[
\begin{array}{ccccc}
{{\polF}}(u)=0 &\quad\quad\stackrel{\rm discretize}{\longrightarrow}\quad\quad  & {{\polF}_h}(u_h)=0\\
& & \\
\qquad \boldsymbol{\Big\downarrow} \text{\scriptsize linearize} & 
&\qquad \Big\downarrow \text{\scriptsize linearize}\\
& &\\
{L}(w)=0& \qquad \stackrel{\rm discretize}{\longrightarrow} \qquad& {L_{h}}(w_h)=0\\
\end{array}
\]
\caption{\label{fig:commute}A commuting diagram connecting
the nonlinear problems and their discretizations.}
\end{figure}

A consistent
discretization of  linear operators in non-divergence form (such as $L$)
 was introduced in \cite{FengNeilan16}.  In the case 
 that the linear problem is given by \eqref{eqn:LinearizationOperator}, the 
discretization is given by
$L_h:V_p\to V^\prime_h$ with
\begin{equation}
\label{eqn:LhDef}
\begin{aligned}
\bl L_h v,w_h\br &= -\sum_{T\in \mct} \int_T \big(\cof D^2 u :D^2 v\big)w_h\diff x \\ &+\sum_{{F}\in \calF_h^I} \int_F \avg{\cof D^2 u}:\jump{\nab v} w_h\diff s,
\end{aligned}
\end{equation}
where $\bl \cdot,\cdot\br$ denotes the dual pairing
between some Banach space and its dual.
The operator $L_h$ is clearly consistent with $L$: If $v\in W^{2,p}(\Omega)\cap W^{1,p}_0(\Omega)$,
then  $\bl L_h v,w_h\br = \bl L v,w_h\br$ for all $w_h\in V_h$.
In addition, the discrete operator is stable as the next lemma shows.  We refer the reader to \cite{FengNeilan16} for a proof.

\begin{lemma}[{stability}]\label{lem:LhStabilityEstimate}
Define the 
 discrete $W^{2,p}$-norm 
\begin{alignat*}{2}
\|v\|_{W^{2,p}_h(\Omega)}^p 
&: = \|D^2_h v\|_{L^p(\Omega)}^p+\sum_{{F}\in \calF_h^I} h_F^{1-p}\big\|\jump{\nab v}\big\|_{L^p({F})}^p\quad &&1<p<\infty,\\
\|v\|_{W^{2,\infty}_h(\Omega)} &:=
\|D^2_h v\|_{L^\infty(\Omega)} + \max_{{F}\in \calF_h^I} h_F^{-1} \big\|\jump{\nab v}\|_{L^\infty({F})},
\end{alignat*}
{where $D^2_hv$ is the piecewise Hessian of $v$}.
Assume that $u\in C^2(\bar\Omega)$
and is strictly convex over $\bar\Omega$.
Then there exists $h_0>0$
depending on the modulus of continuity
of $D^2 u$, such that for $h\in (0, h_0]$, there holds
the following inf-sup condition $(2\le p<\infty$)
\begin{align*}
\|w_h\|_{W^{2,p}_h(\Omega)}
\leq C   \|{L}_h w_h\|_{L^p_h(\Omega)}:=\sup_{v_h\in V_{h}\backslash \{0\}} \frac{\bl {L}_h w_h,v_h\br}{\|v_h\|_{L^{p'}(\Omega)}}\qquad \forall w_h\in V_h,
\end{align*}
where $1/p+1/p'=1$
\end{lemma}

Based on the definition of $L_h$
and the stability results stated in Lemma \ref{lem:LhStabilityEstimate}
we can develop a consistent discretization
for the \MA problem as well as a convergence theory.
Essentially, its construction is based on the observations
that the expressions
$\int_T \big(\cof D^2 u:D^2 v\big)w_h\diff x$
and $\int_{F} \avg{\cof D^2 u }:\jump{\nab v}w_h\,d s$
are the linearizations of $\int_T \big(f-\det D^2 v \big)w_h$
and $\int_{F} \avg{\cof D^2 v}:\jump{\nab v} w_h\diff s$, respectively,
about the solution $u$.
With this in mind, we define the discrete operator
${\polF}_h:V\to V_h^\prime$ via
\begin{align*}
\bl {\polF}_h[v],w\br = \sum_{T\in \calT_h} \int_T \big(f-\det D^2 v \big)w_h\diff x + \sum_{{F}\in \calF_h^I} \int_{{F}} \avg{\cof D^2 v }:\jump{\nab v}w_h\diff s,
\end{align*}
and consider the finite element method: Find $u_h\in V_h$ such that
\begin{align}\label{eqn:C0PenaltyMethod}
\bl {\polF}_h[u_h],v_h\br = 0\qquad \forall v_h\in V_h.
\end{align}
We immediately see that method \eqref{eqn:C0PenaltyMethod}
is consistent: There holds $\jump{\nab u}|_{F}=0$
over all interior faces ${F}$, and therefore $\bl {\polF}_h[u],v_h\br=0$
for all $v_h\in V_h$. Furthermore, the proceeding
discussion implies that ${L}_h$ is the linearization of ${\polF}_h$:
\begin{align*}
L_hw  = \lim_{t\to 0} \frac{{\polF}_h[u+tw]-{\polF}_h[u]}{t}\quad \text{in $V_h^\prime$}.
\end{align*}
In summary the diagram given in Figure \ref{fig:commute} commutes.
We now show that this property (along with the regularity
and convexity assumptions of $u$)
implies that there exists a locally unique solution to  \eqref{eqn:C0PenaltyMethod}
with optimal rates of convergence.

As a first step, we first point out that
Lemma \ref{lem:LhStabilityEstimate}
implies that $L_h\big|_{V_h}$ is bijective.
Therefore, the mapping
$M_h:V_p\to V_h$ given by
\begin{align}\label{eqn:MhDef}
M_h = \big({L}_h\big|_{V_h}\big)^{-1}\big(L_h - {\polF}_h\big)
\end{align}
is well defined. 
The existence of a solution to the
finite element method \eqref{eqn:C0PenaltyMethod} is proven
by showing that ${M}_h$
has a fixed point in a ball centered at $u_{c,h}$,
where $u_{c,h}$ is the elliptic projection of $u$
given by
\begin{align}\label{eqn:uch}
u_{c,h}:=\big({L}_h\big|_{V_h}\big)^{-1} L_h u.
\end{align}
The basis of this argument
is provided in the next lemma.

\begin{lemma}[{$M_h$ is Lipschitz}]\label{lem:ContractionLemma1}
Assume that the convex solution
of the \MA equation
satisfies $u\in C^{k,\alpha}(\bar\Omega)$ with $k+\alpha>2$. 
Then there holds, for all $p\in [2,\infty)$ and all $v_1,v_2\in V_p$,
\begin{align*}
\|{M}_h v_1-{M}_h v_2\|_{W^{2,p}_h(\Omega)}\le C_1\big\|u-\frac12 (v_1+v_2)\|_{W^{2,\infty}_h(\Omega)} \|v_1-v_2\|_{W^{2,p}_h(\Omega)},
\end{align*}
where $C_1>0$ depends on $p$ and $u$, but
is independent of $h$.
\end{lemma}
\begin{proof}
We give the proof of the two-dimensional case $d=2$; 
the arguments in three dimensions are similar and can be found in \cite{BrennerNeilan12}.

We first use Taylor's Theorem and the fact
that $F_h$ is quadratic in two dimensions, to get
\begin{align*}
  {\polF}_h[v] = {\polF}_h[u] + L_h (v-u) + R_h[v-u] = L_h (v-u) + R_h[v-u],
\end{align*}
where $R_h:V\to V_h^\prime$ is quadratic
in its arguments and independent of $u$.
 
Using this expansion into the mapping $M_h$ yields
\begin{align}\label{eqn:MhExpansion}
M_h[v_1] - M_h[v_2] 
&= \big({L}_h\big|_{V_h}\big)^{-1}\big(L_h v_1-L_h v_2  - \big({\polF}_h[v_1] - {\polF}_h[v_2]\big)\\
&\nonumber = \big(L_h\big|_{V_h}\big)^{-1}\big(R_h[v_2-u] - R_h[v_1-u]\big).
\end{align}
Since $R_h$ is quadratic there holds
\begin{align*}
R_h[v_2-u] - R_h[v_1-u] 
&= \int_0^1 DR_h[t(v_2-u)+(1-t)(v_1-u)](v_2-v_1) \diff t\\
&= DR_h(\frac12 (v_2+v_1)-u)(v_2-v_1),
\end{align*}
where by $DR_h$ we denoted the derivative of $R_h$. 
Therefore, by \eqref{eqn:MhExpansion}
and Lemma \ref{lem:LhStabilityEstimate}
we have
\begin{align*}
\|M_h v_1-M_h v_2\|_{W^{2,p}_h(\Omega)}\leq C  \left\| DR_h\left(\frac12 (v_2+v_1)-u\right)(v_2-v_1)\right\|_{L^p_h(\Omega)}.
\end{align*}
Several applications of H\"older's inequality
yields (cf.~\cite[Lemma 4.2]{Neilan13})
\begin{align*}
\|DR_h(w)(q)\|_{L^p_h(\Omega)} \leq C \|w\|_{W^{2,\infty}_h(\Omega)} \|q\|_{W^{2,p}_h(\Omega)},
\end{align*}
and therefore
\begin{align*}
\|M_h v_1-M_h v_2\|_{W^{2,p}_h(\Omega)}
&\leq C  \big\|\frac12 (v_1+v_2)-u\|_{W^{2,\infty}_h(\Omega)} \|v_1-v_2\|_{W^{2,p}_h(\Omega)}.
\end{align*} 
 \end{proof}
 
 \begin{lemma}[{contraction}]
 \label{lem:ContractionLemma2}
 Assume that the hypotheses of Lemma \ref{lem:ContractionLemma1}
 are satisfied. For fixed $\rho>0$ and $p\in [2,\infty)$, define 
 the closed ball
\begin{align*}
B_{\rho,p} = \left\{v_h\in V_h:\ \|u_{c,h}-v_h\|_{W^{2,p}_h(\Omega)}\le \rho\right\},
\end{align*} 
where $u_{c,h}\in V_h$ is defined by \eqref{eqn:uch}. Then, for all $v_1,v_2\in B_{\rho,p}$,
there holds
\begin{align*}
\|M_h v_1 - M_h v_2\|_{W^{2,p}_h(\Omega)}\le C_2  h^{-d/p}\big(h^{\ell+\alpha} + \rho\big)\|v_1-v_2\|_{W^{2,p}_h(\Omega)},
\end{align*}
where $\ell = \min\{r-2,k-2\}$.
 \end{lemma}
 \begin{proof}
 First, {the smoothness assumptions on $u$ allows us to conclude that} the elliptic projection $u_{c,h}$ satisfies \cite[Theorem 3.2]{FengNeilan16}
\begin{align}\label{eqn:C3Def}
\|u-u_{c,h}\|_{W^{2,p}_h(\Omega)}\leq  C_3 h^{\ell+\alpha} \qquad p\in [2,\infty),
\end{align}
where $C_3>0$ depends on $p$ and $\|u\|_{C^{k,\alpha}(\bar\Omega)}$. 
Consequently, there holds by an inverse estimate, for any $w_h\in V_h$,
\begin{align*}
\|u-u_{c,h}&\|_{W^{2,\infty}_h(\Omega)}
\le \|u-w_h\|_{W^{2,\infty}_h(\Omega)}+ Ch^{-d/p} \|u_{c,h}-w_h\|_{W^{2,p}_h(\Omega)}\\
&\le \|u-w_h\|_{W^{2,\infty}_h(\Omega)}+ C h^{-d/p} \big(\|u-u_{c,h}\|_{W^{2,p}_h(\Omega)}+\|u-w_h\|_{W^{2,p}_h(\Omega)}\big).
\end{align*}
Taking {$w_h$} to be the nodal interpolant of $u$ yields
\begin{align}\label{eqn:C4Def}
\|u-u_{c,h}\|_{W^{2,\infty}_h(\Omega)}\le C_4 h^{\ell+\alpha-d/p}.
\end{align}

Applying this result to Lemma \ref{lem:ContractionLemma1}
and using an inverse estimate, we obtain
\begin{align*}
&\|M_h v_1 - M_h v_2\|_{W^{2,p}_h(\Omega)}\\
&\ \ \le C \left(\|u-u_{c,h}\|_{W^{2,\infty}_h(\Omega)} + h^{-d/p}\|u_{c,h}- \frac12 (v_1+v_2)\|_{W^{2,p}_h(\Omega)}\right)\|v_1-v_2\|_{W^{2,p}_h(\Omega)}\\
&\ \ \le C h^{-d/p} \big(h^{\ell+\alpha} + \rho\big)\|v_1-v_2\|_{W^{2,p}_h(\Omega)}
\end{align*} 
for all $v_1,v_2\in B_{\rho,p}$. 
 \end{proof}

 \begin{theorem}[{error estimate}]
\label{thm:AnErrorEstimate}
 Assume that $u\in C^{k,\alpha}(\bar\Omega)$ with $k+\alpha>2$
 and is strictly convex.
Set $\ell = \min\{r-2,k-2\}$.
 There exists $h_1>0$
 such that for $h\le h_1$,
 there exists a solution to \eqref{eqn:C0PenaltyMethod}
 satisfying
 \begin{align}\label{eqn:AnErrorEstimate}
 \|u-u_h\|_{W^{2,p}_h(\Omega)}\leq C  h^{\ell+\alpha}.
 \end{align}
 Moreover, if $\tilde{u}_h$ is another
 solution to \eqref{eqn:C0PenaltyMethod}
 then there holds $\|u-\tilde{u}_h\|_{W^{2,\infty}_h(\Omega)}\ge C$, with the constant
 $C>0$ independent of $h$.
 \end{theorem}
 \begin{proof}
 Fix $p\in [2,\infty)$ such that $\ell+\alpha - d/p>0$, and let 
\[
h_1 = \min\{1/(4C_2),1/(2C_1C_2C_3C_4)\}^{1/(\alpha+\ell-d/p)}.
\]
Then, {for $h\le \min\{h_0,h_1\}$, where $h_0$ was defined in Lemma~\ref{lem:LhStabilityEstimate}}, set 
$\rho_1 = h^{\ell+\alpha}/(4C_2)$.
Lemma \ref{lem:ContractionLemma2}
then shows that, for $v_1,v_2\in B_{\rho_1,p}$,
\begin{align*}
\|{M}_h v_1 - {M}_h v_2 \|_{W^{2,p}_h(\Omega)}&\le C_2 \big(h^{\ell+\alpha-d/p}+ h^{-d/p}\rho_1\big)\|v_1-v_2\|_{W^{2,p}_h(\Omega)}\\
& \le 2C_2 h_1^{\alpha+\ell-d/p} \|v_1-v_2\|_{W^{2,p}_h(\Omega)} 
\le \frac12 \|v_1-v_2\|_{W^{2,p}_h(\Omega)},
\end{align*}
and therefore ${M}_h\big|_{V_h}$ is a contraction mapping on $B_{\rho_1,p}$.
Likewise, we can use Lemma  \ref{lem:ContractionLemma1} and the fact
that $u_{c,h}= M_h u$ to get (cf.~\eqref{eqn:C3Def}--\eqref{eqn:C4Def})
\begin{align*}
\|u_{c,h}-M_h v\|_{W^{2,p}_h(\Omega)} 
&= \|M_h u - M_h v\|_{W^{2,p}_h(\Omega)}\\
&\le \frac{C_1}{2} \|u-u_{c,h}\|_{W^{2,\infty}_h(\Omega)}\|u-u_{c,h}\|_{W^{2,p}_h(\Omega)}\\
&\le \frac{C_1 C_3 C_4 h^{2\alpha+2\ell-d/p}}{2}\le \frac{h^{\ell+\alpha}}{4C_2} = \rho_1.
\end{align*}
Therefore $M_h$ maps $B_{\rho_1,p}$ to itself. By Banach's fixed point theorem,
we conclude that $M_h$ has a fixed point in $B_{\rho_1,p}$, and this fixed point
is a solution to \eqref{eqn:C0PenaltyMethod}.  The error estimate for $\alpha-d/p>0$
\eqref{eqn:AnErrorEstimate} follows from the inclusion $u_h\in B_{\rho_1,p}$ and
the definition of $\rho_1$.  The other cases $\ell+\alpha-d/p\le 0$ then follow from H\"older's inequality. 

Finally, if $\tilde{u}_h\in V_h$ is another solution to \eqref{eqn:C0PenaltyMethod},
then there holds $M_h \tilde{u}_h = \tilde{u}_h$.  Therefore, by Lemma
\ref{lem:ContractionLemma1} we conclude that
\begin{align*}
\|\tilde{u}_h-u_h\|_{W^{2,p}_h(\Omega)} 
&= \|M_h \tilde{u}_h - M_h u_h\|_{W^{2,p}_h(\Omega)}\\
&\le \frac{C_1}2 \big(\|u-u_h\|_{W^{2,\infty}_h(\Omega)}+\|u-\tilde{u}_h\|_{W^{2,\infty}_h(\Omega)}\big)\|u_h-\tilde{u}_h\|_{W^{2,p}_h(\Omega)}.
\end{align*}
Now applying similar arguments as those found in Lemma 
\ref{lem:ContractionLemma1}, we conclude that $\|u-u_h\|_{W^{2,\infty}_h(\Omega)}\le C h^{\alpha-d/p}\to 0$.
Therefore, by dividing by $\|u_h-\tilde{u}_h\|_{W^{2,p}_h(\Omega)}$, 
we get $C\le \|u-\tilde{u}_h\|_{W^{2,\infty}_h(\Omega)}$ for $h$ sufficiently small.
 \end{proof}

\begin{remark}[extensions]
{The proposed method and the conclusion of Theorem~\ref{thm:AnErrorEstimate} deserve the following comments:}
\begin{enumerate}
\item[$\bullet$] As mentioned earlier, 
 the analysis given here slightly improves
 the results given in \cite{BGNS11,Neilan13}.
Namely, the paper \cite{BGNS11} 
requires $d=2$, $r\ge 3$, and $u\in H^s(\Omega)$ for $s>3$ (implying
that $u\in C^{2,\alpha}(\bar\Omega)$ by a Sobolev embedding).
The paper \cite{Neilan13} requires $r\ge 2$ and regularity $u\in W^{3,\infty}(\Omega)$
to carry out the analysis.
\item[$\bullet$] Discontinuous Galerkin methods have also been developed under this methodology in \cite{Neilan13}.
The analysis carried out in this section can be applied to these methods
using the recent results for non-divergence PDEs given in \cite{TooManySelfCitations18}.
\item[$\bullet$]A two--grid method to solve the nonlinear method has recently been proposed in \cite{AwanouLiMalitz}.
\end{enumerate}
\eremk\end{remark}

\subsection{Mixed formulations}\label{sec:MixedFEMs}
In this section we {describe} 
mixed finite element formulations for the \MA equation
proposed in \cite{LakkisPryer11ABC,Neilan14ABC,AwanouLi14,Awanou15E,AwanouOops,KaweckiLakkisPryer18}.
Essentially, the main idea in these approaches is to introduce the Hessian matrix 
of $u$ as an additional auxiliary unknown in the formulation 
of the \MA problem, that is, we write the PDE  \eqref{eqn:MA1}
as
\begin{align}\label{eqn:MixedMotivation1}
\sigma = D^2 u,\qquad \det \sigma  = f\qquad \text{in }\Omega.
\end{align}
As before, assuming regularity $u\in W^{2,d}(\Omega)$ so that $\sigma \in L^d(\Omega)$, we can multiply the second equation by a {smooth} test function
and integrate over the domain:
\begin{align}\label{eqn:MixedMotivation2}
\int_\Omega (f - \det \sigma )v\diff x=0 
\end{align}
for all $v\in L^\infty(\Omega)$.  

The direct analogue of this formulation in the discrete setting
requires $C^1$ finite element spaces by {the same reasons that the method described in Section~\ref{sec:C0FEMs} does. In other words}, to ensure that the discrete
version of \eqref{eqn:MixedMotivation2} is well--defined,
we require that the Hessian of the discrete approximation
$u_h$ has (global) second-order derivatives in $L^d(\Omega)$;
if $u_h$ is a piecewise polynomial, then this restriction
implies that $u\in C^1(\Omega)$.
To relax this restriction on the finite element spaces, one
can instead develop finite element methods that only
employ continuous (or discontinuous) bases
based on this formulation by introducing the notion
of a discrete Hessian (also known as a  finite element Hessian \cite{LakkisPryer11ABC}).
The discrete Hessian is defined globally via an {integration by parts} procedure rather
than a piecewise fashion.  This idea has been carried out for (linear) Kirchhoff plates 
 in \cite{HuangHuanHan10}, and its formulation is reminiscent of the construction of local discontinuous Galerkin 
methods for second order problems \cite{CockburnShu98,UnifiedDG02}.

To motivate the definition of the discrete Hessian, 
we introduce the auxiliary space 
 \[
\Sigma_h = \{\tau_h\in L^\infty(\Omega;\mathbb{R}^{d\times d}):\ \tau_h|_T\in \mathbb{P}_r(T;\mathbb{R}^{d\times d})\ \forall T\in \calT_h\},
\]
and note the following integration by parts
identity
\begin{equation}\label{eqn:IBPMixed}
\begin{aligned}
\sum_{T\in \mct}  \int_T D^2 w:\tau_h\diff x 
& = - \sum_{T\in \mct} \int_T (\nab \cdot \tau_h)\cdot \nab w\diff x \\ &+ \sum_{T\in \mct}\int_{\p T}  (\tau_h \bn_T)\cdot \nab w\diff s,
\end{aligned}
\end{equation}
for all $w\in H^2(\Omega)$ and $\tau_h\in \Sigma_h$.
Here, $\bn_T$ is the outward unit normal of $\p T$, and
the divergence acting on a matrix is performed row-wise.
We may then write the integral boundary
terms in \eqref{eqn:IBPMixed} 
using the jump and average operators.
In addition to \eqref{eqn:JumpDef}--\eqref{eqn:BJumpAvgDef}, we define
the jump of a matrix-valued function $\tau$ across $F = \p T_+\cap T_-\in \calF_h^I$ as
\begin{align*}
\jump{\tau} = \tau_+\bn_++\tau_-\bn_-,
\end{align*}
and define $\jump{\tau} = \tau_+ \bn_+$ if $F = \p T_+\cap \p \Omega\in \calF_h^B$.
We then have
\begin{align*}
\sum_{T\in \mct}\int_{\p T}  (\tau_h \bn_T)\cdot \nab w\diff s
& = \sum_{F\in \calF_h^I} \int_F \avg{\tau_h}:\jump{\nab w}\diff s + \sum_{F\in \calF_h} \int_F \jump{\tau_h} \cdot \avg{\nab w}\diff s\\
& =  \sum_{F\in \calF_h} \int_F \jump{\tau_h}\cdot \avg{\nab w}\diff s,
\end{align*}
where we used that $\jump{\nab w}|_F=0$ for all $F\in \calF_h^I$ due
to the regularity $w\in H^2(\Omega)$. Combining this identity
with \eqref{eqn:IBPMixed}, we arrive at
\begin{align*}
\sum_{T\in \mct}  \int_T D^2 w:\tau_h\diff x 
& = - \sum_{T\in \mct} \int_T (\nab \cdot \tau_h)\cdot \nab w\diff x + \sum_{F\in \calF_h} \int_F \jump{\tau_h}\cdot \avg{\nab w}\diff s.
\end{align*}

This identity leads to the following definitions of the discrete Hessian.
\begin{definition}[{discontinuous discrete Hessian}]
\label{def:DDiscreteHessian}
The {\it discontinuous discrete \\ Hessian}
is the operator $\bbH_h:H^1(\Omega)\cap H^2(\mct)\to \Sigma_h$
uniquely defined by the conditions
\begin{align*}
 \int_\Omega \bbH_h(w):\tau_h\diff x 
& = - \sum_{T\in \mct} \int_T (\nab \cdot \tau_h)\cdot \nab w\diff x + \sum_{F\in \calF_h} \int_F \jump{\tau_h}\cdot \avg{\nab w}\diff s
\end{align*}
for all $\tau_h\in \Sigma_h$.

\end{definition}

\begin{remark}[{characterization through liftings}]
Define the lifting operator
$\Theta:{L^2(\calF_h^I;\R^d)}\to \Sigma_h$ via
\begin{align*}
\int_\Omega \Theta(\bv):\tau_h\diff x = -\sum_{F\in \calF_h^I} \int_F \avg{\tau_h}:\jump{\bv}\diff s\qquad \forall \tau_h\in \Sigma_h.
\end{align*}
{Integrating by parts we obtain}
\begin{align*}
\sum_{T\in \mct}  \int_T \bbH_h(w):\tau_h\diff x 
&= \sum_{T\in \mct} \int_T D^2 w:\tau_h\diff x - \sum_{F\in \calF_h^I} \int_F \avg{\tau_h}:\jump{\nab w}\diff s\\
&= \sum_{T\in \mct} \int_T \big(D^2 w+\Theta(\nab w)\big):\tau_h\diff x.
\end{align*}
Recalling that  $D^2_h w$ denotes the piecewise Hessian of $w$,
and  that $V_h$ is the (scalar) Lagrange space 
of degree $r$, {we} then have $D^2_h V_h \subset \Sigma_h$,
and therefore
\[
\bbH_h(w_h) = D^2_h w_h + \Theta (\nab w_h)\qquad \forall w_h\in V_h.
\]
\eremk\end{remark}

The notion of the discrete Hessian
and the formal identities \eqref{eqn:MixedMotivation1}--\eqref{eqn:MixedMotivation2}
lead to the following scheme introduced in \cite{Neilan14ABC} : Find $u_h\in V_h$
such that
\begin{align}\label{eqn:MNMethod}
\int_\Omega \big(f-\det \bbH_h(u_h)\big)v_h\diff x\qquad \forall v_h\in V_h.
\end{align}
\begin{remark}[{mixed formulation}]
While \eqref{eqn:MNMethod} is written in primal form,
the problem is in fact a mixed finite element method.
Introducing $\sigma_h = \bbH_h(u_h)\in \Sigma_h$, we see 
from the definition of the discrete Hessian
that \eqref{eqn:MNMethod} is equivalent to the system
\begin{subequations}
\label{eqn:MNMixedForm}
\begin{align}
\label{eqn:MNMixedForm1}
\int_\Omega \sigma_h:\tau_h\diff x + \int_\Omega (\nab\cdot \tau_h)\cdot u_h\diff x 
-\sum_{F\in \calF_h} \int_F \jump{\tau_h}\cdot \avg{\nab u_h}\diff s & = 0, \\
%
\label{eqn:MNMixedForm2}
\int_\Omega \big(f - \det \sigma_h \big)v_h\diff x&=0,
\end{align}
\end{subequations}
for all $(\tau_h, v_h) \in \Sigma_h \times V_h$.
Note that the matrix representation of the form
$(\sigma_h,\tau_h)\to \int_\Omega \sigma_h:\tau_h\diff x$
is symmetric positive definite, and more importantly, block-diagonal
because $\Sigma_h$ does not have any continuity constraints.
As a result, the Schur complement (i.e., 
the primal method \eqref{eqn:MNMethod}) represents
a sparse algebraic system of equations.
\eremk\end{remark}

 \begin{theorem}[{error estimate}]\label{thm:DHessErrorEstimate1}
Assume that $d=2$, and that \eqref{eqn:MA} has a unique 
strictly convex solution $u\in C^{r+3,\alpha}(\Omega)$
with $r\ge 3$ and $\alpha>0$. 
Then for $h$ sufficiently small, there exists
a locally unique solution to the finite element
method \eqref{eqn:MNMethod}.  Moreover, there holds
\begin{align}\label{eqn:DHessErrorEstimate1}
\|u-u_h\|_{H^1(\Omega)}+ h \|\sigma-\sigma_h\|_{L^2(\Omega)}\le C h^r.
\end{align}
\end{theorem}
\begin{proof}
See \cite[Theorem 4.2]{Neilan14ABC}.
\end{proof}

 \begin{remark}[{regularity}]
 The regularity assumptions on $u$ in Theorem \ref{thm:DHessErrorEstimate1}
 can be relaxed using the stability 
  analysis for linear non-divergence
form PDES found in \cite{NeilanSelfCitation17}.  There it is {shown} that,
assuming $u\in C^{2}(\bar\Omega)$,
\begin{align*}
\|w_h\|_{W^{2,2}_h(\Omega)}\leq C  \|L_h w_h\|_{L^2_h(\Omega)}\qquad \forall w_h\in V_h,
\end{align*}
with 
\[
\bl L_h w_h,v_h\br = -\int_\Omega \cof D^2  u:\mathbb{H}_h(w_h) v_h\diff x.
\]
By applying the same techniques found in the previous section,
it is simple to show that the solution to \eqref{eqn:MNMethod}
satisfies   
$\|u-u_h\|_{W^{2,2}_h(\Omega)}\leq C  h^{\ell+\alpha}$ 
with $\ell = \min\{r-2,k-2\}$ provided
that $u\in C^{k,\alpha}(\bar\Omega)$
with $k+\alpha>3$, $r\ge 3$, and $h$ is sufficiently small. 
\eremk\end{remark}

To reduce the number of unknowns
in the mixed system \eqref{eqn:MNMixedForm},
continuity constraints can be added 
in the matrix--valued space $\Sigma_h$.
This is the idea of the method proposed 
in \cite{LakkisPryer13}.  There, the auxiliary space
is defined as the matrix-valued Lagrange space, i.e.,
\begin{align*}
 \Sigma^c_h:=\Sigma_h\cap H^1(\Omega;\bbR^{d\times d})
 = \{\tau_h\in H^1(\Omega):\ \tau_h\in \bbP_r(T;\bbR^{d\times d})\ \forall T\in \calT_h\}.
\end{align*}
Restricting Definition \ref{def:DDiscreteHessian} to 
$\Sigma_h^c$ leads
to the following notation
of the discrete Hessian.
\begin{definition}[{continuous discrete Hessian}]\label{def:CDiscreteHessian}
The {\it continuous discrete Hessian}
is the operator $\bbH^c_h:H^1(\Omega)\cap H^2(\calT_h)\to \Sigma^c_h$
uniquely defined by the conditions
\begin{align*}
\int_\Omega \bbH^c_h(w):\tau_h\diff x 
& = - \int_\Omega (\nab \cdot \tau_h)\cdot \nab w\diff x  + {\int_{\p \Omega}} (\tau_h \bn)\cdot {\nab w} \diff s
\end{align*}
for all $\tau_h\in \Sigma^c_h$.
\end{definition}

This definition leads to a finite element method 
proposed in \cite{LakkisPryer13} which
similar to \eqref{eqn:MNMethod}, but with 
the continuous version of the discrete Hessian.  

\begin{alignat}{2}
\label{eqn:LPMethod}
\int_\Omega \big(f-\det \bbH_h^c(u_h) \big)v_h\diff x =0\qquad &&\forall v_h\in V_h.
\end{alignat}
As before, we may set $\sigma_h = \bbH_h^c(u_h)$ as an auxiliary
variable, and deduce from Definition \ref{def:CDiscreteHessian}
that \eqref{eqn:LPMethod} is equivalent to the mixed method
\begin{subequations}
\label{eqn:LPMethodMixed}
\begin{align}
\label{eqn:LPMethodMixed1}
\int_\Omega \sigma_h:\tau_h\diff x + \int_\Omega (\nab \cdot \tau_h)\cdot \nab u_h\diff x - \int_{\p \Omega}  (\tau_h \bn)\cdot \nab u_h\diff s& = 0, \\
\label{eqn:LPMethodMixed2}
 \int_\Omega\big(f-\det \sigma_h \big)v_h\diff x&=0,
\end{align}
\end{subequations}
for all $(\tau_h, v_h) \in \Sigma_h^c \times V_h$
Compared with the formulation using
the discontinuous discrete Hessian, the mixed problem
\eqref{eqn:LPMethodMixed} has significantly less unknowns
than \eqref{eqn:MNMixedForm} due to the continuity
restrictions of $\Sigma_h^c$.  On 
the other hand, the (mass) matrix associated
with the form $(\sigma_h,\tau_h)\to \int_\Omega \sigma_h:\tau_h$
is not block-diagonal, and therefore the Schur complement
of \eqref{eqn:LPMethodMixed} (i.e., the algebraic system representing the primal problem \eqref{eqn:LPMethod})
is dense.

Existence, (local) uniqueness, 
and error estimates for method 
\eqref{eqn:LPMethodMixed}
{are} similar to the statements
given in Theorem \ref{thm:DHessErrorEstimate1}.

 \begin{theorem}[{error estimates}]
 \label{thm:DHessErrorEstimate2}
Assume that $d\in \{2,3\}$, and that \eqref{eqn:MA} has a unique 
strictly convex solution $u\in H^{r+3}(\Omega)$
with $r\ge d$.  Then for $h$ sufficiently small, there exists
a locally unique solution to the finite element
method \eqref{eqn:LPMethodMixed}.  Moreover, there holds
\begin{align}\label{eqn:DHessErrorEstimate2}
\|u-u_h\|_{H^1(\Omega)}+ h \|\sigma-\sigma_h\|_{L^2(\Omega)}\le C h^r.
\end{align}
 \end{theorem}
\begin{proof}
See \cite[Theorem 3.13]{AwanouLi14} and \cite[Theorem 1]{Awanou15E,AwanouOops}.
\end{proof}

\begin{remark}[{extension to optimal transport}]
The mixed finite element method 
\eqref{eqn:LPMethodMixed} has recently been
extended to the optimal transport problem in 
\cite{KaweckiLakkisPryer18}.
\eremk\end{remark}

\begin{remark}[historical remark]
Our presentation follows a reverse chronological order. The first Galerkin-type method based on the concept of discrete Hessians was that of \cite{LakkisPryer13}, where they used the continuous Hessian of Definition~\ref{def:CDiscreteHessian}. The DG version was introduced later.
\eremk\end{remark}

\subsection{Galerkin methods for singular solutions}\label{sec:GalerkinSingular}

The analysis of the Galerkin methods discussed 
thus far require relatively stringent regularity conditions
to carry out the analysis (e.g., $u\in C^{2,\alpha}(\Omega)$).
While numerical experiments indicate that regularity
assumptions can be relaxed somewhat, they also indicate
that some regularity of the solution is required for the 
methods to converge.  For example, the numerical experiments
in \cite{BGNS11} indicate that the $C^0$ penalty method
\eqref{eqn:C0PenaltyMethod} does not converge if $u\not\in H^2(\Omega)$ in two dimensions.
In this section, we discuss various ways to modify the Galerkin methods and the analysis
such that the resulting numerical scheme is robust with respect to the solution's
regularity.

The first approach, introduced in \cite{FengNeilan09Z}, regularizes
the problem at the PDE level by adding a higher-order perturbation, resulting
in a fourth-order, quasi-linear problem.
The motivation of this approach is that solutions of the regularized
problem are defined via variational principles, {so that weak formulations can be obtained via integration by parts},
and therefore the resulting PDE framework is amenable to Galerkin methods.
Applying this methodology to the \MA problem results in
\begin{subequations}
\label{eqn:VMM}
\begin{alignat}{2}
\label{eqn:VMM1}
-\epsilon \Delta^2 u^\epsilon + \det D^2 u^\epsilon
& = f\qquad &&\text{in }\Omega,\\
\label{eqn:VMM2}
u& = 0\qquad &&\text{on }\p\Omega,
\end{alignat}
where $\epsilon>0$ and $\Delta^2 = \Delta \Delta$ denotes the biharmonic operator.
Note that, due to the higher order of the PDE, 
the Dirichlet boundary condition is no longer sufficient to close the system.
In \cite{FengNeilan09Z}, the following additional boundary conditions
are proposed:
\begin{align}
\label{eqn:VMM3}
\Delta u^\epsilon = 0,\quad \text{or}\quad \frac{\p \Delta u^\epsilon}{\p \bn} = 0\qquad \text{on }\p\Omega.
\end{align}
\end{subequations}
These conditions are chosen {so that} the resulting boundary layer is minimized; see \cite{FengNeilan09Z} for details.
For {the sake of illustration}, we take the first boundary condition in  \eqref{eqn:VMM3} in the discussion below.

Since the problem \eqref{eqn:VMM} is quasi--linear and in divergence--form, 
the notion of weak solutions is easily defined.
%
\begin{definition}[{weak solution}]
A function $u\in W^{2,d}(\Omega)\cap W^{1,d}_0(\Omega)$ is 
a weak solution to \eqref{eqn:VMM} provided that 
\begin{align}\label{eqn:VMMVariationalForm}
-\epsilon \int_\Omega \Delta u^\epsilon \Delta v\diff x  +\int_\Omega v \det D^2 u^\epsilon   \diff x = \int_\Omega f v\diff x\qquad
\forall v\in W^{2,d}(\Omega)\cap W^{1,d}_0(\Omega).
\end{align}
The function $u=\lim_{\epsilon\downarrow 0} u^\epsilon$, if it exists,
is called a weak (resp., strong) moment solution to the \MA problem
if convergence holds in a $W^{1,d}$-weak (resp., $W^{2,d}$-weak) topology.
\end{definition}
\begin{remark}[{relation to other solution concepts}]
Except in very simple settings (e.g., radially symmetric solutions \cite{FengNeilan14Radial}),
the existence of moment solutions and their relation with viscosity
and Alexandrov solutions is an open problem.  Nonetheless,
numerical experiments indicate that this methodology leads
to robust numerical methods with respect to regularity of the solution
of the \MA equation.  For example, numerical methods
applied to problem \eqref{eqn:VMM} are able to capture
viscosity/Alexandrov solutions that are merely Lipschitz continuous.
\eremk\end{remark}

Constructing methods 
for the regularized problem \eqref{eqn:VMMVariationalForm}
can be done by applying any of the above Galerkin
methods described above;
one only
needs to tack on a consistent and stable discretization of the biharmonic
operator to the discrete formulation.
For example, the simplest method, at least {in theory}, 
is to restrict the variational formulation \eqref{eqn:VMMVariationalForm}
onto a finite dimensional subspace of $W^{2,d}(\Omega)\cap W^{1,d}_0(\Omega)$.
This results in the method to find $u_h^\epsilon\in X_h$
satisfying 
\begin{align}
\label{eqn:VMMC1}
\epsilon \int_\Omega \Delta u^\epsilon_h \Delta v_h\diff x + \int_\Omega \big(f-\det D^2 u_h^\epsilon \big)v_h\,d x=0\quad \forall v_h\in X_h,
\end{align}
with $X_h\subset C^1(\Omega)\cap W^{1,d}_0(\Omega)$.   A convergence analysis
of this discrete problem has been done in \cite{FengNeilanJSC11}.
There it is shown that, if there exists a moment solution with
sufficient regularity, then there exists a locally unique
solution to the discrete problem \eqref{eqn:VMMC1}.

Analogously, combining the $C^0$ finite element
method \eqref{eqn:C0PenaltyMethod}
with the symmetric $C^0$ interior penalty method
for the biharmonic problem introduced in \cite{EngelEtAl02,BrennerSung05}
results in the method: Find $u^\epsilon_h\in V_h$ satisfying
\begin{multline}
\label{eqn:VMMC0IP}
  \epsilon \sum_{T\in \mathcal{T}_h} \int_T \Delta u_h^\epsilon \Delta v_h\diff x \\
  -\epsilon \sum_{F\in \calF^I_h} \int_F \left(\avg{\Delta u_h^\epsilon }(I :\jump{\nab v_h}) 
  +\avg{\Delta v_h } (I :\jump{\nab u^\epsilon_h}) \right. \\
  \left. -\frac{\sigma}{h_F} \jump{\nab u_h^\epsilon}:\jump{\nab v_h}\right)\diff s
 + \sum_{T\in \calT_h} \int_T \big(f-\det D^2 u_h^\epsilon)v_h\diff x\\
+\sum_{F\in \calF_h^I} \int_F \avg{\cof D^2 u_h^\epsilon }:\jump{\nab u_h ^\epsilon} v_h\diff s=0
\end{multline}
for all $v_h\in V_h$.  Here,
$\sigma>0$ is a penalty parameter, and we recall that $I$ denotes the $d\times d$ identity matrix and $V_h$ is the Lagrange finite element space
of degree $r\ge 2$ with homogeneous Dirichlet boundary conditions.
The method \eqref{eqn:VMMC0IP} can be written succinctly as
\[
\epsilon \langle A_h u_h^\epsilon,v_h\rangle +\langle {\polF}_h[u_h^\epsilon],v_h\rangle =0\qquad \forall v_h\in V_h,
\]
where the operator $\polF_h$ is defined by \eqref{eqn:C0PenaltyMethod}, and 
$A_h$ is a consistent discretization of the biharmonic operator given by
\begin{align*}
\bl A_h w,v_h\br &= \sum_{T\in \mathcal{T}_h} \int_T \Delta w \Delta v_h\diff x
- \sum_{F\in \calF^I_h} \int_F \Big(\avg{\Delta w }(I_d :\jump{\nab v_h}) \\
&\qquad+\avg{\Delta v_h } (I_d :\jump{\nab w_h})-\frac{\sigma}{h_F} \jump{\nab w}:\jump{\nab v_h}\Big)\diff s.
\end{align*}
Arguments given in \cite{EngelEtAl02,BrennerSung05} show that there exists
$\sigma_0>0$, independent of $h$, such that 
$\bl A_h v_h,v_h\br \ge C\|v_h\|_{W^{2,2}_h(\Omega)}^2$ for all $v_h\in V_h$
provided that $\sigma\ge \sigma_0$.  Moreover, there holds $\epsilon \bl A_h u^\epsilon,v_h\br+\bl \polF_h[u^\epsilon],v_h\br=0$
for all $v_h\in V_h$ provided that $u^\epsilon\in H^s(\Omega)$ for some $s>5/2$.  Thus, 
the method \eqref{eqn:VMMC0IP} is consistent.

While a convergence analysis of the regularized PDE \eqref{eqn:VMM}
and the discretization \eqref{eqn:VMMC0IP} is an open problem,
we show, via numerical experiments in the next section,
that the method is able to capture non-smooth solutions for the \MA problem
in a variety of settings.  In addition, as shown in \cite{BGNS11}, Newton's method
is robust for
the regularized solution, which allows a natural way to construct initial guesses
for the (unregularized) problem \eqref{eqn:C0PenaltyMethod}.

\subsubsection{Convergence of interior discretizations}
Recent results given in \cite{Awanou17ABC,AwanouAwi16,Awanou15Z,Awanou16T}
argue that, in certain settings, standard discretizations (both finite element and finite difference)
for the \MA equation converge to the Alexandrov solution
as the discretization parameter tends to zero.  Here, in this section, we summarize
these results and the techniques to obtain them.

{As always, we} assume that $\Omega$ is convex. {More importantly, we assume also} that the Dirichlet boundary conditions
can be extended to a function $\tilde{g}$ that is convex on $\Omega$.
Note that the existence of $\tilde{g}$ is guaranteed if the domain 
is strictly convex.  However, due to our assumption that $u|_{\p\Omega}=0$,
we may simply take $\tilde{g}\equiv 0$ in our setting.
We further assume that  $f\in C(\bar\Omega)$ with $f\ge C>0$ on $\Omega$.
Let $\{f_m\}_{m=0}^\infty\subset C^\infty(\bar\Omega)$ be a sequence 
of approximations of $f$ with $f_m\to f$ uniformly on $\bar\Omega$
and $f_m\ge C>0$ for all $m$.  
We then consider the PDE problem
\begin{subequations}
\label{eqn:SourceRegularized}
\begin{alignat}{2}
\det D^2 u_m  & = f_m\qquad &&\text{in }\Omega,\\
u_m & = 0\qquad &&\text{on }\p\Omega.
\end{alignat}
\end{subequations}
Even though the source data of this problem is smooth, in general
there does not exist smooth solutions to \eqref{eqn:SourceRegularized}
because $\Omega$ is not necessarily strictly convex nor smooth, {see Theorem~\ref{thm:existenceclassical}}.
Nonetheless, there exists a unique (convex) Alexandrov solution $u_m\in C(\bar\Omega)$.

Let $\tilde\Omega\subset \Omega$ be a strict subdomain 
of $\Omega$ that is polyhedral and convex, and let $\tilde{\calT}_h$ be a simplicial triangulation 
of $\tilde\Omega$.  Finally, we denote by $\tilde{X}_h$ 
a $C^1(\overline{\tilde\Omega})$-conforming finite element
space consisting of piecewise polynomials with respect to $\tilde{\calT}_h$.
We then consider the finite element method: Find $\tilde{u}_h\in \tilde X_h$
satisfying $\tilde{u}_{m,h}=u_m$ on $\p\tilde\Omega$ and 
\begin{align}\label{eqn:AwanouC1}
\int_{\tilde\Omega} \big(f_m-\det D^2 \tilde{u}_{m,h} \big)v_h\diff x=0\qquad \forall v_h\in \tilde{X}_h\cap W^{1,d}_0(\tilde\Omega).
\end{align}
This method is similar to \eqref{eqn:NaiveC1FEM}, the  differences
being
\begin{enumerate}
  \item[(i)] the problem is posed on $\tilde{\Omega}$ instead of $\Omega$;
  \item[(ii)] the source function has been regularized; 
  \item[(iii)] the homogeneous Dirichlet boundary conditions have been replaced by 
  \[
    \tilde u_{m,h}|_{\p\tilde\Omega}=u_m|_{\p\tilde\Omega}.
  \]
\end{enumerate}
It is clear that method \eqref{eqn:AwanouC1} is a discretization of the PDE problem
\begin{subequations}
\label{eqn:MAsubDomain}
\begin{alignat}{2}
\det D^2 \tilde u_m  & = f_m\qquad &&\text{in }\tilde\Omega,\\
\tilde u_m & = u_m\qquad &&\text{on }\p\tilde\Omega,
\end{alignat}
\end{subequations}
which, similar to \eqref{eqn:SourceRegularized}, has a unique {Alexandrov}
solution and is generally non-smooth.  In fact, it is simple to see that,
due to the inclusion $\tilde\Omega\subset \Omega$ and the uniqueness
of {Alexandrov} solutions, that $\tilde{u}_m=u_m$ on $\tilde{\Omega}$.

\begin{theorem}[{interior convergence}]
\label{thm:InteriorConv}
There exists $h_0>0$, which depends
on ${\rm dist}\{\p\Omega,\p\tilde\Omega\}$,
such that for $h\le h_0$, there exists
a locally unique solution to \eqref{eqn:AwanouC1}.
In addition, as $h\to 0$, $\tilde u_{m,h}$ converges uniformly
to $\tilde u_m$ (the solution to \eqref{eqn:MAsubDomain})
on compact subsets of $\tilde{\Omega}$.
\end{theorem}
\begin{proof}[Proof (sketch)]
The proof relies on a series of smooth approximations to problem \eqref{eqn:MAsubDomain}. Let $\{\Omega_s\}_{s=0}^\infty$ be a sequence
of strictly convex and smooth domains such that $\Omega_s \subset \Omega_{s+1}\subset \Omega$ 
for all $s$, and $\Omega_s\to \Omega$
as $s\to \infty$;  see Figure \ref{fig:AwanouConstruction}.  Consider the problem
\begin{alignat*}{2}
\det D^2 u_{ms} & = f_m\qquad &&\text{in }\Omega_s,\\
u_{ms} & = 0 \qquad &&\text{on }\p\Omega_s.
\end{alignat*}
Note that, because the data is regular, and since $\Omega_s$
is uniformly convex with smooth boundary, 
the solution to this problem is smooth. 
In particular, interior Schauder estimates \cite[Section 6.1]{GT} show that, for any $D\subset\subset \Omega_s$,
\[
\|u_{ms}\|_{C^{r+1}(D)}\le C_m,
\]
where $C_m>0$ depends on $m$, $f_m$, $D$, and ${\rm dist}\{D,\p \Omega_s\}\le {\rm dist}\{D,\p \Omega\}$.
Moreover, results in \cite{Savin13} show that $u_{ms}$ (up to subsequence) converges uniformly
on compact subsets of $\Omega$ as $s\to \infty$.
Now, because $u_{ms}$ is smooth, and because the derivatives of $u_{ms}$ are  uniformly 
bounded on $\tilde{\Omega}$ (with respect to $s$), arguments similar those given 
in the previous section (see \cite{Bohmer08,Awanou15c}) show that, for $h\le h_0$ with $h_0$ sufficiently small, there exists a locally 
unique and convex solution to the following discrete problem: Find $\tilde u_{ms,h}\in \tilde X_h$ satisfying $\tilde u_{ms,h}|_{\p\tilde\Omega}= u_{ms}|_{\p\tilde\Omega}$ and
\begin{align*}
\int_{\tilde\Omega} \big(f_m-\det D^2 \tilde{u}_{ms,h} \big)v_h\diff x=0\qquad \forall v_h\in \tilde{X}_h\cap W^{1,d}_0(\tilde\Omega).
\end{align*}
Furthermore, there holds $\|u_{ms}-\tilde{u}_{ms,h}\|_{W^{2,2}(\tilde\Omega)}\le C h^{r-1}$ where $C>0$ depends on $\|u_{ms}\|_{C^{r+1}(\tilde\Omega)}$
but is independent of $s$.
Because $\|u_{ms}\|_{C^{r+1}(\tilde\Omega)}$ is uniformly bounded with respect to $s$, it follows from a Sobolev
embedding theorem that $\tilde u_{ms,h}$ is uniformly bounded.  Thus, since $u_{ms,h}$ is convex and uniformly 
bounded, the sequence $\{\tilde u_{ms,h}\}_{s}$ is locally uniformly equicontinuous,
and thus has a pointwise convergent subsequence.  Standard arguments,
along with $u_{ms}\to u_m$ on $\p\tilde\Omega$,
then show that this limit is a solution to the discrete problem \eqref{eqn:AwanouC1}.
\end{proof}

 \begin{remark}[{interior convergence}]
 {Regarding Theorem~\ref{thm:InteriorConv} note that}:
 \begin{enumerate}
 \item  The ideas and techniques given in this section has been applied
 to standard finite difference discretizations of the \MA problem in \cite{Awanou16T}.
 \item While the results and techniques of Theorem \ref{thm:InteriorConv} are interesting,
 it is not immediately clear how to obtain the Dirichlet boundary condition $\tilde u_{m,h} = u_m$
 on $\partial \tilde\Omega$ since $u_m$ is not given data.  One can alternatively 
 use $\tilde u_{m,h}|_{\p\tilde\Omega}=0$, but this condition is not consistent
 with problem \eqref{eqn:SourceRegularized}.  We also point out
 that $h_0$ depends on ${\rm dist}\{\p \tilde\Omega,\p \Omega\}$,
 and therefore Theorem \ref{thm:InteriorConv} suggests
 we cannot take $\tilde\Omega$ to be arbitrarily close to $\Omega$.
 \end{enumerate}
 \eremk\end{remark}

 \pgfplotsset{samples=600}
 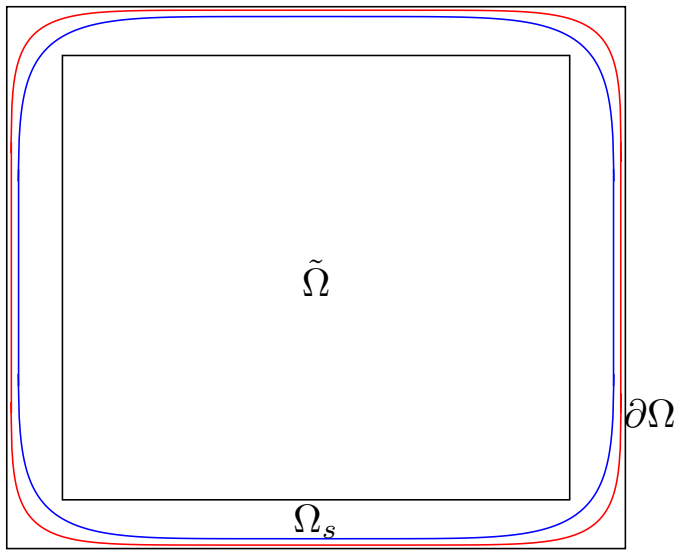
\begin{figure}
 \begin{tikzpicture}[scale=1.5]

  \begin{axis}
  [ axis lines=none,xmin=-1.25,xmax=1.25]			
  \addplot[domain=-1:1.01,blue]{(0.8-\x^6)^(1/6)};
  \addplot[domain=-1:1.01,blue]{-(0.8-\x^6)^(1/6)};
  
    \addplot[domain=-1:1,red]{(0.9-\x^8)^(1/8)};
  \addplot[domain=-1:1,red]{-(0.9-\x^8)^(1/8)};
  
  \addplot[blue] coordinates {(0.8^(1/6)-0.002,-0.4)(0.8^(1/6)-0.002,0.4)}; 
  \addplot[blue] coordinates {(-0.8^(1/6)+0.002,-0.4)(-0.8^(1/6)+0.002,0.4)}; 

  \addplot[red] coordinates {(0.9^(1/8)-0.002,-0.5)(0.9^(1/8)-0.002,0.5)}; 
      \addplot[red] coordinates {(-0.9^(1/8)+0.002,-0.5)(-0.9^(1/8)+0.002,0.5)}; 
 
\draw (-1,-1) rectangle (1,1);
\draw (-.82,-.82) rectangle (.82,.82);
\node at (0,0) {$\tilde\Omega$};
\node  at (0,-0.9) {$\Omega_s$};
\node[right]  at (0.95,-0.5) {$\partial \Omega$};
\end{axis}

 \end{tikzpicture}
 \caption{\label{fig:AwanouConstruction}Pictorial description of the proof of Theorem \ref{thm:InteriorConv}.
 Here, $\tilde\Omega\subset \Omega_s\subset \Omega_{s+1}\subset \Omega$, 
 where $\Omega$ is the physical domain, $\tilde\Omega$ is the computational domain,
 and $\{\Omega_s\}$ are smooth and uniformly convex approximations to $\Omega$.}
 \end{figure}

\definecolor{tableheadcolor}{gray}{0.92}
\definecolor{rulecolor}{RGB}{0,71,171}

\section{Numerical Examples}
\label{sec:Numerics}

\epigraph{\tiny{\it The high point of this classical algorithmic age was perhaps reached in the work of Leonhard Euler [...] Innumerable numerical examples are dispersed in the (so far) seventy volumes of his collected works, showing that Euler always kept foremost in his mind the immediate numerical use of his formulas and algorithms.}}{\tiny{P.~Henrici \cite{Elements}}}

In this section we perform  some simple numerical examples
to show the efficiency and accuracy of 
some of the numerical schemes discussed in the previous sections. 
We consider
three different test problems, each reflecting 
different scenarios of regularity. 
These are computed
using the wide-stencil finite difference scheme \eqref{eq:WSFD}, 
the analogous filtered scheme \eqref{eq:FilterFD}, {Oliker-Prussner method \eqref{OP}}, the
$C^0$ finite element method \eqref{eqn:C0PenaltyMethod}, and its regularized
version using the vanishing moment methodology \eqref{eqn:VMMC0IP}.
We emphasize that
these tests
are not meant to form comparisons, but rather
to highlight their advantages in different situations.

\stepcounter{numexample}
\subsection{Example \arabic{numexample}: smooth solution}
In the first set of experiments, we take 
the data such that the Monge-Amp\`ere equation has a $C^\infty(\Omega)$ solution: $\Omega = (-1,1)^2$, 
\begin{align}
\label{eqn:SmoothSolnExample}
f(x_1,x_2) = (1 + x_1^2 + x_2^2) e^{x_1^2 + x_2^2}, \quad u(x_1,x_2)= e^{\frac{x^2_1 + x_2^2}2}.
\end{align}
In this setting, the Galerkin methods discussed in Sections \ref{sec:C0FEMs}--\ref{sec:MixedFEMs}
are advantageous due to their relative high order.  We implement the the $C^0$ finite element method 
\eqref{eqn:C0PenaltyMethod}  {and the Oliker-Prussner method \eqref{OP}} on a sequence of mesh refinements and report
the resulting errors in Figure \ref{fig:FigureTest1D}.  In agreement with Theorem 
\ref{eqn:AnErrorEstimate} (with $\ell = r-2$ and $\alpha=1$), the plots show optimal order convergence
in $W^{2,p}_h$-norm with respect to the discretization parameter $h$ {for the Galerkin methods}.  In terms of the degrees of freedom (DOFs), 
the errors scale like
\begin{align*}
\|u-u_h\|_{W^{2,p}_h(\Omega)} = \mathcal{O}(DOFs^{(1-r)/2}).
\end{align*}
The errors in $L^\infty$ converge with optimal order 
provided that the polynomial degree is sufficiently high.  Figure \ref{fig:FigureTest1D}
shows that
\begin{alignat*}{2}
\|u-u_h\|_{L^{\infty}(\Omega)} &= \mathcal{O}(DOFs^{(-1-r)/2})\quad &&r=3,4,\\
\|u-u_h\|_{L^\infty(\Omega)} &= \mathcal{O}(DOFs^{-1})\quad &&r=2.
\end{alignat*}
These rates are proven in \cite{Neilan13}. 
For the Oliker-Prussner method and finite difference methods defined on translation invariant meshes, we define its $W^{2,p}$ error on the nodal set as
\[
  \| u - u_h \|_{W^{2,p}_h(\Nh)} = \left( h^d \sum_{x_h \in \Nhi, \be_j \in S } | \Delta_{\be_j} (u - u_h)(x_h) |^p \right)^{1/p}
\]
where $S$ is the $9$-points stencil in two space dimensions and $\Delta_{\be_j} v(x_h)$ denotes the centered second difference, {defined in \eqref{eq:2nddiff}}, of the function $v$ at node $x_h$ in the direction $\be_j$. 
We observe in Figure \ref{fig:FigureTest1D} that, for the Oliker-Prussner method,
\begin{align*}
\|u-u_h\|_{W^{2,p}_h(\Nh)} = \mathcal{O}(DOFs^{-1})
\quad \mbox{and} \quad
\|u-u_h\|_{L^\infty(\dm)} = \mathcal{O}(DOFs^{-1})
\end{align*}
These results on $W^{2,p}$ error are consistent with the theorems proven in \cite{NeilanZhangSINUM}
and Theorem \ref{thm:OPMainEst}.

%
%
\begin{figure}[h] 
\begin{tikzpicture}[scale=0.6]
\begin{loglogaxis}[
	title=Example \arabic{numexample}: $L^\infty$ Error,
	ytick pos=left,
	xtick={100,1000,10000,100000,1000000},
	xticklabels={$10^2$,$10^{3}$,$10^{4}$,$10^{5}$,$10^{6}$},
	legend entries={FE2,FE3,FE4, OP},
	legend pos=south west
]
\addplot table [x=Dim,y=Linferr,col sep=comma]{GalerkinTest1_k2.csv};
\addplot table [x=Dim,y=Linferr,col sep=comma]{GalerkinTest1_k3.csv};
\addplot table [x=Dim,y=Linferr,col sep=comma]{GalerkinTest1_k4.csv};
\addplot table [x=Dim,y=Linferr,col sep=comma]{OPTest1.csv};
\end{loglogaxis}
\end{tikzpicture}
\begin{tikzpicture}[scale=0.6]
\begin{loglogaxis}[
	title=Example \arabic{numexample}: $W^{2,1}$ Error,
	ytick pos=left,
	xtick={100,1000,10000,100000,1000000},
	xticklabels={$10^2$,$10^{3}$,$10^{4}$,$10^{5}$,$10^{6}$},
	legend entries={FE2,FE3,FE4,OP},
	legend pos=south west
]
\addplot table [x=Dim,y=W21err,col sep=comma]{GalerkinTest1_k2.csv};
\addplot table [x=Dim,y=W21err,col sep=comma]{GalerkinTest1_k3.csv};
\addplot table [x=Dim,y=W21err,col sep=comma]{GalerkinTest1_k4.csv};
\addplot table [x=Dim,y=W21err,col sep=comma]{OPTest1.csv};
\end{loglogaxis}
\end{tikzpicture}
\begin{tikzpicture}[scale=0.6]
\begin{loglogaxis}[
	title=Example \arabic{numexample}: ${H^2}$ Error,
	ytick pos=left,
	xtick={100,1000,10000,100000,1000000},
	xticklabels={$10^2$,$10^{3}$,$10^{4}$,$10^{5}$,$10^{6}$},
	legend entries={FE2,FE3,FE4,OP},
	legend pos=south west
]
\addplot table [x=Dim,y=H2err,col sep=comma]{GalerkinTest1_k2.csv};
\addplot table [x=Dim,y=H2err,col sep=comma]{GalerkinTest1_k3.csv};
\addplot table [x=Dim,y=H2err,col sep=comma]{GalerkinTest1_k4.csv};
\addplot table [x=Dim,y=H2err,col sep=comma]{OPTest1.csv};
\end{loglogaxis}
\end{tikzpicture}
\caption{Example \arabic{numexample}: Errors versus degrees of freedom for the $C^0$ finite element method 
 \eqref{eqn:C0PenaltyMethod} with polynomial degrees $r=2,3,4$, {and the Oliker-Prussner method \eqref{OP}} applied to the smooth test problem 
 \eqref{eqn:SmoothSolnExample}.}
\label{fig:FigureTest1D}
\end{figure}
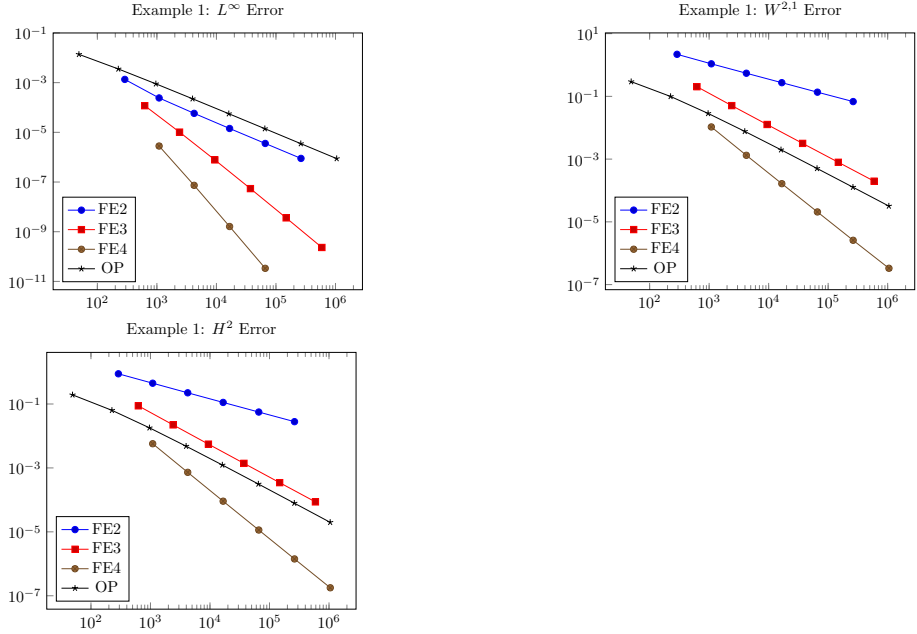


\stepcounter{numexample}
\subsection{Example \arabic{numexample}:  nonclassical solution}
In this set of experiments, we again take $\Omega = (-1,1)^2$,
but choose the data such that the resulting solution 
is not a classical one:
\begin{align*}
 f(x_1,x_2)  =& \begin{dcases}
    16, &|x| \leq 1/2,  \\
    64 -  16 |x|^{-1},  &|x| > 1/2.
      \end{dcases} \\
 u(x_1,x_2)  =& \begin{dcases}
         2 |x|^2,  & |x| \leq 1/2, \\
       2(|x| - 1/2)^2 + 2 |x|^2, & |x| > 1/2.
       \end{dcases}
\end{align*}
One easily finds that {$u\not\in C^{1,1}(\bar\Omega) \setminus C^2(\Omega)$}.
We implement the $C^0$ finite element method
\eqref{eqn:C0PenaltyMethod}, the wide-stencil 
finite difference scheme \eqref{eq:WSFD} with {a stencil size that consists of $33$ grid points}, and the Oliker-Prussner method {\eqref{OP}}.
We also compare the results with the filtered scheme \eqref{eq:FilterFD} 
The errors, depicted in Figure \ref{fig:FigureTest2D}, show that all methods
converge with similar rates, although the finite element scheme and Oliker-Prussner method
have smaller errors with similar DOFs.  While the rate of convergence
in the $L^\infty$ norm is not obvious from the tests, Figure \ref{fig:FigureTest2D}
clearly shows that all three methods converge in the $W^{2,p}$-norms
with rates
\begin{align}\label{eqn:Example2W2pRates}
\|u-u_h\|_{H^2_h(\Omega)}=\mathcal{O}(DOFs^{-1/4}),\quad \|u-u_h\|_{W^{2,1}_h(\Omega)} = \mathcal{O}(DOFs^{-1/2}).
\end{align}
We note that, for the finite element, these rates seem to be the
same rates of interpolation errors.  Indeed, let $\mathcal{T}_h^\Gamma$
denote the set of triangles in $\mathcal{T}_h$ intersect 
the circle $|x|=1/2$.  Likewise, we let $\calF_h^\Gamma$ denote
the set of edges in $\calF_h^I$ that intersect $\Gamma$.  
Finally, we denote by $\calI_h u$ the nodal interpolant of $u$.

Because $u$ is smooth 
on both $\Omega\cap \{x\in \Omega:\ |x|<1/2\}$ 
and $\Omega\cap \{x\in \Omega:\ x>1/2\}$, we have
by standard interpolation estimates,
\begin{multline*}
\|u-\calI_h u\|_{W^{2,p}_h(\Omega)}^p \le C h^{p(r-1)} + \sum_{T\in \mathcal{T}_h^\Gamma}\|D^2 (u-\calI_hu)\|_{L^p(T)}^p \\
+ \sum_{F\in \calF_h^\Gamma} h_F^{1-p} \big\|\jump{\nabla (u-\calI_h u)}\big\|_{L^p(F)}^p\\
\le C\left(h^{p(r-1)} + \sum_{T\in \mathcal{T}_h^\Gamma} \big(\|D^2 (u-\calI_h u)\|_{L^p(T)}^p + h_T^{-p} \|\nab (u-\calI_h u)\|_{L^p(T)}^p\big)\right),
\end{multline*}
where we used a standard trace inequality.
Applying interpolation estimates and H\"older's inequality, noting that $u\in W^{2,\infty}(\Omega)$, yields
\begin{align*}
\|u-\calI_h u\|_{W^{2,p}_h(\Omega)}^p 
&\le C \Big(h^{p(r-1)} + \sum_{T\in \mathcal{T}_h^\Gamma}  \|D^2 u\|_{L^p(T)}^p\Big)\\
&\le C \Big(h^{p(r-1)} + \sum_{T\in \mathcal{T}_h^\Gamma} h^2_T \|D^2 u\|_{L^\infty(T)}^p\Big)\\
&\le C h^{p(r-1)} + C h,
\end{align*}
where we used that the cardinality of $\mathcal{T}_h^\Gamma$ is $\mathcal{O}(h^{-1})$.
We then take the $p$th root of this inequality
to deduce that $\|u-\calI_h u_h\|_{W^{2,p}_h(\Omega)}=\mathcal{O}(h^{1/p}) = \mathcal{O}(DOFs^{-1/(2p)})$,
which is the same rates as \eqref{eqn:Example2W2pRates}.

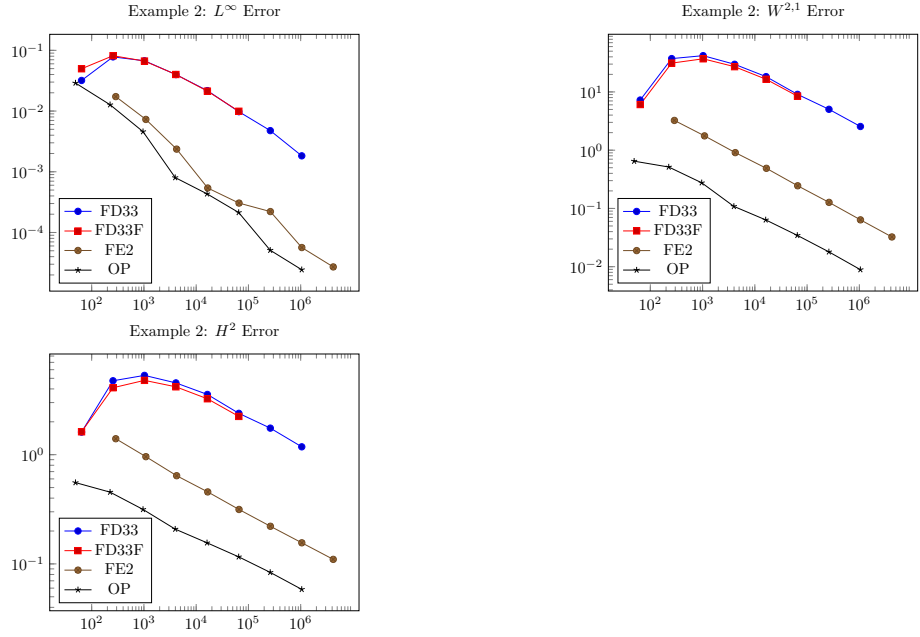
\begin{figure}[h] 
\begin{tikzpicture}[scale=0.6]
\begin{loglogaxis}[
	title=Example \arabic{numexample}: $L^\infty$ Error,
	ytick pos=left,
	xtick={100,1000,10000,100000,1000000},
	xticklabels={$10^2$,$10^{3}$,$10^{4}$,$10^{5}$,$10^{6}$},
	legend entries={FD33,FD33F,FE2, OP},
	legend pos=south west
]
\addplot table [x=Dim,y=Linferr,col sep=comma]{FDTest2_33point.csv};
\addplot table [x=Dim,y=Linferr,col sep=comma]{FDTest2_33Fpoint.csv};
\addplot table [x=Dim,y=Linferr,col sep=comma]{GalerkinTest2_k2.csv};
\addplot table [x=Dim,y=Linferr,col sep=comma]{OPTest2.csv};
\end{loglogaxis}
\end{tikzpicture}
%
%
%
%
%
%
%
\begin{tikzpicture}[scale=0.6]
\begin{loglogaxis}[
	title=Example \arabic{numexample}: $W^{2,1}$ Error,
	ytick pos=left,
	xtick={100,1000,10000,100000,1000000},
	xticklabels={$10^2$,$10^{3}$,$10^{4}$,$10^{5}$,$10^{6}$},
	legend entries={FD33,FD33F,FE2, OP},
	legend pos=south west
]
\addplot table [x=Dim,y=W21err,col sep=comma]{FDTest2_33point.csv};
\addplot table [x=Dim,y=W21err,col sep=comma]{FDTest2_33Fpoint.csv};
\addplot table [x=Dim,y=W21err,col sep=comma]{GalerkinTest2_k2.csv};
\addplot table [x=Dim,y=W21err,col sep=comma]{OPTest2.csv};
\end{loglogaxis}
\end{tikzpicture}
\begin{tikzpicture}[scale=0.6]
\begin{loglogaxis}[
	title=Example \arabic{numexample}: $H^2$ Error,
	ytick pos=left,
	xtick={100,1000,10000,100000,1000000},
	xticklabels={$10^2$,$10^{3}$,$10^{4}$,$10^{5}$,$10^{6}$},
	legend entries={FD33,FD33F,FE2, OP},
	legend pos=south west
]
\addplot table [x=Dim,y=H2err,col sep=comma]{FDTest2_33point.csv};
\addplot table [x=Dim,y=H2err,col sep=comma]{FDTest2_33Fpoint.csv};
\addplot table [x=Dim,y=H2err,col sep=comma]{GalerkinTest2_k2.csv};
\addplot table [x=Dim,y=H2err,col sep=comma]{OPTest2.csv};
\end{loglogaxis}
\end{tikzpicture}
\caption{Example \arabic{numexample}:  Errors versus degrees of freedom for the  $33$-point wide stencil scheme, $33$-point wide stencil filtered scheme, the quadratic $C^0$ finite element method and Oliker-Prussner method.}
\label{fig:FigureTest2D}
\end{figure}


\stepcounter{numexample}
\subsection{Example \arabic{numexample}: Lipschitz and degenerate solution}

In our last set of experiments, we take the domain 
to be $\Omega = (-1,1)^2$ with data
\begin{align*}
 f(x_1,x_2)  &= \begin{dcases}
      36 - 9 x_2^2 x_1^{-6}, & |x_2| \leq |x_1|^3,  \\
      \frac 89 - \frac 59 x_1^2  x_2^{-\frac 23}, &|x_2| > |x_1|^3,
      \end{dcases}\\
 u(x_1,x_2)  =& \begin{dcases}
        |x_1|^4 + \frac {3 x_2^2}{2 x_1^2}, & |x_2| \leq |x_1|^3,  \\
  \frac 12 x_1^2 x_2^{\frac 23} + 2 x_2^{\frac 43}, &|x_2| > |x_1|^3.
       \end{dcases}
\end{align*}
Similar to the previous example, $u$ is not a classical solution 
to \eqref{eqn:MA} {as it} only has regularity $u\in C^{0,1}(\Omega)$
and $u\not\in W^{2,p}(\Omega)$ for any $p > 2$, \cite{Wang95}. 
Moreover, a simple calculation shows that 
$|D^2u(x)| \to \infty$ as $x \to 0$. 
Since the determinant in two dimensions is the product of two eigenvalues of the Hessian and 
$\det D^2u(x) = f(x)$ is bounded in the domain, 
the largest eigenvalue blows up while the other eigenvalue of $D^2 u(x)$ approaches zero  as $x \to 0$.
Hence, the Hessian of the solution degenerates as $x \to 0$.

While the monotone finite difference schemes
presented in Section \ref{sec:FD} are robust
for  problems with low regularity,
Galerkin methods generally fail to capture
solutions
whose second derivatives
are not square integrable; our numerical 
tests show that Newton's method applied to
\eqref{eqn:C0PenaltyMethod} does not
converge for this example even when using very generous initial guesses.
In fact, even for the monotone finite difference schemes
and the Oliker-Prussner method, Newton's method
is very sensitive with respect to  the initial guess 
and the convexity of the iterates
for this problem.  In our implementation, we found that
at each iteration, we require the solution to remain convex.
As Newton's method may not give a convex solution in general,
we applied, if necessary, the algorithm proposed in \cite{Oberman08ConEnv}
to preserve convexity.

In addition to the $33$-point finite difference scheme and Oliker-Prussner method,
we implement the fourth-order regularization
of the $C^0$ finite element method \eqref{eqn:C0PenaltyMethod}
with parameters $\sigma = 100$ and $\epsilon = 0.1 h^2$.
The resulting errors measured in the $L^\infty$ and {$H^1$} 
norms
are {plotted} in Figure \ref{fig:FigureTest3A}.  Similar to the 
previous series of experiments, the plots {show}
that both methods have similar behavior rates.
While the rate in the $L^\infty$ is not clear,
the second plot in Figure \ref{fig:FigureTest3A}
shows that
\[
\|u-u_h\|_{{H^1}(\Omega)} = \mathcal{O}(DOFs^{-1/2}).
\]

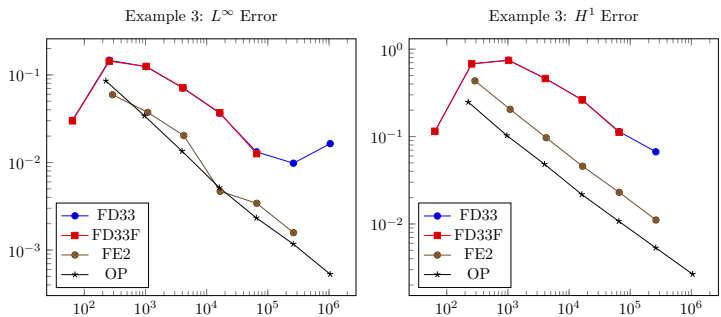
\begin{figure}[h] 
%
\begin{tikzpicture}[scale=0.6]
\begin{loglogaxis}[
	title=Example \arabic{numexample}: $L^\infty$ Error,
	ytick pos=left,
	xtick={100,1000,10000,100000,1000000},
	xticklabels={$10^2$,$10^{3}$,$10^{4}$,$10^{5}$,$10^{6}$},
	legend entries={FD33,FD33F,FE2,OP},
	legend pos=south west
]
\addplot table [x=Dim,y=Linferr,col sep=comma]{FDTest3_33point.csv};
\addplot table [x=Dim,y=Linferr,col sep=comma]{FDTest3_33Fpoint.csv};
\addplot table [x=Dim,y=Linferr,col sep=comma]{GalerkinTest3_k2.csv};
\addplot table [x=Dim,y=Linferr,col sep=comma]{OPTest3.csv};
\end{loglogaxis}
\end{tikzpicture}
\begin{tikzpicture}[scale=0.6]
\begin{loglogaxis}[
	title=Example \arabic{numexample}: ${H^1}$ Error,
	ytick pos=left,
	xtick={100,1000,10000,100000,1000000},
	xticklabels={$10^2$,$10^{3}$,$10^{4}$,$10^{5}$,$10^{6}$},
	legend entries={FD33,FD33F,FE2, OP},
	legend pos=south west
]
\addplot table [x=Dim,y=H1err,col sep=comma]{FDTest3_33point.csv};
\addplot table [x=Dim,y=H1err,col sep=comma]{FDTest3_33Fpoint.csv};
\addplot table [x=Dim,y=H1err,col sep=comma]{GalerkinTest3_k2.csv};
\addplot table [x=Dim,y=H1err,col sep=comma]{OPTest3.csv};
\end{loglogaxis}
\end{tikzpicture}
 \caption{Example \arabic{numexample}: Errors versus
 degrees of freedom for the $33$-point wide-stencil scheme, the $33$-point wide stencil filtered scheme,
 and the quadratic $C^0$ finite element method with regularization and Oliker-Prussner method.}
\label{fig:FigureTest3A}
\end{figure}

\section{Concluding remarks}

\epigraph{\tiny{\it ``And if anyone knows anything about anything'' said Bear to himself, ``it's Owl who knows something about something,'' he said, ``or my name is not Winnie-the-Pooh,'' he said. ``Which it is,'' he added. ``So there you are.''}}{\tiny{Winnie-the-Pooh \cite{TaoPooh}}}

In this work we have reviewed the progress that has been made concerning the approximation and numerical analysis of the \MA problem.  In doing so we highlighted how to develop a convergence analysis of wide stencil finite difference schemes as well as their generalizations, schemes
based on geometric considerations, and finite element methods.
A focus that we have taken, and one of recent development,
is the derivation of rates of convergence for these discretizations.

Despite fundamental advances in only the past decade, there still remain several open problems in the analysis of computational methods for Monge-Amp\`ere problems.  One of these is the derivation of rates of convergence for the Oliker--Prussner scheme on unstructured grids.  Another basic problem is rates of convergence of any of the schemes presented in this work assuming that the solution is not a classical one, i.e., without the assumption $u\in C^{2,\alpha}(\bar\Omega)$. 
In most of the error analyses we have presented, it is assumed that $0 < \lambda I \leq D^2u(x) \leq \Lambda I$ for all $x \in \Omega$. 
However, if the function $f(x)$ is discontinuous, the Hessian of the solution may be degenerate as the third example in the numerics section illustrates.  
The design  and analysis of robust and high order numerical schemes to capture  degenerate solutions remains a challenging problem. 
{\it A posteriori} error estimation, and adaptive methods based on the existing schemes are nonexistent.
Finally let us mention that, as far as we are aware, except for the recent work \cite{Berman18A}, rates of convergence are restricted to the Dirichlet problem \eqref{eqn:MA}; extensions to, e.g., the applications discussed in Section~\ref{sub:Geoapplications} is still unchartered territory.
  
In conclusion, we know something about the numerical analysis of the \MA problem, but there is much more that needs to be developed. It is our hope that this overview will encourage the numerical analysis community to work on the interesting, and challenging, problems found in geometry in general, and those that the \MA equation in particular present to us.

\section*{Acknowledgement}
The work of MJN was supported by NSF Grant DMS--1719829.
The work of AJS was supported by NSF Grant DMS--1720213.
The work of WZ was supported by NSF Grant DMS--1818861.

\bibliographystyle{abbrv}


\begin{thebibliography}{100}

\bibitem{UnifiedDG02}
D.N. Arnold, F.~Brezzi, B.~Cockburn, and L.D. Marini.
\newblock Unified analysis of discontinuous {G}alerkin methods for elliptic
  problems.
\newblock {\em SIAM J. Numer. Anal.}, 39(5):1749--1779, 2001/02.

\bibitem{ArnoldBadAss}
V.I. Arnol'd.
\newblock On the teaching of mathematics.
\newblock {\em Uspekhi Mat. Nauk}, 53(1(319)):229--234, 1998.

\bibitem{Awanou14H}
G.~Awanou.
\newblock Isogeometric method for the elliptic {M}onge-{A}mp\`ere equation.
\newblock In {\em Approximation theory {XIV}: {S}an {A}ntonio 2013}, volume~83
  of {\em Springer Proc. Math. Stat.}, pages 1--13. Springer, Cham, 2014.

\bibitem{Awanou15E}
G.~Awanou.
\newblock Quadratic mixed finite element approximations of the
  {M}onge-{A}mp\`ere equation in 2{D}.
\newblock {\em Calcolo}, 52(4):503--518, 2015.

\bibitem{Awanou15Z}
G.~Awanou.
\newblock Smooth approximations of the {A}leksandrov solution of the
  {M}onge-{A}mp\`ere equation.
\newblock {\em Commun. Math. Sci.}, 13(2):427--441, 2015.

\bibitem{Awanou15F}
G.~Awanou.
\newblock Spline element method for {M}onge-{A}mp\`ere equations.
\newblock {\em BIT}, 55(3):625--646, 2015.

\bibitem{Awanou15c}
G.~Awanou.
\newblock Standard finite elements for the numerical resolution of the elliptic
  {M}onge-{A}mp\`ere equations: classical solutions.
\newblock {\em IMA J. Numer. Anal.}, 35(3):1150--1166, 2015.

\bibitem{Awanou16T}
G.~Awanou.
\newblock On standard finite difference discretizations of the elliptic
  {M}onge-{A}mp\`ere equation.
\newblock {\em J. Sci. Comput.}, 69(2):892--904, 2016.

\bibitem{AwanouOops}
G.~Awanou.
\newblock Erratum to: {Q}uadratic mixed finite element approximations of the
  {M}onge-{A}mp\`ere equation in 2{D} [ {MR}3421667].
\newblock {\em Calcolo}, 54(1):281--297, 2017.

\bibitem{Awanou17ABC}
G.~Awanou.
\newblock Standard finite elements for the numerical resolution of the elliptic
  {M}onge-{A}mp\`ere equation: {A}leksandrov solutions.
\newblock {\em ESAIM Math. Model. Numer. Anal.}, 51(2):707--725, 2017.

\bibitem{AwanouAwi16}
G.~Awanou and R.~Awi.
\newblock Convergence of finite difference schemes to the {A}leksandrov
  solution of the {M}onge-{A}mp\`ere equation.
\newblock {\em Acta Appl. Math.}, 144:87--98, 2016.

\bibitem{AwanouLi14}
G.~Awanou and H.~Li.
\newblock Error analysis of a mixed finite element method for the
  {M}onge-{A}mp\`ere equation.
\newblock {\em Int. J. Numer. Anal. Model.}, 11(4):745--761, 2014.

\bibitem{AwanouLiMalitz}
G.~Awanou, H.~Li, and E.~Malitz.
\newblock A two--grid method for the ${C}^0$ interior penalty discretization of
  the {M}onge-{A}mp\`ere equation.
\newblock preprint, 2018.

\bibitem{MR1305147}
I.J. Bakelman.
\newblock {\em Convex analysis and nonlinear geometric elliptic equations}.
\newblock Springer-Verlag, Berlin, 1994.
\newblock With an obituary for the author by William Rundell, Edited by Steven
  D. Taliaferro.

\bibitem{BarlesSoug91}
G.~Barles and P.~E. Souganidis.
\newblock Convergence of approximation schemes for fully nonlinear second order
  equations.
\newblock {\em Asymptotic Anal.}, 4(3):271--283, 1991.

\bibitem{BCM16}
J.-D. Benamou, F.~Collino, and J.-M. Mirebeau.
\newblock Monotone and consistent discretization of the {M}onge-{A}mp\`ere
  operator.
\newblock {\em Math. Comp.}, 85(302):2743--2775, 2016.

\bibitem{BenamouFroeseObermanTBC}
J.-D. Benamou, B.D. Froese, and A.M. Oberman.
\newblock Numerical solution of the optimal transportation problem using the
  {M}onge-{A}mp\`ere equation.
\newblock {\em J. Comput. Phys.}, 260:107--126, 2014.

\bibitem{Berman18A}
R.J. Berman.
\newblock Convergence rates for discretized {M}onge-{A}mp\`ere equations and
  quantitative stability of optimal transport.
\newblock arXiv:1803.00785 [math.NA], 2018.

\bibitem{Bohmer08}
K.~B{\"o}hmer.
\newblock On finite element methods for fully nonlinear elliptic equations of
  second order.
\newblock {\em SIAM J. Numer. Anal.}, 46(3):1212--1249, 2008.

\bibitem{MR3167449}
A.~Bonito, J.-L. Guermond, and B.~Popov.
\newblock Stability analysis of explicit entropy viscosity methods for
  non-linear scalar conservation equations.
\newblock {\em Math. Comp.}, 83(287):1039--1062, 2014.

\bibitem{MR1100809}
Y.~Brenier.
\newblock Polar factorization and monotone rearrangement of vector-valued
  functions.
\newblock {\em Comm. Pure Appl. Math.}, 44(4):375--417, 1991.

\bibitem{BGNS11}
S.C. Brenner, T.~Gudi, M.~Neilan, and L.-Y. Sung.
\newblock {$C^0$} penalty methods for the fully nonlinear {M}onge-{A}mp\`ere
  equation.
\newblock {\em Math. Comp.}, 80(276):1979--1995, 2011.

\bibitem{BrennerNeilan12}
S.C. Brenner and M.~Neilan.
\newblock Finite element approximations of the three dimensional
  {M}onge-{A}mp\`ere equation.
\newblock {\em ESAIM Math. Model. Numer. Anal.}, 46(5):979--1001, 2012.

\bibitem{BrennerBook}
S.C. Brenner and L.R. Scott.
\newblock {\em The mathematical theory of finite element methods}, volume~15 of
  {\em Texts in Applied Mathematics}.
\newblock Springer, New York, third edition, 2008.

\bibitem{BrennerSung05}
S.C. Brenner and L.-Y. Sung.
\newblock {$C^0$} interior penalty methods for fourth order elliptic boundary
  value problems on polygonal domains.
\newblock {\em J. Sci. Comput.}, 22/23:83--118, 2005.

\bibitem{MR1351007}
L.A. Caffarelli and X.~Cabr{\'e}.
\newblock {\em Fully nonlinear elliptic equations}, volume~43 of {\em American
  Mathematical Society Colloquium Publications}.
\newblock American Mathematical Society, Providence, RI, 1995.

\bibitem{Calabi90}
E.~Calabi.
\newblock Affine differential geometry and holomorphic curves.
\newblock In {\em Complex geometry and analysis ({P}isa, 1988)}, volume 1422 of
  {\em Lecture Notes in Math.}, pages 15--21. Springer, Berlin, 1990.

\bibitem{CiarletBook}
P.G. Ciarlet.
\newblock {\em The finite element method for elliptic problems}, volume~40 of
  {\em Classics in Applied Mathematics}.
\newblock Society for Industrial and Applied Mathematics (SIAM), Philadelphia,
  PA, 2002.
\newblock Reprint of the 1978 original [North-Holland, Amsterdam; MR0520174 (58
  \#25001)].

\bibitem{CockburnShu98}
B.~Cockburn and C.-W. Shu.
\newblock The local discontinuous {G}alerkin method for time-dependent
  convection-diffusion systems.
\newblock {\em SIAM J. Numer. Anal.}, 35(6):2440--2463, 1998.

\bibitem{MR0213764}
R.~Courant, K.~Friedrichs, and H.~Lewy.
\newblock On the partial difference equations of mathematical physics.
\newblock {\em IBM J. Res. Develop.}, 11:215--234, 1967.

\bibitem{CIL}
M.G. Crandall, H.~Ishii, and P.-L. Lions.
\newblock User's guide to viscosity solutions of second order partial
  differential equations.
\newblock {\em Bull. Amer. Math. Soc. (N.S.)}, 27(1):1--67, 1992.

\bibitem{MR1009457}
J.A. Cuesta and C.~Matr\'{a}n.
\newblock Notes on the {W}asserstein metric in {H}ilbert spaces.
\newblock {\em Ann. Probab.}, 17(3):1264--1276, 1989.

\bibitem{DavydovSaeed13}
O.~Davydov and A.~Saeed.
\newblock Numerical solution of fully nonlinear elliptic equations by
  {B}\"ohmer's method.
\newblock {\em J. Comput. Appl. Math.}, 254:43--54, 2013.

\bibitem{MR3324933}
G.~De~Philippis and A.~Figalli.
\newblock Optimal regularity of the convex envelope.
\newblock {\em Trans. Amer. Math. Soc.}, 367(6):4407--4422, 2015.

\bibitem{MR1989280}
E.J. Dean and R.~Glowinski.
\newblock Numerical solution of the two-dimensional elliptic {M}onge-{A}mp\`ere
  equation with {D}irichlet boundary conditions: an augmented {L}agrangian
  approach.
\newblock {\em C. R. Math. Acad. Sci. Paris}, 336(9):779--784, 2003.

\bibitem{MR2111728}
E.J. Dean and R.~Glowinski.
\newblock Numerical solution of the two-dimensional elliptic {M}onge-{A}mp\`ere
  equation with {D}irichlet boundary conditions: a least-squares approach.
\newblock {\em C. R. Math. Acad. Sci. Paris}, 339(12):887--892, 2004.

\bibitem{MR2169156}
E.J. Dean and R.~Glowinski.
\newblock On the numerical solution of a two-dimensional {P}ucci's equation
  with {D}irichlet boundary conditions: a least-squares approach.
\newblock {\em C. R. Math. Acad. Sci. Paris}, 341(6):375--380, 2005.

\bibitem{GlowinskiDean06B}
E.J. Dean and R.~Glowinski.
\newblock An augmented {L}agrangian approach to the numerical solution of the
  {D}irichlet problem for the elliptic {M}onge-{A}mp\`ere equation in two
  dimensions.
\newblock {\em Electron. Trans. Numer. Anal.}, 22:71--96 (electronic), 2006.

\bibitem{GlowinskiDean06A}
E.J. Dean and R.~Glowinski.
\newblock Numerical methods for fully nonlinear elliptic equations of the
  {M}onge-{A}mp\`ere type.
\newblock {\em Comput. Methods Appl. Mech. Engrg.}, 195(13-16):1344--1386,
  2006.

\bibitem{MR3042570}
K.~Debrabant and E.R. Jakobsen.
\newblock Semi-{L}agrangian schemes for linear and fully non-linear diffusion
  equations.
\newblock {\em Math. Comp.}, 82(283):1433--1462, 2013.

\bibitem{EngelEtAl02}
G.~Engel, K.~Garikipati, T.J.R. Hughes, M.G. Larson, L.~Mazzei, and R.L.
  Taylor.
\newblock Continuous/discontinuous finite element approximations of
  fourth-order elliptic problems in structural and continuum mechanics with
  applications to thin beams and plates, and strain gradient elasticity.
\newblock {\em Comput. Methods Appl. Mech. Engrg.}, 191(34):3669--3750, 2002.

\bibitem{FGN13}
X.~Feng, R.~Glowinski, and M.~Neilan.
\newblock Recent developments in numerical methods for fully nonlinear second
  order partial differential equations.
\newblock {\em SIAM Rev.}, 55(2):205--267, 2013.

\bibitem{FengNeilan16}
X.~Feng, L.~Hennings, and M.~Neilan.
\newblock Finite element methods for second order linear elliptic partial
  differential equations in non-divergence form.
\newblock {\em Math. Comp.}, 86(307):2025--2051, 2017.

\bibitem{FengJensen}
X.~Feng and M.~Jensen.
\newblock Convergent semi-{L}agrangian methods for the {M}onge-{A}mp\`ere
  equation on unstructured grids.
\newblock {\em SIAM J. Numer. Anal.}, 55(2):691--712, 2017.

\bibitem{FengNeilan09Z}
X.~Feng and M.~Neilan.
\newblock Vanishing moment method and moment solutions for fully nonlinear
  second order partial differential equations.
\newblock {\em J. Sci. Comput.}, 38(1):74--98, 2009.

\bibitem{FengNeilanJSC11}
X.~Feng and M.~Neilan.
\newblock Analysis of {G}alerkin methods for the fully nonlinear
  {M}onge-{A}mp\`ere equation.
\newblock {\em J. Sci. Comput.}, 47(3):303--327, 2011.

\bibitem{FengNeilan14Radial}
X.~Feng and M.~Neilan.
\newblock Convergence of a fourth-order singular perturbation of the
  {$n$}-dimensional radially symmetric {M}onge-{A}mp\`ere equation.
\newblock {\em Appl. Anal.}, 93(8):1626--1646, 2014.

\bibitem{TooManySelfCitations18}
X.~Feng, M.~Neilan, and S.~Schnake.
\newblock Interior penalty discontinuous {G}alerkin methods for second order
  linear non-divergence form elliptic {PDE}s.
\newblock {\em J. Sci. Comput.}, 74(3):1651--1676, 2018.

\bibitem{MR3617963}
A.~Figalli.
\newblock {\em The {M}onge-{A}mp\`ere equation and its applications}.
\newblock Zurich Lectures in Advanced Mathematics. European Mathematical
  Society (EMS), Z\"{u}rich, 2017.

\bibitem{FroeseTransportBC}
B.D. Froese.
\newblock A numerical method for the elliptic {M}onge-{A}mp\`ere equation with
  transport boundary conditions.
\newblock {\em SIAM J. Sci. Comput.}, 34(3):A1432--A1459, 2012.

\bibitem{Froese18MF}
B.D. Froese.
\newblock Meshfree finite difference approximations for functions of the
  eigenvalues of the {H}essian.
\newblock {\em Numer. Math.}, 138(1):75--99, 2018.

\bibitem{FroeseOberman11}
B.D. Froese and A.M. Oberman.
\newblock Convergent finite difference solvers for viscosity solutions of the
  elliptic {M}onge-{A}mp\`ere equation in dimensions two and higher.
\newblock {\em SIAM J. Numer. Anal.}, 49(4):1692--1714, 2011.

\bibitem{MR2745457}
B.D. Froese and A.M. Oberman.
\newblock Fast finite difference solvers for singular solutions of the elliptic
  {M}onge-{A}mp\`ere equation.
\newblock {\em J. Comput. Phys.}, 230(3):818--834, 2011.

\bibitem{MR3033017}
B.D. Froese and A.M. Oberman.
\newblock Convergent filtered schemes for the {M}onge-{A}mp\`ere partial
  differential equation.
\newblock {\em SIAM J. Numer. Anal.}, 51(1):423--444, 2013.

\bibitem{GT}
D.~Gilbarg and N.S. Trudinger.
\newblock {\em Elliptic partial differential equations of second order}.
\newblock Classics in Mathematics. Springer-Verlag, Berlin, 2001.
\newblock Reprint of the 1998 edition.

\bibitem{GuanSpruck93}
B.~Guan and J.~Spruck.
\newblock Boundary-value problems on {$S^n$} for surfaces of constant {G}auss
  curvature.
\newblock {\em Ann. of Math. (2)}, 138(3):601--624, 1993.

\bibitem{MR3860123}
J.-L. Guermond, M.~Nazarov, B.~Popov, and I.~Tomas.
\newblock Second-order invariant domain preserving approximation of the {E}uler
  equations using convex limiting.
\newblock {\em SIAM J. Sci. Comput.}, 40(5):A3211--A3239, 2018.

\bibitem{MR3204837}
J.-L. Guermond and R.~Pasquetti.
\newblock Entropy viscosity method for high-order approximations of
  conservation laws.
\newblock In {\em Spectral and high order methods for partial differential
  equations}, volume~76 of {\em Lect. Notes Comput. Sci. Eng.}, pages 411--418.
  Springer, Heidelberg, 2011.

\bibitem{MR2787948}
J.-L. Guermond, R.~Pasquetti, and B.~Popov.
\newblock Entropy viscosity method for nonlinear conservation laws.
\newblock {\em J. Comput. Phys.}, 230(11):4248--4267, 2011.

\bibitem{Gutierrez01}
C.E. Guti{\'e}rrez.
\newblock {\em The {M}onge-{A}mp\`ere equation}.
\newblock Progress in Nonlinear Differential Equations and their Applications,
  44. Birkh\"auser Boston, Inc., Boston, MA, 2001.

\bibitem{FroeseGauss}
B.F. Hamfeldt.
\newblock Convergent approximation of non-continuous surfaces of prescribed
  {G}aussian curvature.
\newblock {\em Commun. Pure Appl. Anal.}, 17(2):671--707, 2018.

\bibitem{Elements}
P.~Henrici.
\newblock {\em Elements of numerical analysis}.
\newblock John Wiley \& Sons, Inc., New York-London-Sydney, 1964.

\bibitem{MR1972219}
M.~Hinterm{\"u}ller, K.~Ito, and K.~Kunisch.
\newblock The primal-dual active set strategy as a semismooth {N}ewton method.
\newblock {\em SIAM J. Optim.}, 13(3):865--888 (2003), 2002.

\bibitem{TaoPooh}
B.~Hoff.
\newblock {\em The {T}ao of {P}ooh}.
\newblock Penguin Books, 1982.

\bibitem{HuangHuanHan10}
J.~Huang, X.~Huang, and W.~Han.
\newblock A new {$C^0$} discontinuous {G}alerkin method for {K}irchhoff plates.
\newblock {\em Comput. Methods Appl. Mech. Engrg.}, 199(23-24):1446--1454,
  2010.

\bibitem{JensenNonConvex}
M.~Jensen.
\newblock Numerical solution of the simple {M}onge--{A}mp{\`e}re equation with
  nonconvex {D}irichlet data on nonconvex domains.
\newblock In Dante Kalise, Karl Kunisch, and Zhiping Rao, editors, {\em
  Hamilton-Jacobi-Bellman Equations: Numerical Methods and Applications in
  Optimal Control}, volume~21 of {\em Radon Series on Computational and Applied
  Mathematics}, pages 129--142. De Gryuter, Boston, Berlin, 2018.

\bibitem{JensenSmearsComparison}
M.~Jensen and I.~Smears.
\newblock On the notion of boundary conditions in comparison principles.
\newblock In Dante Kalise, Karl Kunisch, and Zhiping Rao, editors, {\em
  Hamilton-Jacobi-Bellman Equations: Numerical Methods and Applications in
  Optimal Control}, volume~21 of {\em Radon Series on Computational and Applied
  Mathematics}, pages 143--154. De Gryuter, Boston, Berlin, 2018.

\bibitem{SuliBook}
Bo\v{s}ko~S. Jovanovi\'{c} and Endre S\"{u}li.
\newblock {\em Analysis of finite difference schemes}, volume~46 of {\em
  Springer Series in Computational Mathematics}.
\newblock Springer, London, 2014.
\newblock For linear partial differential equations with generalized solutions.

\bibitem{MR2117877}
L.V. Kantorovich.
\newblock On a problem of {M}onge.
\newblock {\em Zap. Nauchn. Sem. S.-Peterburg. Otdel. Mat. Inst. Steklov.
  (POMI)}, 312(Teor. Predst. Din. Sist. Komb. i Algoritm. Metody. 11):15--16,
  2004.

\bibitem{KaweckiLakkisPryer18}
E.~Kawecki, O.~Lakkis, and T.~Pryer.
\newblock A finite element method for the {M}onge-{A}mp\`ere equation with
  transport boundary conditions.
\newblock arXiv:1807.03535, 2018.

\bibitem{MR3416386}
I.~Kossaczk{\'y}, M.~Ehrhardt, and M.~G{\"u}nther.
\newblock On the non-existence of higher order monotone approximation schemes
  for {HJB} equations.
\newblock {\em Appl. Math. Lett.}, 52:53--57, 2016.

\bibitem{MR0126722}
M.A. Krasnosel'ski\u{\i} and Ja.B. Ruticki\u{\i}.
\newblock {\em Convex functions and {O}rlicz spaces}.
\newblock Translated from the first Russian edition by Leo F. Boron. P.
  Noordhoff Ltd., Groningen, 1961.

\bibitem{KrylovBook88}
N.V. Krylov.
\newblock {\em Nonlinear elliptic and parabolic equations of the second order},
  volume~7 of {\em Mathematics and its Applications (Soviet Series)}.
\newblock D. Reidel Publishing Co., Dordrecht, 1987.
\newblock Translated from the Russian by P. L. Buzytsky [P. L. Buzytski{\u\i}].

\bibitem{LakkisPryer11ABC}
O.~Lakkis and T.~Pryer.
\newblock A finite element method for second order nonvariational elliptic
  problems.
\newblock {\em SIAM J. Sci. Comput.}, 33(2):786--801, 2011.

\bibitem{LakkisPryer13}
O.~Lakkis and T.~Pryer.
\newblock A finite element method for nonlinear elliptic problems.
\newblock {\em SIAM J. Sci. Comput.}, 35(4):A2025--A2045, 2013.

\bibitem{LiNochettoZhangMA}
W.~Li and R.H. Nochetto.
\newblock Optimal pointwise error estimates for two-scale methods for the
  {M}onge-{A}mp\`ere equation.
\newblock {\em SIAM J. Numer. Anal.}, 56(3):1915--1941, 2018.

\bibitem{LiNochettoCE}
W.~Li and R.H. Nochetto.
\newblock Two-scale methods for convex envelopes.
\newblock arXiv:1812.11519 [math.NA], 2018.

\bibitem{MR1314597}
P.-L. Lions and P.E. Souganidis.
\newblock Convergence of {MUSCL} and filtered schemes for scalar conservation
  laws and {H}amilton-{J}acobi equations.
\newblock {\em Numer. Math.}, 69(4):441--470, 1995.

\bibitem{Mirebeau15}
J.-M. Mirebeau.
\newblock Discretization of the 3{D} {M}onge-{A}mpere operator, between wide
  stencils and power diagrams.
\newblock {\em ESAIM Math. Model. Numer. Anal.}, 49(5):1511--1523, 2015.

\bibitem{MR3504992}
J.-M. Mirebeau.
\newblock Minimal stencils for discretizations of anisotropic {PDE}s preserving
  causality or the maximum principle.
\newblock {\em SIAM J. Numer. Anal.}, 54(3):1582--1611, 2016.

\bibitem{MotzkinWasow53}
T.S. Motzkin and W.~Wasow.
\newblock On the approximation of linear elliptic differential equations by
  difference equations with positive coefficients.
\newblock {\em J. Math. Physics}, 31:253--259, 1953.

\bibitem{Neilan13}
M.~Neilan.
\newblock Quadratic finite element approximations of the {M}onge-{A}mp\`ere
  equation.
\newblock {\em J. Sci. Comput.}, 54(1):200--226, 2013.

\bibitem{Neilan14ABC}
M.~Neilan.
\newblock Finite element methods for fully nonlinear second order {PDE}s based
  on a discrete {H}essian with applications to the {M}onge-{A}mp\`ere equation.
\newblock {\em J. Comput. Appl. Math.}, 263:351--369, 2014.

\bibitem{Neilan14}
M.~Neilan.
\newblock A unified analysis of three finite element methods for the
  {M}onge-{A}mp\`ere equation.
\newblock {\em Electron. Trans. Numer. Anal.}, 41:262--288, 2014.

\bibitem{NeilanSelfCitation17}
M.~Neilan.
\newblock Convergence analysis of a finite element method for second order
  non-variational elliptic problems.
\newblock {\em J. Numer. Math.}, 25(3):169--184, 2017.

\bibitem{NSWActa}
M.~Neilan, A.J. Salgado, and W.~Zhang.
\newblock Numerical analysis of strongly nonlinear {PDE}s.
\newblock {\em Acta Numer.}, 26:137--303, 2017.

\bibitem{NeilanZhangSINUM}
M.~Neilan and W.~Zhang.
\newblock Rates of convergence in {$W^2_p$}-norm for the {M}onge-{A}mp\`ere
  equation.
\newblock {\em SIAM J. Numer. Anal.}, 56(5):3099--3120, 2018.

\bibitem{NochettoNtogkasFilter}
R.H. Nochetto and D.~Ntogkas.
\newblock Convergent two-scale filtered schemme for the {M}onge-{A}mp\`{e}re
  equation.
\newblock {\em arXiv:1807.04866}, 2018.

\bibitem{NochettoNtogkasZhangfirsttwoscale}
R.H. Nochetto, D.~Ntogkas, and W.~Zhang.
\newblock Two-scale method for the {M}onge-{A}mp\`ere equation: {C}onvergence
  to the viscosity solution.
\newblock {\em Math. Comp.}, 88(316):637--664, 2019.

\bibitem{NochettoNtogkasZhang}
R.H. Nochetto, D.~Ntogkas, and W.~Zhang.
\newblock Two-scale method for the {M}onge-{A}mp\`{e}re equation: pointwise
  error estimates.
\newblock {\em IMA J. Numer. Anal.}, 2019.
\newblock to appear.

\bibitem{NochettoZhang16}
R.H. Nochetto and W.~Zhang.
\newblock Discrete {ABP} estimate and convergence rates for linear elliptic
  equations in non-divergence form.
\newblock {\em Found. Comput. Math.}, 18(3):537--593, 2018.

\bibitem{NochettoZhangMA}
R.H. Nochetto and W.~Zhang.
\newblock Pointwise rates of convergence for the {O}liker-{P}russner method for
  the {M}onge-{A}mp\`ere equation.
\newblock {\em Numer. Math.}, 2019.
\newblock to appear.

\bibitem{NorrisWestcoff76}
A.P. Norris and B.S. Westcott.
\newblock Computation of reflector surfaces for bivariate beamshaping in the
  elliptic case.
\newblock {\em Journal of Physics A: Mathematical and General}, 9(12):2159,
  1976.

\bibitem{Oberman06}
A.M. Oberman.
\newblock Convergent difference schemes for degenerate elliptic and parabolic
  equations: {H}amilton-{J}acobi equations and free boundary problems.
\newblock {\em SIAM J. Numer. Anal.}, 44(2):879--895 (electronic), 2006.

\bibitem{Oberman08ConEnv}
A.M. Oberman.
\newblock Computing the convex envelope using a nonlinear partial differential
  equation.
\newblock {\em Math. Models Methods Appl. Sci.}, 18(5):759--780, 2008.

\bibitem{Oberman08}
A.M. Oberman.
\newblock Wide stencil finite difference schemes for the elliptic
  {M}onge-{A}mp\`ere equation and functions of the eigenvalues of the
  {H}essian.
\newblock {\em Discrete Contin. Dyn. Syst. Ser. B}, 10(1):221--238, 2008.

\bibitem{MR3621814}
A.M. Oberman and Y.~Ruan.
\newblock A partial differential equation for the rank one convex envelope.
\newblock {\em Arch. Ration. Mech. Anal.}, 224(3):955--984, 2017.

\bibitem{Oliker93}
V.~Oliker and E.~Newman.
\newblock The energy conservation equation in the reflector mapping problem.
\newblock {\em Appl. Math. Lett.}, 6(1):91--95, 1993.

\bibitem{OlikerWaltman86}
V.~Oliker and P.~Waltman.
\newblock Radially symmetric solutions of a {M}onge-{A}mp\`ere equation arising
  in a reflector mapping problem.
\newblock In {\em Differential equations and mathematical physics
  ({B}irmingham, {A}la., 1986)}, volume 1285 of {\em Lecture Notes in Math.},
  pages 361--374. Springer, Berlin, 1987.

\bibitem{Oliker84}
V.I. Oliker.
\newblock Hypersurfaces in {${\bf R}^{n+1}$} with prescribed {G}aussian
  curvature and related equations of {M}onge-{A}mp\`ere type.
\newblock {\em Comm. Partial Differential Equations}, 9(8):807--838, 1984.

\bibitem{Oliker87Reflect}
V.I. Oliker.
\newblock Near radially symmetric solutions of an inverse problem in geometric
  optics.
\newblock {\em Inverse Problems}, 3(4):743--756, 1987.

\bibitem{Oliker88}
V.I. Oliker and L.D. Prussner.
\newblock On the numerical solution of the equation {$(\partial^2z/\partial
  x^2)(\partial^2z/\partial y^2)-((\partial^2z/\partial x\partial y))^2=f$} and
  its discretizations. {I}.
\newblock {\em Numer. Math.}, 54(3):271--293, 1988.

\bibitem{MR3289212}
G.~P{\'o}lya.
\newblock {\em How to solve it}.
\newblock Princeton Science Library. Princeton University Press, Princeton, NJ,
  2014.
\newblock A new aspect of mathematical method, With a foreword by John H.
  Conway, Reprint of the second (2004) edition.

\bibitem{MR1035606}
L.~R\"{u}schendorf and S.T. Rachev.
\newblock A characterization of random variables with minimum {$L^2$}-distance.
\newblock {\em J. Multivariate Anal.}, 32(1):48--54, 1990.

\bibitem{MR1062553}
L.~R\"{u}schendorf and S.T. Rachev.
\newblock Corrigendum: ``{A} characterization of random variables with minimum
  {$L^2$}-distance''.
\newblock {\em J. Multivariate Anal.}, 34(1):156, 1990.

\bibitem{RussellBook}
B.~Russell.
\newblock {\em The Scientific Outlook}.
\newblock Routledge, 1931.

\bibitem{Savin13}
O.~Savin.
\newblock Pointwise {$C^{2,\alpha}$} estimates at the boundary for the
  {M}onge-{A}mp\`ere equation.
\newblock {\em J. Amer. Math. Soc.}, 26(1):63--99, 2013.

\bibitem{MR2425007}
E.~Schmutz.
\newblock Rational points on the unit sphere.
\newblock {\em Cent. Eur. J. Math.}, 6(3):482--487, 2008.

\bibitem{TrudWang05}
N.S. Trudinger and X.-J. Wang.
\newblock The affine {P}lateau problem.
\newblock {\em J. Amer. Math. Soc.}, 18(2):253--289, 2005.

\bibitem{TrudWang08}
N.S. Trudinger and X.-J. Wang.
\newblock Boundary regularity for the {M}onge-{A}mp\`ere and affine maximal
  surface equations.
\newblock {\em Ann. of Math. (2)}, 167(3):993--1028, 2008.

\bibitem{Wang95}
X.-J. Wang.
\newblock Some counterexamples to the regularity of monge-amp\'ere equations.
\newblock {\em Proc. Amer. Math. Soc.}, 123(3):841--845, 1995.

\bibitem{XWang96}
X.-J. Wang.
\newblock On the design of a reflector antenna.
\newblock {\em Inverse Problems}, 12(3):351--375, 1996.

\bibitem{ZTZ}
O.C. Zienkiewicz and R.L. Taylor.
\newblock {\em The finite element method. {V}ol. 1}.
\newblock Butterworth-Heinemann, Oxford, fifth edition, 2000.
\newblock The basis.

\end{thebibliography}

\end{document}